\documentclass[11pt]{article}
\usepackage{amssymb,amsfonts}
\usepackage{amsmath,latexsym}
\usepackage{mathrsfs}
\usepackage[dvips]{graphicx}

\newcommand{\R}{{\mathbb R}}
\newcommand{\N}{{\mathbb N}}

\newcommand{\Z}{{\mathbb Z}}
\newcommand{\C}{{\mathbb C}}
\newcommand{\B}{{\mathbb B}}

\newcommand{\bM}{\mathbf{M}}
\newcommand{\bU}{\mathbf{U}}

\renewcommand{\AA}{\mathcal{A}}
\newcommand{\BB}{\mathcal{B}}
\newcommand{\CC}{\mathscr{C}}
\newcommand{\FF}{\mathcal{F}}
\newcommand{\HH}{\mathcal{H}}
\newcommand{\OO}{\mathcal{O}}
\newcommand{\NN}{\mathcal{N}}
\newcommand{\PP}{\mathcal{P}}
\newcommand{\QQ}{\mathcal{Q}}
\newcommand{\RR}{\mathcal{R}}

\newcommand{\co}{\mathrm{conv}}
\newcommand{\dist}{\mathrm{dist}}
\newcommand{\diag}{\mathrm{diag}}
\newcommand{\dd}{\,\mathrm{d}}
\newcommand{\supp}{\mathrm{supp}}
\newcommand{\singsupp}{\mathrm{sing\thinspace supp}}
\newcommand{\pv}{\mathrm{p.v.}}
\newcommand{\fp}{\mathrm{f.p.}}
\newcommand{\reg}{\mathrm{reg}}
\newcommand{\Var}{\mathrm{Var}}
\newcommand{\Tr}{\mathop{\mathrm{Tr}}}
\newcommand{\sign}{\mathop{\mathrm{sign}}}
\renewcommand{\Re}{\mathop{\mathrm{Re}}}
\renewcommand{\Im}{\mathop{\mathrm{Im}}}
\newcommand{\bepr}{\noindent{\bf Proof. }} 
\newcommand{\enpr}{\hfill \rule{.5em}{.5em}}

\newtheorem{thm}{Theorem}[section]

\newtheorem{prop}[thm]{Proposition}
\newtheorem{lemma}[thm]{Lemma}
\newtheorem{cor}[thm]{Corollary}
\newtheorem{rque}[thm]{Remark}

\begin{document}

\title{The numerical measure of a complex matrix}

\author{Thierry Gallay\thanks{Universit\'e Joseph Fourier, Institut
Fourier (UMR CNRS 5582), BP 74, 38402 St-Martin-d'H\`eres, FRANCE.} 
\and Denis Serre\thanks{\'Ecole Normale Sup\'erieure de Lyon, 
UMPA (UMR CNRS 5669), 46, all\'ee d'Italie, 69364 Lyon, FRANCE.}}

\date{\today}

\maketitle

\begin{abstract}
We introduce and carefully study a natural probability measure 
over the numerical range of a complex matrix $A \in \bM_n(\C)$. 
This {\em numerical measure} $\mu_A$ can be defined as the law of 
the random variable $\langle AX,X\rangle \in \C$ when the vector 
$X \in \C^n$ is uniformly distributed on the unit sphere. 
If the matrix $A$ is normal, we show that $\mu_A$ has a 
piecewise polynomial density $f_A$, which can be identified 
with a multivariate $B$-spline. In the general (nonnormal) case, 
we relate the Radon transform of $\mu_A$ to the spectrum of 
a family of Hermitian matrices, and we deduce an explicit 
representation formula for the numerical density which is 
appropriate for theoretical and computational purposes. As 
an application, we show that the density $f_A$ is polynomial
in some regions of the complex plane which can be characterized
geometrically, and we recover some known results about lacunas 
of symmetric hyperbolic systems in $2+1$ dimensions. Finally,
we prove under general assumptions that the numerical measure of 
a matrix $A \in \bM_n(\C)$ concentrates to a Dirac mass as the 
size $n$ goes to infinity.
\end{abstract}

\section{Introduction}\label{s:intro}

If $A \in \bM_n(\C)$ is a complex square matrix of size $n \in \N^*$, 
the {\em numerical range} of $A$ is the compact subset of the complex 
plane defined by
\[
  W(A) \,=\, \Bigl\{\langle Ax,x\rangle \in \C\,\Big|\, x \in \C^n\,,~
  \|x\| = 1\Bigr\}~,
\]
where $\langle x,y\rangle = y^*x$ is the usual scalar product in
$\C^n$ and $\|x\| = \langle x,x\rangle^{1/2}$. It is quite obvious that
$W(A) \supset \sigma(A)$, where $\sigma(A)$ (the {\em spectrum} of
$A$) is the collection of all eigenvalues of $A$, and that $W(A) =
W(U^*AU)$ for any unitary matrix $U \in \bU_n(\C)$. Moreover, a
celebrated result due to Toeplitz \cite{Toe} and Hausdorff \cite{Hau}
asserts that $W(A)$ is always a convex subset of the complex plane. In
particular, $W(A)$ contains the convex hull of $\sigma(A)$, and it is
easy to verify that $W(A) = \co(\sigma(A))$ if the matrix $A$ is
normal, namely $AA^* = A^*A$. The interested reader is referred to 
Chapter~1 of \cite{HJ} for a detailed discussion of the various properties 
of the numerical range, including complete proofs. 

Let $\partial\B^n = \{x \in \C^n\,|\, \|x\| = 1\}$ be the unit sphere 
in $\C^n$, considered as a real manifold of dimension $2n-1$. 
By definition, the numerical range $W(A)$ is the image of the {\em 
numerical map} $\Phi_A : \partial\B^n\to \C$ defined by
\[
  \Phi_A(x) \,=\, \langle Ax,x\rangle~, \qquad x \in \partial\B^n~.
\]
The algebraic and geometric properties of the map $\Phi_A$ have been
extensively studied, see \cite{Kip,YB,BY,Fie,JAG,JS}. In particular,
the set of all critical values of $\Phi_A$, which we denote by
$\Sigma_A \subset \C$, has received a lot of attention, because this
is an interesting object which contains a lot of information on the
matrix $A$. For instance, it is known that $\partial W(A) \subset
\Sigma_A$ and $W(A) = \co(\Sigma_A)$. In addition, there exists a real
algebraic curve $C_A \subset \C \simeq \R^2$ with the property that
$\Sigma_A = C_A \cup C_A'$, where $C_A'$ denotes the set of all line
segments joining pairs of points of $C_A$ at which $C_A$ has the same
tangent line \cite{JS}. Under generic assumptions on $A$, the
bitangent set $C_A'$ is empty, and the critical set $\Sigma_A$ is
therefore the union of a finite number of closed curves, one of which
is the boundary of the numerical range $W(A)$. This distinguished
curve is smooth, and encloses all the other ones in its interior.  We
refer to Section~\ref{s:geom} below for more details on the geometry
of the singular set, and to Section~\ref{s:examples} for a few
concrete examples.

Our purpose in this paper is to introduce another mathematical quantity 
which is naturally related to the numerical map $\Phi_A$. Given 
$A \in \bM_n(\C)$, the {\em numerical measure} of $A$ is the 
probability measure $\mu_A$ on $\C$ defined by the formula
\begin{equation}\label{def:muA}
  \int_\C \phi(z)\dd \mu_A(z) \,=\, \int_{\partial\B^n}\phi(\langle Ax,x\rangle)
  \dd \bar\sigma(x)~,
\end{equation}
for all continuous functions $\phi : \C \to \C$. Here $\bar\sigma$
denotes the Euclidian measure on the unit sphere $\partial\B^n$,
normalized as a probability measure. In words, the numerical measure
is thus the image under the numerical map of the normalized Euclidean
measure on the unit sphere. Equivalently, if $X$ is a random variable
that is uniformly distributed on $\partial\B^n$, the numerical measure
$\mu_A$ is just the distribution of the random variable $\langle
AX,X\rangle \in \C$. This probabilistic intepretation will be useful
later, especially in Section~\ref{s:stat}.

Our first goal is to establish a few general properties of the
numerical measure $\mu_A$. It is clear by construction that $\mu_A$ is
invariant under unitary conjugations of $A$, namely $\mu_{\,U^*\!AU} =
\mu_A$ for all $U\in \bU_n(\C)$. This is precisely the reason why we
used the Euclidean measure on $\partial\B^n$ in the definition
\eqref{def:muA}. It is also easy to verify that the support of $\mu_A$
is exactly the numerical range $W(A)$, see Section~\ref{s:general}
below. Less obvious, perhaps, is the fact that $\mu_A$ is absolutely
continuous with respect to the Lebesgue measure $\lambda$ on $W(A)$, so
that we can define the {\em numerical density} $f_A$ as the
Radon-Nikodym derivative of $\mu_A$ with respect to $\lambda$ (in the
particular situation where $W(A)$ reduces to a line segment $\Gamma$,
we understand $\lambda$ as the one-dimensional Lebesgue measure on
$\Gamma$, see Section~\ref{s:general}.) We also prove that the
numerical density $f_A$ is strictly positive in the interior of $A$, a
property that can be interpreted as a strong version of Hausdorff's
theorem \cite{Hau}. Finally, we shall see that the singular
support of $\mu_A$ is contained in the critical set $\Sigma_A$, which
means that the numerical density $f_A$ is smooth outside $\Sigma_A$.
In fact, we conjecture that $\singsupp(\mu_A) = \Sigma_A$ for all $A
\in \bM_n(\C)$, but this has not been proved yet.

After these general properties have been established, our next
goal is to give a more precise description of the numerical 
density $f_A$. For this purpose, it is convenient to distinguish 
between various cases:

\smallskip\noindent{\bf 1.} ({\em The scalar case}) If $W(A)$ is reduced 
to a single point $\{z\}$, then $A = z I_n$ (where $I_n \in 
\bM_n(\C)$ denotes the identity matrix) and $\mu_A = \delta_z$. 
In this trivial situation, there is of course no need to introduce
a numerical density. 

\smallskip\noindent{\bf 2.} ({\em The Hermitian case}) Assume that $n
\ge 2$ and that $W(A) \subset \C$ is a line segment. Then we can find
$z \in \C$, $\theta \in [0,2\pi]$, and a Hermitian matrix $H$ such
that $A = zI_n + e^{i\theta}H$. If $\lambda_1 \le \lambda_2 \le \dots
\le \lambda_n$ are the eigenvalues of $H$, we shall see in
Section~\ref{s:normal} that the numerical measure $\mu_H$ is
absolutely continuous with respect to Lebesgue's measure on $\R$, and
that the corresponding density $f_H$ is exactly the normalized {\em
$B$-spline} of degree $n-2$ with knots $\lambda_1, \dots,\lambda_n$
\cite{dB}.  In particular, $f_H$ is polynomial of degree $n-2$ on each
interval $[\lambda_i,\lambda_{i+1}]$, vanishes identically outside
$[\lambda_1,\lambda_n]$, and is continuous at each point $\lambda_i$
together with its derivatives up to order $d_i = n-2-m_i$, where $m_i
\ge 1$ is the multiplicity of $\lambda_i$ as an eigenvalue of $H$ (if
$d_i < 0$, then $f_H$ is discontinuous at $\lambda_i$.) This gives an
explicit representation of the numerical measure $\mu_H$, and the
measure $\mu_A$ is the image of $\mu_H$ under the affine isometry $w
\mapsto z + e^{i\theta}w$.

\smallskip\noindent{\bf 3.} ({\em The normal case}) Suppose now that
$n \ge 3$ and that $A \in \bM_n(\C)$ is a normal matrix whose spectrum
$\sigma(A)$ is not contained in a line segment. Then $W(A) =
\co(\sigma(A))$ is a convex polygon with nonempty interior, and it
turns out that the numerical density $f_A$ is the {\em bivariate
  B-spline} of degree $n-3$ whose knots are the eigenvalues of
$A$. Here we refer to the work of W.~Dahmen \cite{Dah} for the
definition and the main properties of multivariate $B$-splines. In
this particular case, the critical set $\Sigma_A$ is thus the
collection of all line segments joining pairs of eigenvalues of $A$,
and the density $f_A$ is polynomial of degree $n-3$ in each connected
component of $\C \setminus \Sigma_A$. In the generic situation where
no straight line contains more than two eigenvalues of $A$, one can
show that $f_A$ is continuous together with its derivatives up to
order $n-4$ (and is discontinuous on $\partial W(A)$ if $n = 3$.)

\smallskip\noindent{\bf 4.} ({\sl The nonnormal case}) Finally, we
consider the most interesting situation where the matrix $A \in
\bM_n(\C)$ is not normal. In that case, there is no explicit formula
for the numerical density, but the problem can be reduced in some
sense to the Hermitian case by the following simple observation. For
any $\theta \in S^1 = \R/(2\pi\Z)$, let $H(\theta)$ be the Hermitian 
matrix defined by
\begin{equation}\label{def:Htheta}
 H(\theta) \,=\, \frac{1}{2}(e^{-i\theta}A + e^{i\theta}A^*) 
 \,=\, A_1 \cos(\theta) + A_2 \sin(\theta)~,
\end{equation}
where $A_1 = (A+A^*)/2$ and $A_2 = (A-A^*)/(2i)$. Then
$\Re\,(e^{-i\theta}\langle Ax,x\rangle) = \langle H(\theta)x,x\rangle$
for all $x \in \partial\B^n$. Now, if the random variable $X$ is
uniformly distributed on $\partial\B^n$, the distribution of $\langle
H(\theta)X,X\rangle$ is by definition the numerical measure
$\mu_{H(\theta)}$, whereas the distribution of $\Re\,(e^{-i\theta}\langle 
AX,X\rangle)$ is easily identified as the two-dimensional {\em Radon 
transform} of the numerical measure $\mu_A$, evaluated at $\theta \in
[0,2\pi]$. We thus have
\begin{equation}\label{eq:Radon}
  \RR \mu_A (\theta) \,=\,  \mu_{H(\theta)}~, \quad \theta \in S^1~,
\end{equation}
where $\RR$ denotes the two-dimensional Radon transformation. Since
the numerical density of $H(\theta)$ is known to be the $B$-spline
based on the eigenvalues $\lambda_1(\theta),\dots,\lambda_n(\theta)$
of $H(\theta)$, we can reconstruct the numerical measure $\mu_A$ by
inverting the Radon transformation in \eqref{eq:Radon}, using the the
well-known {\em backprojection} method which plays an important role
in tomography \cite{Hel}. This provides a useful representation
formula for the numerical density, as well as an efficient algorithm
for numerical calculations, see Section~\ref{s:Radon} for more
details.

\smallskip 
It is worth mentioning here that the critical set $\Sigma_A$ can be  
conveniently characterized using the family of Hermitian matrices 
$H(\theta)$ associated with $A$. Indeed, if we define the
eigenvalues $\lambda_j(\theta)$ in such a way that they depend
analytically on $\theta$, one can shown that the algebraic curve
$C_A$ which generates $\Sigma_A$ is given by
\[
  C_A \,=\, \Bigl\{e^{i\theta}(\lambda_j(\theta) + i\lambda_j'(\theta)) 
 \,\Big|\, j \in \{1,\dots,n\}\,,~ \theta \in [0,\pi]\Bigr\}~,
\]
see \cite{YB,JAG,JS} and Section~\ref{s:geom} below. In the 
generic case where $\lambda_1(\theta) < \lambda_2(\theta) < \dots
< \lambda_n(\theta)$ for all $\theta \in [0,\pi]$, the bitangent
set $C_A'$ is empty and $\Sigma_A = C_A$. 

\begin{figure}
\begin{tabular}{cc}
\includegraphics[width=8.0cm,height=5cm]{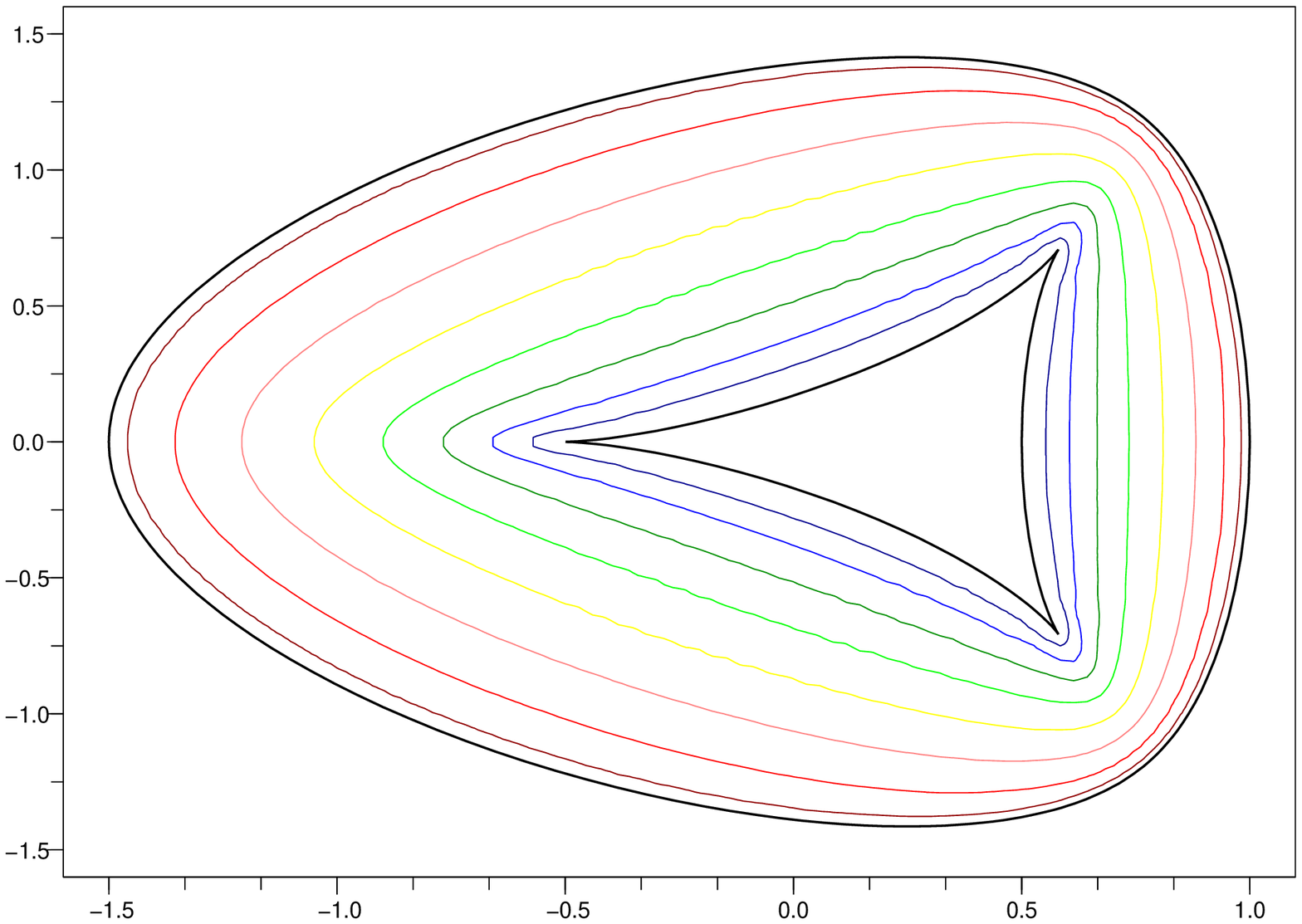} & 
\includegraphics[width=8.0cm,height=5cm]{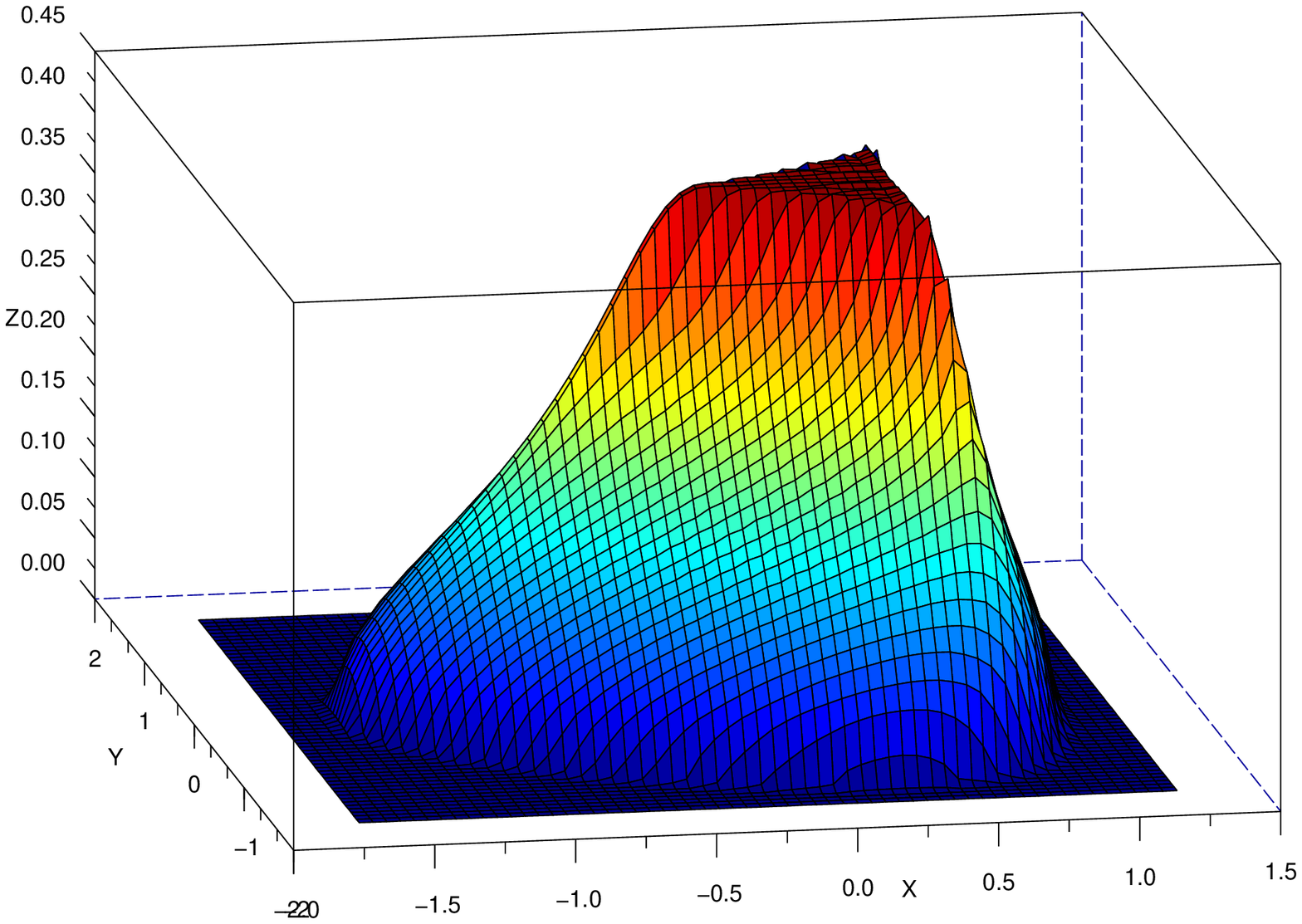} \\
\end{tabular}
\caption{\small The numerical density $f_A$ is represented for a typical 
matrix $A \in \bM_3(\R)$, given by \eqref{def:A3gen} below. In the 
contour plot (left), the exterior ovate curve is the boundary of 
the numerical range $W(A)$, and the other component of the of the 
critical set $\Sigma_A$ is the interior cuspidal triangle. The 
three-dimensional plot (right) confirms that the numerical density 
is continuous, positive inside the numerical range, and constant over 
the cuspidal triangle, in agreement with the results of Section~\ref{s:main}.}
\label{Fig1}
\end{figure}

As was already mentioned, the numerical density $f_A$ is smooth (in
fact, real-analytic) on each connected component of $\C \setminus
\Sigma_A$. The regularity across $\Sigma_A$ is more difficult to
study, but we shall show in Section~\ref{ss:regular} that $f_A$ is
everywhere of class $C^{n-3}$ if $n \ge 3$ and $A \in \bM_n(\C)$
satisfies some generic hypotheses, which exclude in particular the
case of normal matrices.  In addition, for an arbitrary matrix of size
$n$, we shall prove that all derivatives of $f_A$ of order $n-2$
vanish identically in some distinguished regions, which have the
following geometric characterization. For any $z \in \C \setminus
\Sigma_A$, let $N(z)$ be the number of straight lines containing $z$
which are tangent to the curve $C_A$, see \eqref{def:N} below for a
precise definition where possible multiplicities are taken into
account. It is easy to verify that $N(z)$ is constant in each
connected component of $\C \setminus \Sigma_A$, and that $N(z) \le n$
\cite{YB}. The distinguished regions where $f_A$ is polynomial of
degree $n-3$ (if $n \ge 3$) or $f_A \equiv 0$ (if $n = 2$) are exactly
those connected components of $\C \setminus \Sigma_A$ on which $N(z)$
{\em takes its maximal value} $n$. This remarkable property of the
numerical density, which is one of our main results, will be
established in Section~\ref{ss:poly}. The geometric condition $N(z) =
n$ is always satisfied in the complement of the numerical range, where
$f_A$ vanishes identically, but for many matrices of size $n \ge 3$ is
it also met in some regions inside $W(A)$. For instance, in the
three-dimensional example represented in Fig.\,\ref{Fig1}, it is easy
to verify that $N(z) = 3$ if $z$ is outside $W(A)$ or inside the
cuspidal triangle, and $N(z) = 1$ in the intermediate region where the
numerical density is not constant.

At this point, it is necessary to make a connection with the 
theory of {\em lacunas} of symmetric hyperbolic systems of partial 
differential equations in $2+1$ variables. Given $A \in M_n(\C)$, we 
consider the following system of linear PDE's in $\R_t \times \R^2_x$:
\begin{equation}\label{eq:hypsys}
  \partial_t u + A_1 \partial_{x_1}u + A_2 \partial_{x_2}u  \,=\, 0~,
\end{equation}
where $A_1,A_2$ are as in \eqref{def:Htheta} and $u = (u_1,\dots,u_n)^\top
: \R\times\R^2 \to \R^n$. The {\em fundamental solution} of 
\eqref{eq:hypsys} is the unique (matrix-valued) distribution $E$ 
supported in the half-space $\R_+ \times \R^2$ which satisfies 
\begin{equation}\label{eq:hypsys2}
  \partial_t E + A_1 \partial_{x_1}E + A_2 \partial_{x_2}E  \,=\,
  I_n\,\delta_{t=0} \otimes \delta_{x = 0}~.
\end{equation}
One can show that $E(t,x)$ is homogeneous of degree $-2$ in $t$ and $x$, 
it is thus sufficient to consider the time-one trace $E_* = E(1,\cdot)$, 
which is a distribution on $\R^2$. Due to the finite speed of propagation, 
it is well-known that $E_*$ is zero outside a compact set of $\R^2$, but 
it may also happen that $E_*$ vanishes identically in some regions 
inside the domain of influence of the origin. Such regions are called 
{\em lacunas} of the hyperbolic system \eqref{eq:hypsys}. 

The properties of the fundamental solution of symmetric hyperbolic
systems have been studied by many authors, see
e.g. \cite{Pet,YB,BY,ABG1,ABG2}. In the particular case of system
\eqref{eq:hypsys}, J. Bazer and D. Yen have shown that, if one
identifies $\C$ with $\R^2$, the singular support of the distribution
$E_*$ is contained in the critical set $\Sigma_A$, and the (stable)
lacunas of system \eqref{eq:hypsys} are exactly the regions described
above where the numerical density $f_A$ is polynomial of degree
$n-3$. This remarkable coincidence is of course not fortuitous. In
Section~\ref{ss:fundam}, we explain it by showing that the fundamental
solution $E$ can be expressed as a linear combination of derivatives
of order $n-1$ of a homogeneous extension of the numerical density
$f_A$. This connection allows us to recover some of the main results
of \cite{BY}, and therefore confirms that the numerical measure is a
natural quantity attached to the matrix $A$. One might even argue that
$\mu_A$ contains more information than $E_*$, since for instance
$\supp(\mu_A) = W(A)$ while $\supp(E_*)$ is in general strictly
smaller and not necessarily convex, see
Section~\ref{s:examples}. Similarly, we believe that
$\singsupp(\mu_A)$ always coincide with $\Sigma_A$, while
$\singsupp(E_*)$ is usually smaller.

A final question that is worth investigating is the behavior of the
numerical measure $\mu_A$ when the size of the matrix $A$ goes to
infnity. Here of course, specific assumptions have to be made in order
to obtain convergence results. Suppose for instance that $\{A_n\}_{n
  \ge 1}$ is a sequence of complex matrices with $A_n \in \bM_n(\C)$,
$\Tr(A_n) = 0$, and $\|A_n\| \le C$ for all $n \ge 1$. If, for each $n
\ge 1$, $X_n$ is a random variable that is uniformly distributed on
$\partial\B^n$, we show in Section~\ref{s:stat} that the complex
variable $\langle A_n X_n,X_n\rangle$ converges almost surely to zero
as $n \to \infty$. This is reminiscent of the strong law of large
numbers in probability theory. Under slightly stronger assumptions, we
also establish the analog of the central limit theorem in this
context. Our convergence results mean that $\mu_{A}$ is very close to
$\delta_{z}$ when $\dim(A)$ is large, where $z$ is the barycenter of
$\sigma(A)$. This explains why plotting $\langle Ax,x\rangle$ for
randomly chosen points $x \in
\partial\B^n$ is a very unefficient algorithm for determining
the numerical range $W(A)$ if $A$ is a large matrix! 

The rest of the paper is organized as follows. In
Section~\ref{s:general}, we establish some general properties of the
numerical measure.  Section~\ref{s:normal} is devoted to the
particular situations where the matrix $A$ is Hermitian or normal. The
nonnormal case is treated in Sections~\ref{s:Radon}--\ref{s:main},
which constitute the core of the paper. In Section~\ref{s:Radon}, we
derive a representation formula for the numerical density using the
inversion of the Radon transformation. Section~\ref{s:geom} collects a
few results on the geometry of the critical set $\Sigma_A$, which are
mainly borrowed from \cite{Kip,YB,JS}. These informations are used in
Section~\ref{s:main} to derive an explicit formula for the derivatives
of order $n-2$ of the numerical density, which allows us to obtain
generic regularity results and to express the fundamental solution of
the hyperbolic system \eqref{eq:hypsys} in terms of derivatives of the
numerical density. To illustrate our results, a few explicit examples
are treated in Section~\ref{s:examples}. Finally, we investigate 
in Section~\ref{s:stat} the concentration properties of the numerical 
density for large matrices, and we discuss in Section~\ref{s:persp} a 
possible extension of our results to hyperbolic polynomials with an 
arbitrary number of variables.

\bigskip\noindent
{\bf Acknowledgements.} This work has benefited of stimulating
discussions with several of our colleagues, including Y. Colin de 
Verdi\`ere, F. Faure, and A.~Joye. 

\section{General properties of the numerical measure}
\label{s:general}

In this section, we establish a few general properties of the
numerical measure of a complex matrix. In particular, we show
that $\mu_A$ is absolutely continuous with respect to the 
Lebesgue measure on $W(A)$, and we prove a direct sum formula
which will be useful later. 

\subsection{Support and regularity properties}

We first show that the support of the numerical measure always coincides 
with the numerical range of the matrix.

\begin{lemma}\label{lem:supp}
For any $A \in \bM_n(\C)$, one has $\supp(\mu_A) = W(A)$. 
\end{lemma}

\bepr If $V = \C\setminus W(A)$, then $\Phi_A^{-1}(V) = \emptyset$,
hence $\mu_A(V) = \bar \sigma(\Phi_A^{-1}(V)) = 0$. This shows that
$\supp(\mu_A) \subset W(A)$. Conversely, if $V \subset \C$ is any open
set such that $V \cap W(A) \neq \emptyset$, then $\Phi_A^{-1}(V)$ is a
nonempty open subset of $\partial\B^n$, hence $\mu_A(V) = \bar
\sigma(\Phi_A^{-1}(V)) > 0$. Thus $W(A) \subset \supp(\mu_A)$. 
\enpr

\medskip
Our next goal is to locate the singular support of $\mu_A$. We recall 
that $x \in \partial\B^n$ is a {\em regular point} of $\Phi_A$ if the 
differential map $\mathrm{d}_x \Phi_A : T^*_x \partial\B^n \to \C$ is 
onto. Otherwise, we say that $x$ is a {\em critical point}. 
The following characterization will be useful:

\begin{lemma}\label{lem:singset} {\bf \cite{JAG,JS}} 
Let $A \in \bM_n(\C)$.\\
1) A point $x \in \partial\B^n$ is a critical point of the numerical 
map $\Phi_A$ if and only if $x$ is an eigenvector of the Hermitian 
matrix $H(\theta)$ defined in  \eqref{def:Htheta} for some $\theta
\in [0,\pi]$. \\[1mm]
2) The differential of $\Phi_A$ vanishes at $x \in \partial\B^n$ if 
and only if $x$ is an eigenvector of both $A$ and $A^*$. 
\end{lemma}

In other words, the range of the differential $\mathrm{d}_x \Phi_A$
has (real) dimension $1$ if and only if $x$ is an eigenvector of
$H(\theta)$ for a {\em unique} $\theta \in [0,\pi)$, and is reduced to
$\{0\}$ if and only if $x$ is an eigenvector of $H(\theta)$ for 
{\em all} $\theta \in [0,\pi]$. The proof is neither new nor difficult, 
but we shall repeat it here in order to introduce some notation that 
will be needed later on.

\medskip\bepr
Since $\Phi_A(e^{i\theta}x) = \Phi_A(x)$ for all $\theta \in [0,2\pi]$, 
we can consider the numerical map as acting on the quotient space
$\partial\B^n/S^1 \simeq \C P^{n-1}$ \cite{JAG}. Thus, to detect the 
critical points of $\Phi_A$, we study the reduced map $\tilde \Phi_A : 
\partial\B^n/S^1 \to \C$ defined by
\[
  \tilde \Phi_A([x]) \,=\, \Phi_A(x) \,=\, \langle A_1x,x\rangle + 
  i\langle A_2x,x\rangle~, \qquad x \in \partial\B^n~,
\]
where $[x] = \{e^{i\theta}x\,|\,\theta \in S^1\}$ and $A_1,A_2$ are the 
Hermitian matrices introduced in \eqref{def:Htheta}. 

If $x \in \partial\B^n$, the tangent space to $\partial\B^n/S^1$ at $[x]$
is just the $(2n-2)$-dimensional affine subspace $\{[x+y]\,|\, y \in 
\C^n\,,~\langle x,y\rangle = 0\}$. Thus, using the definition above 
of $\tilde \Phi_A$, it is straightforward to verify that, for all 
$y \in \C^n$ with $\langle x,y\rangle = 0$, one has
\begin{equation}\label{eq:DPhiA}
  \frac12\,\mathrm{d}_{[x]}\tilde\Phi_A(y) \,=\,  (A_1x | y) + 
  i(A_2x | y) \,=\, (v_1(x) | y) +i (v_2(x) | y)~, 
\end{equation}
where $(x|y) = \Re\,\langle x,y\rangle = (\Re x)^t (\Re y) + (\Im x)^t 
(\Im y)$ denotes the real scalar product in $\C^n \simeq \R^{2n}$, and 
\[
  v_1(x) \,=\, A_1x - \langle A_1x,x\rangle x~, \qquad
  v_2(x) \,=\, A_2x - \langle A_2x,x\rangle x~. 
\]
Of course, replacing $A_1x$, $A_2x$ with $v_1(x)$, $v_2(x)$ has no effect in
\eqref{eq:DPhiA} since $\langle x,y\rangle = 0$, but after this 
substitution we can let $y$ run over the whole of $\C^n$ without 
increasing the range. So our task is reduced to computing the rank
of the $\R$-linear map $y \mapsto (v_1(x) | y) +i (v_2(x) | y)$, which is
just the rank of the $2\times2$ matrix
\begin{equation}\label{def:D}
  D(x) \,=\, \begin{pmatrix} 
  (v_1(x)|v_1(x)) & (v_1(x)|v_2(x)) \\
  (v_2(x)|v_1(x)) & (v_2(x)|v_2(x)) \end{pmatrix}~.
\end{equation}
By the Cauchy-Schwarz inequality, the positive matrix $D(x)$ is singular 
if and only if there exists $\theta \in [0,\pi]$ such that $v_1(x)
\cos\theta + v_2(x)\sin\theta = 0$, which exactly means that $x$ is an 
eigenvector of $H(\theta)$. Moreover, $D(x) = 0$ if and only if 
$v_1(x) = v_2(x) = 0$, which is equivalent to saying that $x$ is an 
eigenvector of both $A_1$ and $A_2$, hence of both $A$ and $A^*$. \enpr

\medskip
Let $\Sigma_A \subset \C$ denote the set of all {\em critical values} 
of $\Phi_A$, namely $\Sigma_A = \Phi_A(\Gamma(A))$ where $\Gamma(A) 
\subset \partial\B^n$ is the set of all critical points of $\Phi_A$.
Our next result is:
 
\begin{lemma}\label{lem:singsupp} If $A \in \bM_n(\C)$, 
then $\singsupp(\mu_A) \subset \Sigma_A$. 
\end{lemma}

\bepr
If the numerical range $W(A)$ is reduced to a line segment or to 
a single point, then $\Sigma_A = W(A)$, hence $\singsupp(\mu_A) \subset 
\supp(A) = \Sigma_A$ by Lemma~\ref{lem:supp}. Thus, we assume from
now on that $W(A)$ has nonempty interior. By Sard's lemma, 
the critical set $\Sigma_A$ is then a compact subset of $W(A)$ with 
zero Lebesgue measure. We have to show that there exists a smooth
density function $f_A \ge 0$ such that $\mathrm{d}\mu_A(z) = 
f_A(z)\dd z$ on $\C\setminus \Sigma_A$. Clearly, we must have 
$f_A = 0$ on $\C\setminus W(A)$. 

If $z \in W(A) \setminus \Sigma_A$, then $\NN_z := \Phi_A^{-1}(z)$ 
is a compact submanifold of $\partial\B^n$ of codimension $2$, which
depends smoothly on $z$. Using classical arguments, involving 
a partition of unity and the Implicit Function Theorem, it is 
not difficult to verify that, for any continuous function $\phi$ 
with $\supp(\phi) \subset W(A) \setminus \Sigma_A$, one has 
\[
  \int_{\partial\B^n}\phi(\langle Ax,x\rangle)\dd \bar\sigma(x) \,=\,
  \frac{1}{\omega_n} \int_{\C}\phi(z)\left\{\int_{\NN_z}
  \frac{\dd \nu(x)}{2\Delta(x)^{1/2}}\right\}\dd z~,
\]
where $\omega_n = 2\pi^n/((n-1)!)$ is the total measure of 
$\partial\B^n$, $\nu$ is the $(2n-3)$-dimensional Euclidean 
measure on the submanifold $\NN_z$, and $\Delta(x) = \det 
D(x)$ where $D(x)$ is the $2\times2$ matrix defined in \eqref{def:D}. 
Remark that $2\Delta(x)^{1/2} = \bar \lambda_1(x)\bar\lambda_2(x)$, 
where $\bar \lambda_1(x), \bar\lambda_2(x)$ are the singular values
of the differential map $\mathrm{d}_x \Phi_A$. In view of \eqref{def:muA}, 
we conclude that $\mathrm{d}\mu_A(z) = f_A(z)\dd z$ on $\C\setminus 
\Sigma_A$, where
\begin{equation}\label{def:fA}
  f_A(z) \,=\, \frac{1}{\omega_n}\int_{\NN_z} \frac{\dd \nu(x)}{2\Delta(x)^{1/2}}~, 
  \qquad z \in \C\setminus \Sigma_A~.
\end{equation}
It is easily verified that the density $f_A$ is smooth and strictly
positive on $W(A) \setminus \Sigma_A$. 
\enpr

\medskip
The results obtained so far are summarized in the following 
proposition, which also asserts that the numerical measure is 
absolutely continuous with respect to Lebesgue's measure on
$W(A)$. 

\begin{prop}\label{prop:singsupp}
Let $A \in \bM_n(\C)$. \\
1) If the numerical range $W(A)$ has nonempty interior, the numerical
measure $\mu_A$ is absolutely continuous with respect to the 
(two-dimensional) Lebesgue measure on $\C \simeq \R^2$. The numerical
density $f_A = \mathrm{d}\mu_A/\mathrm{d}z$ is smooth outside 
the critical set $\Sigma_A$. \\[1 mm]
2) If $A$ is a nonscalar Hermitian matrix, then $W(A) \subset \R$ 
and the numerical measure is absolutely continuous with respect to the  
(one-dimensional) Lebesgue measure on $\R$. The numerical
density $f_A = \mathrm{d}\mu_A/\mathrm{d}x$ is smooth outside 
the spectrum $\sigma(A)$. 
\end{prop}

\begin{rque}
As is explained in the introduction, Proposition~\ref{prop:singsupp}
covers all interesting cases. Indeed, if the numerical range has empty 
interior, then either $W(A)$ is reduced to a single point, in which 
case $A$ is a scalar matrix and $\mu_A$ is just a Dirac mass, 
or $W(A)$ is a line segment of nonzero length, in which case $A$ 
can be reduced to a nonscalar Hermitian matrix by a simple affine
transformation. 
\end{rque}

\bepr
Using the same notations as in Lemmas~\ref{lem:singset} and 
\ref{lem:singsupp}, we observe that $\Gamma(A) = \{x \in \partial\B^n\,|\, 
\Delta(x) = 0\}$, where $\Delta(x) = \det D(x)$ is a polynomial 
in the $2n$ variables $\Re x_i, \Im x_i$ ($i = 1,\dots,n$). Thus 
one of the following two situations must occur: 

\smallskip\noindent 1) $\Gamma(A)$ is an algebraic submanifold of
$\partial\B^n$ of codimension at least $1$. By Sard's lemma, this is 
the case if and only if $W(A)$ has nonempty interior. In that 
situation, since we already know that $\mu_A$ has a smooth 
density outside the critical set $\Sigma_A$, we only need  
to show that $\mu_A(\Sigma_A) = 0$. Given $\epsilon > 0$, let 
$\Gamma_\epsilon(A) = \{x \in \partial\B^n\,|\, \dist(x,\Gamma(A)) \le 
\epsilon\}$, where ``$\dist$'' denotes here the geodesic distance 
on the unit sphere. We decompose
\[
  \Phi_A^{-1}(\Sigma_A) \,=\, \Bigl(\Phi_A^{-1}(\Sigma_A) \cap
  \Gamma_\epsilon(A)\Bigr) \cup \Bigl(\Phi_A^{-1}(\Sigma_A) \cap
  \Gamma_\epsilon(A)^c\Bigr) \,=\, E_1(\epsilon) \cup E_2(\epsilon)~,
\]
where $\Gamma_\epsilon(A)^c = \partial\B^n\setminus\Gamma_\epsilon(A)$. 
Since $\bar\sigma(\Gamma(A)) = 0$, we have $\bar \sigma(E_1(\epsilon)) 
\le \bar\sigma(\Gamma_\epsilon(A)) \to 0$ as $\epsilon \to 0$. 
Moreover, the proof of Lemma~\ref{lem:singsupp} shows $E_2(\epsilon)$
is a codimension two submanifold of $\partial\B^n$, so that $\bar\sigma
(E_2(\epsilon)) = 0$ for any $\epsilon > 0$. Using the definition of 
the numerical measure, we conclude that $\mu_A(\Sigma_A) = \bar\sigma
(\Phi_A^{-1}(\Sigma_A)) = 0$. 

\smallskip\noindent 2) $\Gamma(A) = \partial\B^n$. This is the case
if and only if $W(A)$ has empty interior, and without loss 
of generality we can then assume that the matrix $A$ is Hermitian. 
Since $W(A) \subset \R$, it is more natural here to consider 
$\Phi_A$ as a map from $\partial\B^n$ into $\R$. If we do that, 
then repeating the proofs of Lemmas~\ref{lem:singset} and 
\ref{lem:singsupp} we easily find that the critical points
of $\Phi_A$ are exactly the eigenvectors of $A$. Moreover, the 
numerical measure has a smooth density on $\R\setminus \sigma(A)$, 
and is absolutely continuous with respect to the Lebesgue measure
$\mathrm{d}x$ if $A$ is not a scalar matrix. We skip the details
here, because the Hermitian case will be treated in full details
in Section~\ref{s:normal} below. 
\enpr

\medskip
To conclude this section, we show that the numerical density is
strictly positive in the interior of $W(A)$. A little care is
needed in the formulation of that result, because as we shall 
see in Section~\ref{ss:rad} the numerical density need not be a 
continuous function. 

\begin{prop}\label{prop:pos}
If $z_0 \in \C$ is an interior point of $W(A)$, then
\begin{equation}\label{eq:muApos}
  \liminf_{\epsilon \to 0} \frac{1}{\epsilon^2}\,\mu_A
  (\{z \in \C\,|\, |z-z_0| \le \epsilon\}) \,>\, 0~.
\end{equation}
\end{prop}

\bepr
If $z_0$ is an interior point of $W(A)$, it is shown in 
\cite[Proposition 2.11]{JS} that the preimage $\Phi_A^{-1}(z_0)$ 
contains at least one regular point $x_0$. Let $V$ be an open
geodesic ball centerd at $x_0 \in \partial\B^n$ whose closure does not
intersect $\Gamma(A)$. Proceeding as in the proof of 
Lemma~\ref{lem:singsupp}, we find
\[
  \,\mu_A(\{z \in \C\,|\, |z-z_0| \le \epsilon\})
  \,\ge\, \frac{1}{\omega_n} \int_{|z-z_0| \le \epsilon}
  \left\{\int_{\NN_z \cap V}\frac{\dd \nu(x)}{2\Delta(x)^{1/2}}
  \right\}\dd z~.
\]
If $\epsilon > 0$ is sufficiently small, the integral inside 
the curly brackets is a smooth and positive function of $z$, and
\eqref{eq:muApos} follows. 
\enpr

\subsection{The direct sum formula}\label{ss:direct}

Let $p,q \in \N^*$ and $n = p+q$. Given $A\in \bM_p(\C)$ and 
$B\in \bM_q(\C)$, the {\em direct orthogonal sum} of $A$ and $B$ 
is the matrix $A\oplus B \in \bM_n(\C)$ defined by 
\[
  A\oplus B \,=\, \begin{pmatrix} A & 0 \\ 0 & B \end{pmatrix}~.
\]
In this situation, we have a formula for the numerical measure 
$\mu_{A\oplus B}$ in terms of $\mu_A$ and $\mu_B$.

\begin{prop}\label{prop:AplusB}
For any $\phi\in C^0(\C)$, we have
\begin{equation}\label{eq:formsum}
  \int_\C \phi(z)\dd\mu_{A\oplus B}(z) \,=\, \frac{1}{\BB(p,q)}
  \int_\C\int_\C\int_0^1\phi(tz'+(1-t)z'')\,t^{p-1}(1-t)^{q-1}\dd t
  \dd\mu_A(z')\dd\mu_B(z'')~,
\end{equation}
where $\BB(p,q)$ is Euler's beta function
\[
  \BB(p,q) \,=\, \int_0^1 t^{p-1}(1-t)^{q-1}\dd t \,=\, \frac{\Gamma(p)
  \Gamma(q)}{\Gamma(p+q)}~.
\]
\end{prop}

\bepr
Any unit vector $x \in \partial\B^n$ can be written as
\[ 
  x \,=\, \begin{pmatrix} \sqrt t\,u \\[1mm] \sqrt{1-t}\,v \end{pmatrix}~,
\]
where $u \in \partial\B^p$, $v \in \partial\B^q$, and $t\in[0,1]$. 
Up to negligible sets, the map $(u,v,t)\mapsto x$ defines
a diffeomorphism from $\partial\B^p \times \partial\B^q \times [0,1]$ 
onto $\partial\B^n$. With this parametrization it is not difficult 
to verify that the Euclidean measure on $\partial\B^n$ has the following
expression
\[
  \dd\sigma_n(x) \,=\, \frac12 \,t^{p-1}(1-t)^{q-1}\dd t \dd\sigma_p(u) 
  \dd\sigma_q (v)~.
\]
Equivalently, since $\sigma_n(\partial\B^n) = \omega_n = 2\pi^n/\Gamma(n)$, 
the normalized Euclidean measure satisfies
\[
  \dd\bar\sigma_n(x) \,=\, \frac{1}{\BB(p,q)}\,t^{p-1}(1-t)^{q-1}\dd t 
  \dd\bar\sigma_p(u)\dd\bar\sigma_q (v)~.
\]
Thus, using definition \eqref{def:muA} and the fact that
$\langle (A\oplus B)x,x\rangle = t\langle Au,u\rangle + (1-t) 
\langle Bv,v\rangle$, we easily obtain
\begin{align*}
  \int_\C \phi(z)\dd\mu_{A\oplus B}(z) \,&=\, \frac{1}{\BB(p,q)}\int_{\partial\B^p}\!
  \int_{\partial\B^q}\!\int_0^1 \phi\Bigl(t\langle Au,u\rangle + (1-t) \langle 
  Bv,v\rangle\Bigr)\,t^{p-1}(1-t)^{q-1}\dd t \dd\bar\sigma_q(v) 
  \dd\bar\sigma_p(u) \\
  \,&=\, \frac{1}{\BB(p,q)}\int_\C\int_\C\int_0^1 \phi(tz'+(1-t)z'')
  \,t^{p-1}(1-t)^{q-1}\dd t \dd\mu_A(z')\dd \mu_B(z'')~,
 \end{align*}
which is the desired result. 
\enpr

\bigskip
As an application, if we choose $\phi(z) = e^{-i\xi\cdot z}$ in 
Proposition~\ref{prop:AplusB}, we obtain the following relation 
between the Fourier transforms of the measures $\mu_A$, $\mu_B$ and 
$\mu_{A\oplus B}$ :
\begin{equation}\label{eq:FAplusB}
 \hat\mu_{A\oplus B}(\xi) \,=\,  \frac{1}{\BB(p,q)}\int_0^1 \hat\mu_A(t\xi)
 \,\hat\mu_B((1-t)\xi)\,t^{p-1}(1-t)^{q-1}\dd t~, \quad \xi \in \R^2~.
\end{equation}

\medskip
The formula given in Proposition~\ref{prop:AplusB} can be generalized 
in a straightforward way to a direct sum with an arbitrary number of 
terms. Assume that $A = A_1 \oplus \dots \oplus A_k$ where $A_j \in 
\bM_{p_j}(\C)$, so that $A \in \bM_n(\C)$ with $n = p_1 + \dots + p_k$.
Setting $p = (p_1,\dots,p_k)$, we denote
\[
  \BB(p) \,=\, \int_{D_{k-1}}t_1^{p_1-1}\cdots \,t_k^{p_k-1}\dd t_1 
  \dots \dd t_{k-1} \,=\, \frac{\Gamma(p_1)\cdots\Gamma(p_k)}
  {\Gamma(p_1 + \dots + p_k)}~,
\]
where $t_k = 1 - (t_1 + \dots + t_{k-1})$ and $D_{k-1}$ denotes the 
$(k{-}1)$-dimensional simplex
\begin{equation}\label{def:Dk}
  D_{k-1} \,=\, \Bigl\{(t_1,\dots,t_{k-1}) \in \R_+^{k-1} \,\Big|\, 
  t_1 + \dots + t_{k-1} \le 1\Bigr\}~.
\end{equation}
Using \eqref{eq:formsum} and proceeding by induction over $k$,
we easily obtain the general formula
\begin{align}\nonumber
  \int_\C \phi(z)\dd\mu_A(z) \,=\, \frac{1}{\BB(p)}\int_{\C^k}
  \int_{D_{k-1}}&\phi(t_1 z_1 + \dots + t_k z_k)\,t_1^{p_1-1}\cdots\, 
  t_k^{p_k-1} \\ \label{eq:formsum2}
 &\dd t_1 \dots \dd t_{k-1} \dd\mu_{A_1}(z_1)\dots\dd\mu_{A_k}(z_k)~,
\end{align}
where it is understood again that $t_k = 1 - (t_1 + \dots + t_{k-1})$.

\section{The numerical density of a normal matrix}\label{s:normal}

If $A \in \bM_n(\C)$ is a normal matrix, the numerical measure $\mu_A$
is entirely determined by the spectrum $\sigma(A) = \{\lambda_1,\dots,
\lambda_n\}$. Indeed, we know that $A$ is unitarily equivalent to the
diagonal matrix $\diag(\lambda_1,\dots,\lambda_n)$, and that a unitary
conjugation does not affect the numerical measure. Using this observation
and the direct sum formula of Section~\ref{ss:direct}, we shall prove
that the numerical density of $A$ is a piecewise polynomial
function, which can be characterized as a multivariate $B$-spline whose
knots are the eigenvalues of $A$.  We begin with the important
particular case where all eigenvalues of $A$ are colinear.

\subsection{The Hermitian case}\label{ss:Hermitian}

If $A \in \bM_n(\C)$ is a Hermitian matrix, then $W(A) \subset 
\R$ and the numerical measure $\mu_A$ is therefore supported 
on the real axis. Assuming that $A$ is not a multiple of the identity 
matrix, we show in this section that $\mu_A$ is absolutely continuous
with respect to Lebesgue's measure on $\R$, and we give a simple
characterization of the numerical density $f_A = \mathrm{d}\mu_A /
\mathrm{d}x$. The result is:

\begin{prop}\label{prop:Bspline}
If $A \in \bM_n(\C)$ is a nonscalar Hermitian matrix, the 
numerical density $f_A : \R \to \R_+$ is the normalized 
$B$-spline of degree $n-2$ whose knots are the eigenvalues 
of $A$.
\end{prop}

To make the statement clear, we briefly recall the definition 
and some elementary properties of the classical $B$-splines 
\cite{dB}. If $\lambda_1,\dots,\lambda_n \in \R$ are pairwise
distinct, the {\em $(n{-}1)^{th}$ divided difference} of a continuous 
function $g$ at the points $\lambda_1,\dots,\lambda_n$ is the quantity
\begin{equation}\label{def:dd1}
  \delta^{n-1}[\lambda_1,\dots,\lambda_n]g \,=\, \sum_{j=1}^n 
  \frac{g(\lambda_j)}{\prod_{k\neq j}(\lambda_j-\lambda_k)}~.
\end{equation}
This is the leading coefficient of the unique polynomial of degree at
most $n{-}1$ which agrees with $g$ at the points $\lambda_1,\dots,\lambda_n$.
It is easy to verify that the right-hand side of \eqref{def:dd1} is 
a completely symmetric function of the variables $\lambda_j$. If 
$g \in C^{n-1}(\R)$, the divided difference can be extended by continuity 
to arbitrary (not necessarily distinct) values of $\lambda_1,\dots,
\lambda_n$, and we have the integral formula:
\begin{equation}\label{def:dd2}
  \delta^{n-1}[\lambda_1,\dots,\lambda_n]g \,=\, \int_{D_{n-1}}
  g^{(n-1)}(t_1\lambda_1 + \dots + t_n\lambda_n)\dd t_1\dots \dd t_{n-1}~,
\end{equation}
where $D_{n-1}$ is the $(n{-}1)$-dimensional simplex defined in \eqref{def:Dk}
and $t_n = 1 - (t_1 + \dots + t_{n-1})$. In what follows, we shall 
always assume that the set $S = \{\lambda_1,\dots,\lambda_n\}$ is not
reduced to a single point, so that the $(n{-}1)^{th}$ divided 
difference is well-defined as soon as $g$ is of class $C^{n-2}$ in a 
neighborhood of $S$.  

With these notations, the {\em normalized $B$-spline of degree $n-2$ 
with knots $\lambda_1,\dots,\lambda_n$} is the function $B : \R 
\to \R$ defined by
\begin{equation}\label{def:Bspline}
  B(x) \,\equiv\, B[\lambda_1,\dots,\lambda_n](x) \,=\, (n-1)
  \,\delta^{n-1}[\lambda_1,\dots,\lambda_n](\cdot - x)_+^{n-2}~, \qquad
  x \in \R~,
\end{equation}
where $(\cdot - x)_+^{n-2}$ denotes the map $y \mapsto\max(0,y-x)^{n-2}$.
If $\lambda_1 \le \dots \le \lambda_n$, it is not
difficult to show that $B(x)$ vanishes identically outside
$[\lambda_1,\lambda_n]$, and coincides with a polynomial of degree at
most $n-2$ on each nonempty interval $(\lambda_j,\lambda_{j+1})$.
Moreover, if $m_j$ denotes the multiplicity of $\lambda_j$ in $S$, 
one can verify that $B(x)$ is continuous at $x = \lambda_j$ together
with its derivatives up to order $d_j = n-2-m_j$, provided $d_j \ge 0$. 
If $d_j = -1$, then $B(x)$ is discontinuous at $\lambda_j$. Finally, 
we shall see below that $B(x)$ is positive on $(\lambda_1,\lambda_n)$ 
and that $\int_\R B(x)\dd x = 1$.  

\bigskip\noindent{\bf Proof of Proposition~\ref{prop:Bspline}.}
Let $A \in \bM_n(\C)$ be a normal matrix with eigenvalues $\lambda_1,
\dots,\lambda_n$. To compute the numerical density, we can assume 
without loss of generality that $A$ is diagonal, namely $A = A_1 \oplus 
\dots \oplus A_n$ with $A_j = \lambda_j \in \bM_1(\C)$. Thus we can 
use the direct sum formula \eqref{eq:formsum2} with $k = n$ and 
$p = (1,\dots,1)$. Since $\BB(p)^{-1} = (n-1)!$ and $\mu_{A_j} = 
\delta_{\lambda_j}$ for $j = 1,\dots,n$, we obtain the relation
\begin{equation}\label{eq:normrep}
  \int_\C \phi(z)\dd\mu_A(z) \,=\, (n-1)! \int_{D_{n-1}} 
  \phi(t_1 \lambda_1 + \dots + t_n \lambda_n) \dd t_1 \dots 
  \dd t_{n-1}~,
\end{equation}
for any continuous function $\phi : \C \to \C$. 

Assume now that $A$ is Hermitian, so that $\lambda_1,\dots,\lambda_n 
\in \R$, and that $\sigma(A)$ is not reduced to a single point. 
If $B = B[\lambda_1,\dots\lambda_n]$ is the normalized $B$-spline 
defined by \eqref{def:Bspline}, we claim that
\begin{equation}\label{eq:normrep2}
   \int_\R \phi(x) B[\lambda_1,\dots,\lambda_n](x)\dd x \,=\, 
  (n-1)! \int_{D_{n-1}} \phi(t_1 \lambda_1 + \dots + t_n \lambda_n) 
  \dd t_1 \dots \dd t_{n-1}~.
\end{equation}
Indeed, it is clearly sufficient to prove \eqref{eq:normrep2} for
compactly supported functions $\phi \in C^0(\R)$. Moreover, since
both members of \eqref{eq:normrep2} depend continuously on
$\lambda_1,\dots,\lambda_n$, we can also assume that the eigenvalues
of $A$ are all distinct. In that case, it follows immediately from
\eqref{def:dd1}, \eqref{def:Bspline} that
\[
  \int_\R \phi(x) B[\lambda_1,\dots,\lambda_n](x)\dd x \,=\, 
  (n-1)\,\delta^{n-1}[\lambda_1,\dots,\lambda_n] \Phi~, 
\]
where $\Phi(y) = \int_\R \phi(x)(y - x)_+^{n-2}\dd x$. Since $\Phi^{(n-1)} 
= (n-2)!\,\phi$, we can use the integral formula \eqref{def:dd2} to evaluate
the divided difference in the right-hand side, and we obtain
\eqref{eq:normrep2}. 

Now, comparing \eqref{eq:normrep} and \eqref{eq:normrep2}, we conclude
that $\mu_A$ is absolutely continuous with respect to Lebesgue's 
measure on $\R$, and that the numerical density $f_A = \mathrm{d}
\mu_A /\mathrm{d}x$ is precisely the normalized $B$-spline 
$B[\lambda_1,\dots,\lambda_n]$.  Incidentally, the argument above
shows that $B$ is positive on its support (see Proposition~\ref{prop:pos}) 
and that $\int_\R B(x)\dd x = 1$.  
\enpr

\subsection{The quasi-Hermitian case}\label{ss:quasi-Hermitian}

We say that a matrix $A \in \bM_n(\C)$ is {\em quasi-Hermitian}
if the numerical range $W(A) \subset \C$ has empty interior, i.e.
$W(A)$ is a single point or a line segment. In such a case, there
exist $z \in \C$, $\theta \in S^1$, and a Hermitian matrix
$H$ such that $A = zI_n + e^{i\theta}H$ (in particular, $A$ is 
normal). Indeed, if we choose $z,\theta$ such that $W(A) \subset
z + e^{i\theta}\R$, the matrix $H = e^{-i\theta}(A-zI_n)$ satisfies
$W(H) \subset \R$ and is therefore Hermitian. Since $\langle Ax,x
\rangle = z + e^{i\theta}\langle Hx,x\rangle$ for all $x \in 
\partial\B^n$, it is clear from \eqref{def:muA} that the numerical
measure $\mu_A$ is just the image of $\mu_H$ under the affine isometry
$w \mapsto z + e^{i\theta}w$. Combining this remark with 
Proposition~\ref{prop:Bspline}, we thus obtain a precise characterization 
of the numerical measure of any quasi-Hermitian matrix.

\subsection{The normal case}\label{ss:normal}

Finally, we consider the case of a normal matrix $A \in \bM_n(\C)$ 
whose numerical range $W(A)$ has nonempty interior. This, of course, 
is possible only if $n \ge 3$. By Proposition~\ref{prop:singsupp}, 
the numerical measure $\mu_A$ is absolutely continuous with respect
to Lebesgue's measure, and we have the following characterization 
of the numerical density $f_A = \mathrm{d}\mu_A /\mathrm{d}z$:

\begin{prop}\label{prop:Bspline2}
If $A \in \bM_n(\C)$ is a normal matrix whose numerical range $W(A)
\subset \C$ has nonempty interior, the numerical density $f_A : \C 
\to \R_+$ is the bivariate $B$-spline of degree $n-3$ whose knots are 
the eigenvalues of $A$.
\end{prop}

The reader is referred here to the work of W.~Dahmen \cite{Dah}, 
where multivariate $B$-splines are defined and studied in detail. 
To make the connection with the numerical density of a normal 
matrix, we use the relation \eqref{eq:normrep}, which corresponds 
to formula~(2.2) in \cite{Dah}. In the rest of this section, we
assume that $A \in \bM_n(\C)$ is a normal matrix whose eigenvalues
$\lambda_1,\dots,\lambda_n$ are not colinear, and we often identify 
the complex plane $\C$ with $\R^2$. 

\bigskip\noindent{\bf Proof of Proposition~\ref{prop:Bspline2}.}
Assume first that $n = 3$. Then $W(A) \subset \C \simeq \R^2$ is the 
$2$-simplex with vertices $\lambda_1,\lambda_2, \lambda_3$, and 
using the change of variables $w = t_1\lambda_1 + t_2\lambda_2 + 
(1-t_1-t_2)\lambda_3$ in \eqref{eq:normrep} we easily obtain
\[
  \int_\C \phi(z)\dd\mu_A(z) \,=\, \frac{1}{|W(A)|}
  \int_{W(A)} \phi(w)\dd w~.
\]
This shows that the numerical measure $\mu_A$ is uniformly distributed
on $W(A)$. The numerical density is thus a multiple of the characteristic
function of $W(A)$, which (by definition) is the bivariate $B$-spline 
of degree zero with knots $\lambda_1, \lambda_2,\lambda_3$.

\medskip
We now assume that $n \ge 4$. Then we can choose $n$ vectors 
$v_1,\dots,v_n \in \R^{n-1}$ such that\\[2pt]
1) The $(n{-}1)$-simplex $S \subset \R^{n-1}$ with vertices 
$v_1,\dots,v_n$ has unit volume; \\[2pt]
2) $Pv_i = \lambda_i$ for $i = 1,\dots,n$, where $P : \R^{n-1} \to \R^2$ 
is defined by $P(x_1,\dots,x_{n-1}) = (x_1,x_2)$.\\[2pt]
This (elementary) claim is proved in \cite[Section~2]{Dah}. The 
simplex $S \subset \R^{n-1}$ is not uniquely defined, but any choice
satisfies $PS = W(A) = \co(\sigma(A))$. Returning to \eqref{eq:normrep}, 
we have
\[
  \int_\C \phi(z)\dd\mu_A(z) \,=\, (n-1)! \int_{D_{n-1}} 
  \phi(P(t_1 v_1 + \dots + t_n v_n)) \dd t_1 \dots 
  \dd t_{n-1} \,=\, \int_S \phi(Pw)\dd w~,
\]
where the second equality is obtained by applying the change of 
variables $(t_1,\dots,t_{n-1}) \mapsto w = t_1 v_1 + \dots + t_n v_n$, 
with $t_n = 1-(t_1+\dots+t_{n-1})$. This shows that the numerical 
measure is the image under the projection $P$ of the Lebesgue measure 
on the simplex $S \subset \R^{n-1}$. Given $z \in \C \simeq \R^2$, 
the numerical density $f_A(z)$ is thus the $(n{-}3)$-dimensional 
measure of the simplex $P^{-1}z \cap S$. This is precisely the 
definition given in \cite{Dah} of a bivariate $B$-spline of degree
$n-3$ with knots $\lambda_1,\dots,\lambda_n$. 
\enpr

\bigskip We conclude by listing a few properties of the numerical
density $f_A$, which follow from \cite[Theorem~4.1]{Dah}. Let
$\Sigma_A \subset \C$ be the union of all line segments joining pairs
of eigenvalues of $A$. Using Lemma~\ref{lem:singset}, it is
straightforward to verify that $\Sigma_A$ is exactly the set of
critical values of the numerical map $\Phi_A$. Then $f_A$ is a
polynomial of total degree $n-3$ in each connected component of $\C
\setminus \Sigma_A$. Moreover, if $\ell \subset \C$ is a straight line
passing through a pair of eigenvalues of $A$, the numerical density is
continuous across $\ell \cap W(A)$ together with its derivatives up to
order $n-2-m \ge 0$, where $m \ge 2$ is the number of eigenvalues of
$A$ (counted with multiplicities) which belong to $\ell$. If $m =
n-1$, then $f_A$ is discontinuous on $\ell \cap W(A)$. In particular,
in the generic case where no straight line contains more than two
eigenvalues of $A$, the numerical density $f_A$ is of class $C^{n-4}$ 
if $n \ge 4$.

\paragraph{Remark.} In the Hermitian and normal cases, the fact that
the numerical density $f_A$ is a projection of the characteristic 
function of a convex set, which is a log-concave function, implies 
that the density is itself a log-concave function. We thus have 
the following inequality
\begin{equation}\label{eq:logccve}
  f_A(\lambda z + (1-\lambda)z') \,\ge\, f_A(z)^\lambda 
  \,f_A(z')^{1-\lambda}~, 
\end{equation}
for all $z,z'\in \C$ and all $\lambda\in(0,1)$. This follows from the 
Pr\'ekopa--Leindler inequality, see e.g. \cite[Section~9]{Grd}. 
We warn the reader that this property does not extend to nonnormal 
matrices, as we can see already from the two-dimensional case
considered in Section~\ref{ss:twod}.

\section{The Radon transform of the numerical measure}
\label{s:Radon}

Our purpose in this section is to derive a representation formula
for the numerical density of a nonnormal matrix $A \in \bM_n(\C)$. 
Our approach is based on a natural expression of the {\em Radon 
transform} of the numerical measure $\mu_A$ in terms of the 
Hermitian matrices $H(\theta)$ defined in \eqref{def:Htheta}. 
By definition, the Radon transform of $\mu_A$ is the family $\RR\mu_A 
= \{\RR\mu_A(\theta)\,|\,\theta \in S^1\}$, where $\RR\mu_A(\theta)$ 
denotes the Borel measure on $\R$ defined by
\[
  (\RR\mu_A(\theta))(I) \,=\, \mu_A(\{z \in \C\,|\, 
  \Re\,(e^{-i\theta}z) \in I\})~, 
\]
for any open set $I \subset \R$. In other words, $\RR\mu_A(\theta)$
is the image of the measure $\mu_A$ under the orthogonal projection 
in $\C \simeq \R^2$ onto the line $e^{i\theta}\R$. The fundamental 
observation is:

\begin{prop}\label{prop:Radon}
For any $\theta \in S^1 = \R/(2\pi\Z)$, one has $\RR\mu_A(\theta) = 
\mu_{H(\theta)}$. 
\end{prop}

\bepr
The definition \eqref{def:Htheta} implies that $\Re\,(e^{-i\theta}
\langle Ax,x\rangle) = \langle H(\theta)x,x\rangle$ for any 
$x \in \partial\B^n$. Thus, for any open set $I \subset \R$, we have
\begin{align*}
  (\RR\mu_A(\theta))(I) \,&=\, \bar\sigma(\{x \in \partial\B^n\,|\, 
  \Re\,(e^{-i\theta}\langle Ax,x\rangle) \in I\}) \\ 
  \,&=\, \bar\sigma(\{x \in \partial\B^n\,|\, 
  \langle H(\theta)x,x\rangle \in I\}) \,=\, \mu_{H(\theta)}(I)~,
\end{align*}
which proves the claim. 
\enpr

\medskip Proposition~\ref{prop:Radon} shows that the numerical measure
of an arbitrary matrix $A \in \bM_n(\C)$ is entirely determined by the
one-dimensional measures associated with the Hermitian matrices
$\{H(\theta)\,|\,\theta \in S^1\}$. If we assume that $W(A)$ has
nonempty interior, which is always the case if $A$ is nonnormal, the
matrix $H(\theta)$ is nonscalar for every $\theta \in S^1$ and
Proposition~\ref{prop:Bspline} show that its numerical density is the
$B$-spline $B[\lambda_1(\theta),\dots,\lambda_n(\theta)]$, where
$\lambda_1(\theta), \dots,\lambda_n(\theta)$ are the eigenvalues of
$H(\theta)$. In that case, the result of Proposition~\ref{prop:Radon}
can be stated in the following equivalent form
\begin{equation}\label{eq:Radon1}
  \int_\R f_A(e^{i\theta}(x+iy))\dd y \,=\, B[\lambda_1(\theta),\dots,
  \lambda_n(\theta)](x)~, 
\end{equation}
where equality holds for all $\theta \in S^1$ and almost all $x \in \R$
since the numerical density $f_A$ belongs to $L^1(\C)$. 

Our goal is to invert the Radon transform \eqref{eq:Radon1} 
to obtain a representation formula for the numerical density $f_A$. 
The general results established in \cite{Hel} show that
\begin{equation}\label{eq:Radon2}
  (J f_A)(x+iy) \,=\, \frac{1}{4\pi} \int_{S^1} 
  B[\lambda_1(\theta),\dots,\lambda_n(\theta)](x\cos\theta + y\sin\theta)
  \dd\theta~, 
\end{equation}
where $J = (-\Delta)^{-1/2}$ is the Riesz potential defined by
$(Jf)(z) = \frac{1}{2\pi}\int_\C |z-z'|^{-1}f(z')\dd z'$. The idea is
thus to apply the nonlocal operator $(-\Delta)^{1/2}$ to both sides of
\eqref{eq:Radon2}, but since we are not dealing with smooth functions
we have to differentiate in the sense of distributions.  As a
preliminary remark, if $f(x,y) = g(x\cos\theta + y\sin\theta)$ for
some test function $g : \R \to \R$ and some fixed $\theta \in S^1$,
a direct calculation shows that $(\Delta f)(x,y) =
g''(x\cos\theta + y\sin\theta)$, and a standard interpolation argument
allows us to conclude that $(-\Delta)^{1/2}f(x,y) = \HH g'(x\cos\theta
+ y\sin\theta)$, where $\HH g'$ denotes the {\em Hilbert transform} of
the derivative $g'$.  Using this observation, we easily obtain the
representation formula
\begin{equation}\label{eq:Radon3}
  f_A(x+iy) \,=\, \frac{1}{4\pi} \int_{S^1} \HH B'[\lambda_1(\theta),\dots,
  \lambda_n(\theta)](x\cos\theta + y\sin\theta)\dd\theta~,
\end{equation}
where both sides define integrable functions of $z = x+iy \in \C$, 
and equality holds almost everywhere. So we have shown:

\begin{prop}\label{prop:RadonInv}
If the numerical range of a matrix $A \in \bM_n(\C)$ has 
nonempty interior, the numerical density of $A$ can be represented
as in \eqref{eq:Radon3} for almost all $(x,y) \in \R^2$. 
\end{prop}

In the rest of this section, we shall apply Proposition~\ref{prop:RadonInv} 
to compute the numerical density of a two-dimensional nonnormal matrix. 
We shall also consider the interesting particular situation where the 
numerical density is radially symmetric, in which case the representation 
formula takes a simpler form. Proposition~\ref{prop:RadonInv} will be 
used again in Section~\ref{s:main} to derive some important properties of 
the numerical measure in the general case. 

\subsection{The two-dimensional case}\label{ss:twod}

Let $A \in \bM_2(\C)$, and assume that the numerical range 
$W(A)$ has nonempty interior. As is well-known \cite{HJ}, 
$W(A)$ is then a filled ellipse, and without loss of generality
we can assume that this ellipse is centered at the origin and that its 
major axis is aligned with the real axis of the complex plane. 
In that case, up to a unitary conjugation, the matrix $A$ has the 
following form
\begin{equation}\label{def:A2}
  A \,=\,\begin{pmatrix} -c & 2b \\ 0 & c\end{pmatrix}~, \qquad
  W(A) \,=\, \Bigl\{x+iy \in \C \,\Big|\, \frac{x^2}{a^2} + 
  \frac{y^2}{b^2} \le 1 \Bigr\}~,
\end{equation}
where $b > 0$, $c \ge 0$, and $a = (b^2 + c^2)^{1/2}$. For 
$\theta \in S^1 = \R/(2\pi\Z)$, let
\[
  H(\theta) \,=\, \frac{1}{2}(e^{-i\theta}A + e^{i\theta}A^*)
  \,=\, \begin{pmatrix} -c\cos\theta & b\,e^{-i\theta} \\ 
  b\,e^{i\theta} & c\cos\theta \end{pmatrix}~.
\]
The eigenvalues of $H(\theta)$ are $\pm \lambda(\theta)$, where 
$\lambda(\theta) = (b^2 + c^2\cos^2\theta)^{1/2}$. It follows that
\[
  B[-\lambda(\theta),\lambda(\theta)](s) \,=\, \frac{1}{2\lambda(\theta)}
  \mathbf{1}_{[-\lambda(\theta),\lambda(\theta)]}(s)~.
\]
According to \eqref{eq:Radon3}, we have to differentiate this expression
with respect to $s$, which yields a linear combination of Dirac masses, and 
to apply the Hilbert transformation, which is the convolution with the 
distribution $\pv \frac{1}{\pi s}$. We obtain
\[
  \HH B'[-\lambda(\theta),\lambda(\theta)](s) \,=\, 
  \frac{1}{2\pi\lambda(\theta)}\Bigl(\pv\frac{1}{s+\lambda(\theta)}
  - \pv\frac{1}{s-\lambda(\theta)}\Bigr) \,=\,
  \frac{1}{\pi}\,\pv\frac{1}{\lambda(\theta)^2-s^2}~,
\]
hence
\begin{equation}\label{eq:twod1}
  f_A(x+iy) \,=\, \frac{1}{4\pi^2}\,\pv\int_{S^1} \frac{1}{b^2 + 
  c^2\cos^2\theta - (x\cos\theta + y\sin\theta)^2}\dd\theta~.
\end{equation}

It remains to compute the right-hand side of \eqref{eq:twod1}, 
which is a simple exercise in complex analysis. Setting $z = x+iy$ 
and $w = e^{2i\theta}$, we first observe that 
$$
  b^2 + c^2\cos^2\theta -  (x\cos\theta + y\sin\theta)^2 \,=\, \frac{1}{4w}
  \Bigl((c^2-z^2) + 2w(a^2+b^2-|z|^2) + w^2(c^2-\bar z^2)\Bigr)~,
$$
hence
\begin{equation}\label{eq:twod2}
  f_A(z) \,=\, \frac{1}{i\pi^2}\,\pv\oint_{|w| = 1} \frac{1}{(c^2-z^2) 
  + 2w(a^2+b^2-|z|^2) + w^2(c^2-\bar z^2)} \dd w~. 
\end{equation}
In \eqref{eq:twod2}, the roots of the denominator are
\begin{align}\nonumber
  w_\pm \,&=\, \frac{1}{c^2-\bar z^2}\Bigl(-a^2-b^2+|z|^2 \mp
  \sqrt{(a^2+b^2-|z|^2)^2-(c^2-z^2)(c^2-\bar z^2)}\Bigr) \\ 
  \label{eq:sqroot}
  \,&=\, \frac{1}{c^2-\bar z^2}\Bigl(-a^2-b^2+|z|^2 \mp
  \sqrt{4(a^2b^2 - b^2x^2 - a^2y^2)}\Bigr)~.
\end{align}
We can therefore distinguish between two cases:

\smallskip\noindent{\bf 1.}
The point $z = x+iy$ belongs to the interior of $W(A)$. Then the 
expression under the square root is positive, and it is easy to 
verify that $|w_+| > 1$, $|w_-| < 1$ (in the limiting case where
$z = \pm c$ is a focus of the ellipse, one can set $w_- = 0$ and 
$w_+ = \infty$.) Thus the principal value in \eqref{eq:twod2}
is not needed, and the residue theorem shows that
\begin{equation}\label{eq:twod3}  
  f_A(z) \,=\, \frac{1}{\pi}\frac{1}{(a^2+b^2-|z|^2) + (c^2-\bar z^2)w_-}
  \,=\, \frac{1}{2\pi}\frac{1}{\sqrt{a^2b^2 - b^2x^2 - a^2y^2}}~. 
\end{equation}

\smallskip\noindent{\bf 2.} The point $z = x+iy$ lies outside $W(A)$. 
Then the expression under the square root in \eqref{eq:sqroot} is
negative, and one verifies that $|w_\pm| = 1$. Thus the integrand in
\eqref{eq:twod2} is holomorphic outside the unit circle $\{|w| = 1\}$
and decreases like $1/|w|^2$ at infinity. It follows that
\[
  I_z(r) \,:=\, \frac{1}{i\pi^2}\oint_{|w| = r} \frac{1}{(c^2-z^2) 
  + 2w(a^2+b^2-|z|^2) + w^2(c^2-\bar z^2)} \dd w \,=\, 0~,
\]
for all $r \neq 1$, hence $f_A(z) = \frac12(I_z(1+) + I_z(1-)) = 0$. 

\smallskip
Summarizing, we have shown that the numerical density of the 
matrix \eqref{def:A2} is the function $f_A \in L^1(\C)$ defined by 
\eqref{eq:twod3} inside $W(A)$, and vanishing identically 
outside $W(A)$. Note that $\singsupp(\mu_A) = \partial W(A) = 
\Sigma_A$, in agreement with Lemma~\ref{lem:singsupp}, 
and that $f_A(z)$ blows up when $z$ converges to the boundary 
of $W(A)$ from inside. In particular $f_A$ is not log-concave, 
in contrast to what happens when $A$ is normal. 

\subsection{The radially symmetric case}\label{ss:rad}

It sometimes happens that the numerical range of a matrix $A \in
\bM_n(\C)$ is a disk in the complex plane and that the numerical
density is radially symmetric about the center. In such a case, it is
possible to obtain a representation formula which is simpler than 
\eqref{eq:Radon3}. Indeed, assume that the disk $W(A)$ is centered 
at the origin, and let $R > 0$ denote the numerical radius of $A$ 
(or any larger positive number). We set
\begin{equation}\label{def:FA}
  f_A(z) \,=\, F_A(R^2 - |z|^2)~, \qquad z \in \C~, \quad |z| \le R~,
\end{equation}
where $F_A : [0,R^2] \to \R_+$ has to be determined. Since the 
numerical measure is invariant under rotations about the origin, 
it follows from Proposition~\ref{prop:Radon} that the projected 
measure $\mu_{H(\theta)}$ does not depend on $\theta$. In analogy 
with \eqref{def:FA}, if $f_H$ denotes the numerical density of 
$H(\theta)$ for any $\theta \in S^1$, we set 
\begin{equation}\label{def:FH}
  f_H(x) \,=\, F_H(R^2 - x^2)~, \qquad x \in \R~, \quad |x| \le R~.
\end{equation}
We then have the following result:

\begin{prop}\label{prop:radial}
If the numerical density of a nonscalar matrix $A \in \bM_n(\C)$ is 
radially symmetric, the functions $F_A$, $F_H$ defined in 
\eqref{def:FA}, \eqref{def:FH} satisfy the relations
\begin{equation}\label{eq:FAFH}
  F_H(t) \,=\, \int_0^t F_A(t-s)\frac{1}{\sqrt s}\dd s~, \qquad
  F_A(s) \,=\, \frac{1}{\pi}\,\frac{\mathrm{d}}{\mathrm{d}s} 
  \int_0^s F_H(s-t)\frac{1}{\sqrt t}\dd t~,
\end{equation}
for $t,s \in [0,R^2]$. 
\end{prop}

\bepr
By the definition of the Radon transformation, we have $f_H(x) = \int_\R 
f_A(x+iy)\dd y$ for all $x \in \R$. Using \eqref{def:FA}, \eqref{def:FH} 
and the support property, we thus find 
\[ 
  F_H(R^2-x^2) \,=\, \int_{y^2 \le R^2-x^2}F_A(R^2-x^2-y^2)\dd y
  \,=\, \int_0^{R^2-x^2}F_A(R^2-x^2-s)\frac{1}{\sqrt s}\dd s~,
\]
for all $x \in [-R,R]$. Setting $t = R^2-x^2$, we obtain the first 
relation in \eqref{eq:FAFH}. So far, we have shown that $F_H = \pi^{1/2}
I F_A$, where $I$ is the Riesz potential
\[
  (If)(t) \,=\, \frac{1}{\sqrt{\pi}}\int_0^tf(t-s)\frac{1}
  {\sqrt s}\dd s~, \qquad t > 0~.
\]
Now, it is well known that $(I^2 f)(t) = \int_0^t f(s)\dd s$, see 
e.g. \cite[Chapter~V.5]{Hel}, thus the second relation in \eqref{eq:FAFH}
follows from the first one. \enpr

\paragraph{Examples.} Let $A = (a_{ij})$ be a complex matrix, and assume 
that there exists a nonzero integer $k$ such that $a_{ij} = 0$ 
whenever $j-i \neq k$. Then $W(A)$ is a disk centered at the origin, 
and the numerical density $f_A$ is radially symmetric. Indeed, 
given any $\theta \in S^1$, the map
\[
  x \,=\, \begin{pmatrix} x_1 \\ \vdots \\ x_n\end{pmatrix}
  ~\mapsto~ x_\theta \,=\, \begin{pmatrix} e^{i\theta }x_1 \\ \vdots 
  \\ e^{in\theta}x_n\end{pmatrix}
\]
is a measure-preserving isomorphism of the unit sphere $\partial\B^n$, 
and our assumption on $A$ implies that $\langle Ax_\theta,x_\theta\rangle
= e^{ik\theta}\langle Ax,x\rangle$ for all $x \in \partial\B^n$. This 
proves that a rotation of angle $k\theta$ about the origin does 
not affect the numerical measure $\mu_A$. Since $\theta$ is arbitrary, 
the numerical density $f_A$ is necessarily radially symmetric. 

The simplest example in this category is the $2\times2$ Jordan
block
\begin{equation}\label{eq:defA2}
   A_2 \,=\, \begin{pmatrix} 0 & 2 \\ 0 & 0 \end{pmatrix}~.
\end{equation}
Here $W(A_2) = \{z \in \C\,|\, |z| \le 1\}$, and applying \eqref{eq:twod3}
with $a = b = 1$ we find
\begin{equation}\label{eq:fA2}
  f_{A_2}(z) \,=\, \frac{1}{2\pi}\,\frac{1}{\sqrt{1-|z|^2}}
  \,\mathbf{1}_{\{|z| < 1\}}~.
\end{equation}
Alternatively, since the eigenvalues of $H_2 = \frac12(A_2+A_2^*)$ are
$\pm 1$, we have $f_{H_2} = \frac12 \,\mathbf{1}_{[-1,1]}$ and 
applying Proposition~\ref{prop:radial} we easily obtain \eqref{eq:fA2}. 

As a more interesting application, consider the $3\times3$ matrix
\begin{equation}\label{eq:defA3}
  A_3 \,=\, \begin{pmatrix} 0 & a & 0 \\ 0 & 0 & b \\ 0 & 0 & 0 
  \end{pmatrix}~,
\end{equation}
where $a,b \in \C$ and $|a| + |b| > 0$. Multiplying $A_3$ with a 
positive constant, we can assume that $|a|^2 + |b|^2 = 4$. Then 
$W(A_3) = \{z \in \C\,|\, |z| \le 1\}$, and 
\begin{equation}\label{eq:fA3}
  f_{A_3}(z) \,=\, \frac{1}{\pi}\,\log\frac{1+\sqrt{1-|z|^2}}{|z|}
  \,\mathbf{1}_{\{0 < |z| < 1\}}~.
\end{equation}
Indeed, the eigenvalues of $H_3 = \frac12(A_3+A_3^*)$ are $0$ and $\pm
1$, hence $f_{H_3}(x) = B[-1,0,1](x) = (1{-}|x|)_+$. Applying
Proposition~\ref{prop:radial}, we obtain \eqref{eq:fA3} by a
straightforward calculation. Alternatively, if we take $(a,b) = 
(0,2)$, we can use the direct sum formula \eqref{eq:formsum}
to deduce \eqref{eq:fA3} from \eqref{eq:fA2}. Here again, we have
$\singsupp(\mu_{A_3}) = \partial W(A_3) \cup \{0\} \,=\,
\Sigma_{A_3}$, in agreement with Lemma~\ref{lem:singsupp}.  Remark
that $f_{A_3}(z) = \OO((1{-}|z|)^{1/2})$ as $|z| \to 1-$, so that the
singularity of the numerical measure at the boundary is weaker than it
was for $A_2$. This reflects the fact that $f_{H_3}(x) = (1-|x|)_+$ is
Lipschitz continuous, whereas $ f_{H_2}(x) = \mathbf{1}_{[-1,+1]}(x)$
had jump discontinuities.  However, we observe that the (logarithmic)
singularity of $f_{A_3}$ at the origin is much stronger than the
(square root) singularity at the boundary. As we shall see in
Section~\ref{s:geom}, this is a nongeneric concentration phenomenon
due to the fact that the component of the critical set $\Sigma_{A_3}$
associated with the eigenvalue $0$ of $H(\theta)$ is reduced to a single
point. We have here the rare instance of an unbounded numerical
density for a matrix of size $n \ge 3$.

Another interesting conclusion that can be drawn from this example is
that the numerical density $\mu_A$ {\em does not} determine the matrix
$A \in \bM_n(\C)$ up to unitary conjugations if $n \ge 3$.  Indeed, if
we set $(a,b) = (0,2)$ and $(a,b) = (\sqrt{2},\sqrt{2})$ in the
definition \eqref{eq:defA3} of $A_3$, the resulting matrices are not
even similar, yet they have the same numerical density, given by
\eqref{eq:fA3}.

\section{The geometry of the singular set}\label{s:geom}

We know from Lemma~\ref{lem:singsupp} that the numerical density of 
a matrix $A \in \bM_n(\C)$ is smooth outside the set $\Sigma_A$ 
of all critical values of the numerical map $\Phi_A$. In this
section, we describe a few geometrical properties of the singular
set $\Sigma_A$ which will be needed to formulate our main results
in Section~\ref{s:main}. We do not claim much originality here: 
the material of this section is essentially borrowed from 
\cite{Kip,BY,JAG,JS}, and is reproduced below for the reader's 
convenience. 

As was shown by Kippenhahn \cite{Kip}, the singular set $\Sigma_A$ has
a natural description in terms of the eigenvalues $\lambda_1(\theta),
\dots,\lambda_n(\theta)$ of the Hermitian matrices $H(\theta)$ defined 
in \eqref{def:Htheta}. To see that, we first recall that these eigenvalues 
can be numbered in such a way that they are real-analytic functions of 
$\theta \in \R$, see \cite{YB}. By analyticity, any two eigenvalues
either coincide for all $\theta \in \R$ or cross at most a finite 
number of times on each compact interval. As a consequence, there exists
an integer $m \le n$ such that $H(\theta)$ has exactly $m$ distinct
eigenvalues for all $\theta \in [0,\pi)\setminus\Theta$, where
$\Theta \subset [0,\pi)$ is a finite set. Since $H(\theta+\pi) = 
-H(\theta)$, it follows that we can number the eigenvalues in such 
a way that
\[
  \sigma(H(\theta)) \,=\, \{\lambda_1(\theta),\dots,\lambda_m(\theta)\}
  ~, \qquad \hbox{for all }\theta \in \R~,
\]
where $\lambda_1(\theta),\dots,\lambda_m(\theta)$ are pairwise 
distinct and have constant multiplicities outside the crossing 
set $\Theta + \pi\Z$. Let $\tau : \{1,\dots,m\} \to \{1,\dots,m\}$
be the permutation defined by
\begin{equation}\label{def:tau}
  \{\lambda_{\tau(1)}(\theta),\dots,\lambda_{\tau(m)}(\theta)\} 
  \,=\, \{-\lambda_1(\theta+\pi),\dots,-\lambda_m(\theta+\pi)\}~,
\end{equation}
for any $\theta \in [0,\pi)\setminus\Theta$ (hence for all $\theta \in
\R$). If we decompose $\tau$ into disjoint cycles $\CC_1,\dots,\CC_k$,
we can associate to each cycle $\CC_J$ its {\em length} $\ell_J$ and
its {\em multiplicity} $m_J$, the latter being defined as the 
multiplicity of $\lambda_j(\theta)$ as an eigenvalue of $H(\theta)$ 
for any $j \in \CC_J$. By construction, we have
\[
  \ell_1 + \dots + \ell_k \,=\, m~, \qquad \hbox{and} \quad
  \ell_1 m_1 + \dots + \ell_k m_k \,=\, n~.
\]
Moreover, if $j \in \CC_J$, then $\lambda_j(\theta+\ell_J\pi) = 
(-1)^{\ell_J}\lambda_j(\theta)$, hence $\lambda_j(\theta)$ is 
periodic with period $\ell_J \pi$ if $\ell_J$ is even and $2\ell_J\pi$ 
if $\ell_J$ is odd. Note however that these periods are not necessarily 
minimal. 

Now, we associate to each cycle $\CC_J$ of the permutation $\tau$ 
a closed curve $C_J \subset \C$ defined by
\begin{equation}\label{def:CJ}
  C_J \,=\, \Bigl\{e^{i\theta}(\lambda_j(\theta) + i\lambda_j'(\theta)) 
 \,\Big|\, \theta \in \R\Bigr\}~,
\end{equation}
where $j \in \{1,\dots,m\}$ is any element of the cycle $\CC_J$. 
Since $\lambda_j(\theta+2\ell_J\pi) = \lambda_j(\theta)$, it is 
clear that $C_J$ is indeed a closed curve, and the definition 
of the permutation $\tau$ shows that the right-hand side of 
\eqref{def:CJ} does not depend on the choice of $j \in \CC_J$. 
Equivalently, we can define $C_J$ as the union over all $j \in 
\CC_J$ of the curve segments $\{e^{i\theta}(\lambda_j(\theta) + 
i\lambda_j'(\theta))\,|\, \theta \in [0,\pi]\}$. Let
\[
  C_A \,=\, C_1 \cup \ldots \cup C_k \,\subset\, \C~, 
\]
and let $C_A' \subset \C$ be the {\em bitangent set} of 
$C_A$, namely the set of all line segments joining pairs of points 
of $C_A$ at which $C_A$ has the same tangent line. With these 
definitions, we have the following useful characterization of 
the singular set $\Sigma_A$: 

\begin{prop}\label{prop:singset} {\rm \cite[Theorem~3.5]{JS}} $~\Sigma_A = 
C_A \cup C_A'$. 
\end{prop}

\bepr
According to Lemma~\ref{lem:singset}, $\Sigma_A$ is the set of all 
complex numbers of the form $\langle Ax,x\rangle$ where $x \in \partial
\B^n$ is a normalized eigenvector of the Hermitian matrix $H(\theta)$
for some $\theta \in [0,\pi]$. To describe that set, fix $j \in 
\{1,\dots,m\}$, $\theta_0 \in \R$, and for $\theta$ in a 
neighborhood of $\theta_0$ let
\[
  z_j(\theta) \,=\, \langle Ax_j(\theta),x_j(\theta)\rangle~, 
\]
where $x_j(\theta)$ is a normalized eigenvector of $H(\theta)$ associated 
with the eigenvalue $\lambda_j(\theta)$ and depending smoothly on $\theta$. 
General results in perturbation theory imply that such an eigenvector
indeed  exists \cite{YB}. Using the definition \eqref{def:Htheta} of 
$H(\theta)$, we find
\begin{align*}
  \Re\,(e^{-i\theta}z_j(\theta)) \,&=\, \Re\,\langle e^{-i\theta} Ax_j(\theta),
  x_j(\theta)\rangle \,=\, \langle H(\theta) x_j(\theta),x_j(\theta)\rangle
  \,=\, \lambda_j(\theta)~, \\ 
  \Im\,(e^{-i\theta}z_j(\theta)) \,&=\, \Im\,\langle e^{-i\theta} Ax_j(\theta),
  x_j(\theta)\rangle \,=\, \langle H'(\theta) x_j(\theta),x_j(\theta)\rangle
  \,=\, \lambda_j'(\theta)~,
\end{align*}
since $\langle H(\theta) x_j'(\theta),x_j(\theta)\rangle + \langle 
H(\theta) x_j(\theta),x_j'(\theta)\rangle = 2\lambda_j(\theta)\Re\,
\langle x_j(\theta),x_j'(\theta)\rangle = 0$ due to the normalization 
condition. Thus
\begin{equation}\label{eq:zj}
  z_j(\theta) \,=\, e^{i\theta}(\lambda_j(\theta) + i\lambda_j'(\theta))~,
  \qquad \hbox{hence} \qquad 
  z_j'(\theta) \,=\, i e^{i\theta}(\lambda_j(\theta) + \lambda_j''(\theta))~.
\end{equation}
These relations show that the curve $C_J = \{z_j(\theta)\,|\,\theta \in 
\R\}$ is tangent, at each point $z_j(\theta)$, to the straight line
\begin{equation}\label{def:Lj}
  L_j(\theta) \,=\, \{e^{i\theta}(\lambda_j(\theta) + i\alpha) \,|\, 
  \alpha \in \R\} \,=\, \{z \in \C \,|\, \Re\,(z e^{-i\theta}) = 
  \lambda_j(\theta)\}~.
\end{equation}
In other words $C_J$ is the {\em envelope} of the family of straight 
lines $L_j(\theta)$, for any $j \in \CC_J$. Since $C_J \subset \Sigma_A$ 
by construction, we have shown that $\Sigma_A$ contains the curve 
$C_A = C_1 \cup \ldots \cup C_k$. 

However, it is important to realize that $\Sigma_A$ can be larger than
$C_A$ if the crossing set $\Theta \subset [0,\pi)$ defined above is
nonempty. Indeed, assume that $\lambda_j(\theta_0) = \lambda_p(\theta_0)$ 
for some $\theta_0 \in \Theta$ and some $j,p \in \{1,\dots,m\}$ with
$j \neq p$. For $\theta$ in a neighborhood of $\theta_0$, let
$x_j(\theta), x_p(\theta)$ be smooth, normalized eigenvectors of
$H(\theta)$ corresponding to $\lambda_j(\theta), \lambda_p(\theta)$
respectively. Using the same notations as above, we have
$z_j(\theta_0) \neq z_p(\theta_0)$ in general, because
$\lambda_j'(\theta_0) \neq \lambda_p'(\theta_0)$.  Now, since $\langle
x_j(\theta),x_p(\theta)\rangle = 0$ whenever $\lambda_j(\theta) \neq
\lambda_p(\theta)$, we also have $\langle
x_j(\theta_0),x_p(\theta_0)\rangle = 0$ by continuity.  In particular,
if $\alpha,\beta \in \C$ satisfy $|\alpha|^2 + |\beta|^2 = 1$, then $x
= \alpha x_j(\theta_0) + \beta x_p(\theta_0)$ is a normalized
eigenvector of $H(\theta_0)$, and a direct calculation yields
\[
  \Im\,\langle e^{-i\theta_0} A x,x\rangle \,=\, \langle H'(\theta_0) 
  x,x \rangle \,=\, |\alpha|^2 \lambda_j'(\theta_0) + 
  |\beta|^2 \lambda_p'(\theta_0)~,
\]
whereas $\Re\,\langle e^{-i\theta_0} A x,x\rangle = \langle H(\theta_0) 
x,x \rangle = \lambda_j(\theta_0) = \lambda_p(\theta_0)$. This 
shows that $\Sigma_A$ contains the line segment $[z_j(\theta_0),
z_p(\theta_0)]$, which by construction is tangent to the curve 
$C_J$ at $z_j(\theta_0)$ and to the curve $C_P$ at $z_p(\theta_0)$.
Repeating the same argument for all eigenvalue crossings, we conclude
that $\Sigma_A$ contains the whole bitangent set $C_A'$. Finally, it 
is clear from Lemma~\ref{lem:singset} that all points of $\Sigma_A$ 
either belong to $C_A$ or to $C_A'$.
\enpr

\bigskip\noindent{\bf Examples:}\\[1 mm] 
{\bf 1.} ({\em The generic case})
As is well known, within the space of all Hermitian matrices
of size $n \ge 2$, the set of matrices having a multiple eigenvalue
is a finite union of submanifolds of codimension at least three. 
This implies that, for a generic matrix $A \in \bM_n(\C)$, the 
Hermitian matrices $H(\theta)$ defined by \eqref{def:Htheta} will
have simple eigenvalues for all $\theta \in S^1$ \cite{JAG}. In 
that situation, we denote by $\lambda_1(\theta) < \lambda_2(\theta)
< \ldots < \lambda_n(\theta)$ the eigenvalues of $H(\theta)$. Using
the notations introduced above, we have $m = n$ and the permutation
$\tau$ defined in \eqref{def:tau} is simply
\[
  \tau \,=\, \begin{pmatrix} 1 & 2 & \dots & n \\
  n & n-1 & \dots & 1\end{pmatrix}~.
\]
If $n = 2k$ is even, $\tau$ has $k$ cycles with $\ell_1 = \dots \ell_k 
= 2$; if $n = 2k-1$ is odd, $\tau$ has one fixed point and $k-1$ 
cycles of length $2$. In all cases, the multiplicities $m_1,\dots,
m_k$ are all equal to $1$, and the bitangent set $C_A'$ is empty 
by assumption. Thus $\Sigma_A = C_A$ is the union of $k = [\frac{n+1}{2}]$
closed curves. It is not difficult to prove that the curve $C_1$ 
associated with the cycle $\CC_1 = (1~n)$ is smooth, strictly convex, 
and contains all the other curves $C_2,\dots,C_k$ in its interior
\cite{Kip}. In particular, $C_1 = \partial W(A)$. Consider now the
curve $C_J = \{z_j(\theta)\,|\,\theta \in \R\}$ for some $j \neq 
1,n$. If $\delta_j(\theta) = \lambda_j(\theta) + \lambda_j''(\theta)$
is not identically zero, the formulas \eqref{eq:zj} show that the
curvature of $C_J$ at any regular point $z_j(\theta)$ is 
strictly positive:
\[
  \kappa_j(\theta) \,=\, \frac{1}{|z_j'(\theta)|} \,=\, 
  \frac{1}{|\lambda_j(\theta) + \lambda_j''(\theta)|} \,>\, 0~.
\]
This means that the tangent vector $z_j'(\theta)$ always rotates 
counterclockwise when $\theta$ is increased. Nevertheless, the whole 
curve $C_J$ is not convex in general, because it may have a finite 
number of singular points corresponding to zeros of $\delta_j(\theta)$.
As is easily verified, simple zeros of $\delta_j(\theta)$ correspond
to {\em cusp points} of the curve $C_J$, see e.g. Fig.\,\ref{Fig1} where
a generic example with $n = 3$ is represented.  On the other hand, 
if $\delta_j(\theta)$ vanishes identically, then $z_j'(\theta) 
\equiv 0$ and the curve $C_J$ reduces to a single point. Under our 
generic assumptions, this can happen only if $n$ is odd and $j = 
(n+1)/2$. As an example, the singular set of the matrix $A_3$ defined 
in \eqref{eq:defA3} includes the isolated point $\{0\}$. 

\smallskip\noindent 
{\bf 2.} ({\em The normal case}) In contrast with the previous example, 
we now consider the particular case where the matrix $A$ is normal. 
If $\lambda_1,\dots,\lambda_m$ denote the distinct eigenvalues of 
$A$, it is straightforward to verify that the eigenvalues of the Hermitian 
matrices $H(\theta)$ defined in \eqref{def:Htheta} are simply 
$\lambda_j(\theta) = \Re\,(\lambda_j\,e^{-i\theta})$, $j = 1,\dots,m$. 
In view of \eqref{eq:zj}, this means that $z_j(\theta) = \lambda_j$ 
for all $\theta \in \R$, hence the curve $C_j$ is reduced to the 
single point $\{\lambda_j\}$ for all $j = 1,\dots,m$. Needless to 
say, the permutation $\tau$ defined by \eqref{def:tau} is the identity 
here. It follows that $C_A = \sigma(A) = \{\lambda_1,\dots,\lambda_m\}$, 
and proceeding as in the proof of Proposition~\ref{prop:singset} we 
easily see that $C_A'$ is the set of all line segments joigning
pairs of eigenvalues of $A$. This is in full agreement with the 
conclusions of Section~\ref{ss:normal}. 

\bigskip\noindent{\bf Remark:} We have seen that the curve $C_A$ which
generates the singular set $\Sigma_A$ is the envelope of the family of
straight lines $L_j(\theta)$ defined in \eqref{def:Lj}. This geometric
construction can be formulated in an equivalent way \cite{At}, which
is more conceptual and worth mentioning here. Assume for simplicity
that the matrix $A \in \bM_n(\C)$ is ``generic'' in the sense of
Example~1 above, and let $a(\xi)$ be the homogeneous polynomial of
degree $n$ defined by
\[
  a(\xi) = \det\Bigl(\xi_0 I_n + \xi_1 A_1 + \xi_2 A_2\Bigr)~,
  \qquad \xi = (\xi_0,\xi_1,\xi_2) \in \R^3~,
\]
where $A_1, A_2$ are as in \eqref{def:Htheta}. Since $a(\xi)$ has real
coefficients, the equation $a(\xi) = 0$ defines, in homogeneous
coordinates, an algebraic curve $\Gamma$ in the projective plane $\R
P^2$. Moreover, our genericity assumption on $A$ implies that this
curve is nonsingular: for each $\bar\xi \in \Gamma$, the tangent line
to $\Gamma$ at $\bar \xi$ is uniquely defined and satisfies, in
homogeneous coordinates, an equation of the form $x_0\xi_0 + x_1\xi_1
+ x_2\xi_2 = 0$. The set of all $x = (x_0,x_1,x_2)$ obtained in this
way is again a curve $\Gamma'$ in $\R P^2$ (the {\em dual} curve of
$\Gamma$) given by the equation $b(x) = 0$ for some homogeneous
polynomial $b$. As is shown in \cite{Fie}, the degree of $b$ does not
exceed $n(n-1)$ if $n \ge 2$. Now, it is rather straightforward to
verify that the curve $C_A$ defined as the envelope of the family of
straight lines \eqref{def:Lj} is nothing but the restriction of the
projective curve $\Gamma'$ to the subspace $x_0 = 1$, namely $z = x+iy
\in C_A$ if and only if $b(1,x,y) = 0$. Thus $C_A$ is a real algebraic
curve in $\C \simeq \R^2$ of degree at most $n(n-1)$ if $n \ge 2$. In
the language of partial differential equations, the algebraic variety
$\Gamma$ is the {\em characteristic variety} of the symmetric
hyperbolic system \eqref{eq:hypsys}, and we shall see in
Section~\ref{s:main} that the dual variety $\Gamma'$ is related to the
singular support of the fundamental solution of \eqref{eq:hypsys}.
Note that our genericity assumption on $A$ precisely means that system
\eqref{eq:hypsys} is {\em strictly hyperbolic}.

\bigskip
We conclude this section with a brief discussion of the number 
of tangent lines to the algebraic curve $C_A$ which can be drawn
from a given point. We recall that $C_A = C_1 \cup \ldots \cup 
C_k$, where each $C_J$ is a closed curve associated with the 
cycle $\CC_J$ of the permutation \eqref{def:tau}. For all 
$J \in \{1,\dots,k\}$ and all $z \in \C\setminus C_J$, 
we denote by $N_J(z)$ the number of straight lines that are 
tangent to the curve $C_J$ and contain the point $z$. Note 
that, since $C_J$ was itself defined as the envelope of a 
family of straight lines, the tangent line to $C_J$ is well defined
even at singular points. In the degenerate case where $C_J$ 
reduces to a single point $\{z_J\}$, the set of tangent lines should 
be understood as the pencil of all straight lines through $z_J$. 
Now, if $z \in \C\setminus C_A$, we denote by $N(z)$ the total 
number of tangents to the curve $C_A = C_1 \cup \ldots \cup 
C_k$ that can be drawn from the point $z$, with multiplicities 
taken into account:
\begin{equation}\label{def:N}
  N(z) \,=\, m_1 N_1(z) + \dots + m_k N_k(z)~, \qquad 
  z \in \C\setminus C_A~.
\end{equation}
The following elementary properties of $N_J(z)$ and $N(z)$ will 
be useful. 

\begin{prop}\label{prop:tangents}
For each $J \in \{1,\dots,k\}$, the number $N_J(z)$ is 
constant in each connected component of $\C\setminus C_J$. 
Moreover $N(z) \le n$ for all $z \in \C \setminus C_A$. 
\end{prop}

\bepr
Fix $J \in \{1,\dots,k\}$ and pick $j \in \CC_J$. For any 
$z \in \C\setminus C_J$, $N_J(z)$ is the number of zeros 
of the function
\[
  f_j(\theta,z) \,=\, \lambda_j(\theta) - \Re\,(z e^{-i\theta})
\]
for $\theta$ in the interval $[0,\ell_J\pi)$. When $z$ is varied, 
this number can only change if $f_j(\theta,z)$ has a double 
zero for some $\theta$, but this would mean that $z = z_j(\theta) 
= e^{i\theta}(\lambda_j(\theta)+i\lambda_j'(\theta)) \in C_J$, 
thus contradicting our assumption. Therefore $N_J(z)$ is necessarily 
constant in each connected component of $\C\setminus C_J$. 
On the other hand, we have the identity
\begin{equation}\label{eq:detfac}
  \det\Bigl(H(\theta) - \Re\,(z e^{-i\theta})I_n\Bigr)
  \,=\, \prod_{J=1}^k \prod_{j \in \CC_J} f_j(\theta,z)^{m_J}~.
\end{equation}
Fix $z \in \C\setminus C_A$, and consider both sides of
\eqref{eq:detfac} as functions of $\theta$. If multiplicities are
taken into account, the number of zeros of the right-hand side on the
interval $[0,\pi)$ is precisely $N(z)$. But the left-hand side, being
a trigonometric polynomial of degree at most $n$, cannot have more
than $n$ zeros in $[0,\pi)$. This proves the claim.  
\enpr

\begin{rque}\label{rem:prod}
As was mentioned in the proof of Proposition~\ref{prop:singset}, 
if $z \in \C\setminus C_A$ it is possible that $f_j(\theta,z) = 
f_p(\theta,z) = 0$ for some $\theta \in [0,\pi)$ and some $j \neq p$. 
This is the case, in particular, whenever $z \in C_A'$. However, 
if $z \in \C\setminus \Sigma_A$, we have
\[
  \partial_\theta f_j(\theta,z) \cdot \partial_\theta f_p(\theta,z)
  \,>\, 0 \qquad \hbox{whenever}\quad f_j(\theta,z) = 
  f_p(\theta,z) = 0~.
\]
Indeed, replacing $A$ with $A - zI_n$, we can assume without loss of
generality that $z = 0$. If $\lambda_j(\theta) = \lambda_p (\theta) =
0$ for some $\theta \in [0,\pi)$ and some $j \neq p$, then
$\lambda_j'(\theta)$ and $\lambda_p'(\theta)$ have necessarily the
same sign, otherwise \eqref{eq:zj} would imply that the origin belongs
to the line segment $[z_j(\theta),z_p(\theta)]$, thus contradicting
our assumption that $z \notin C_A'$.
\end{rque}

\begin{figure}
\begin{tabular}{cc}
\includegraphics[width=8.0cm,height=5cm]{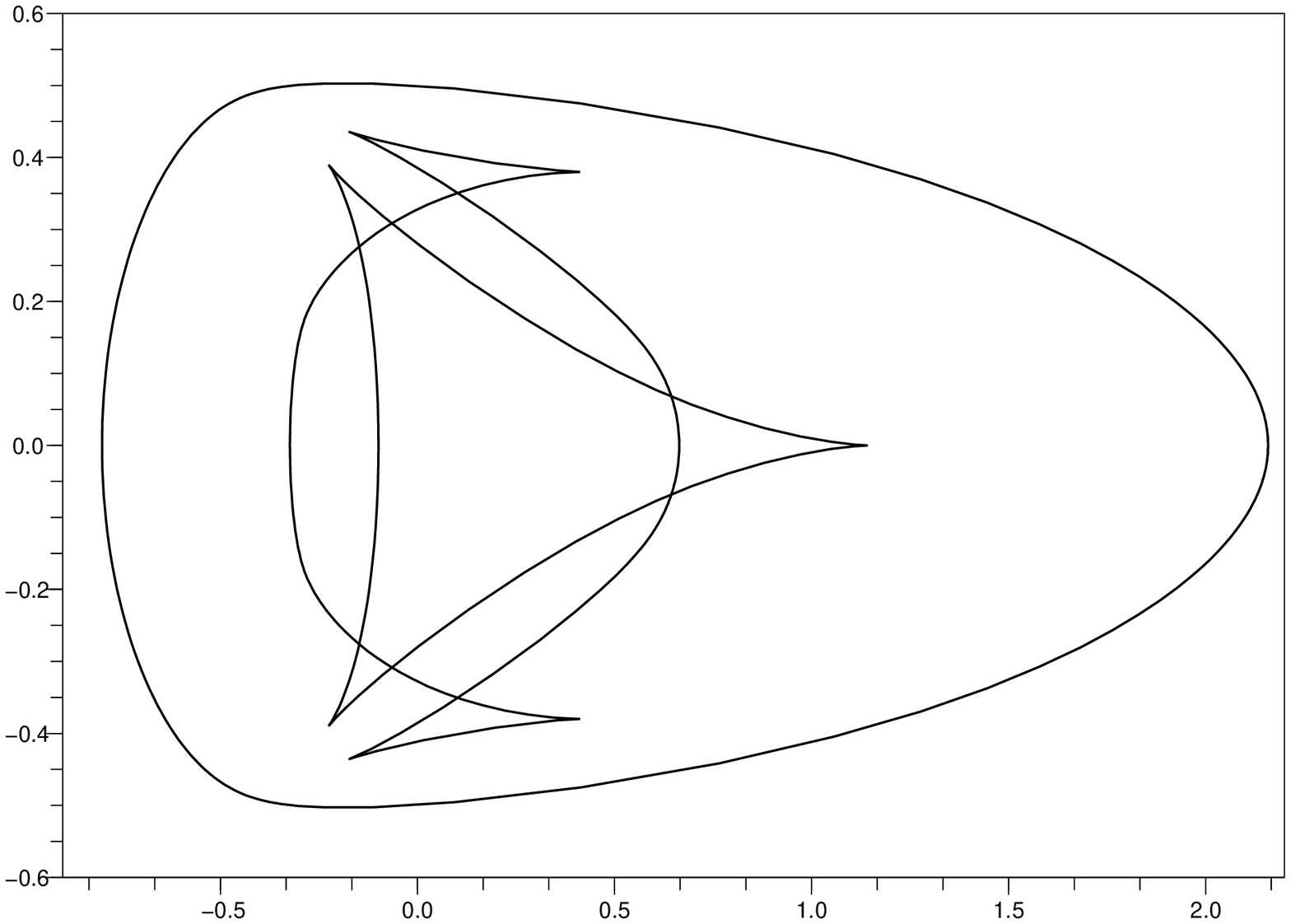} & 
\includegraphics[width=8.0cm,height=5cm]{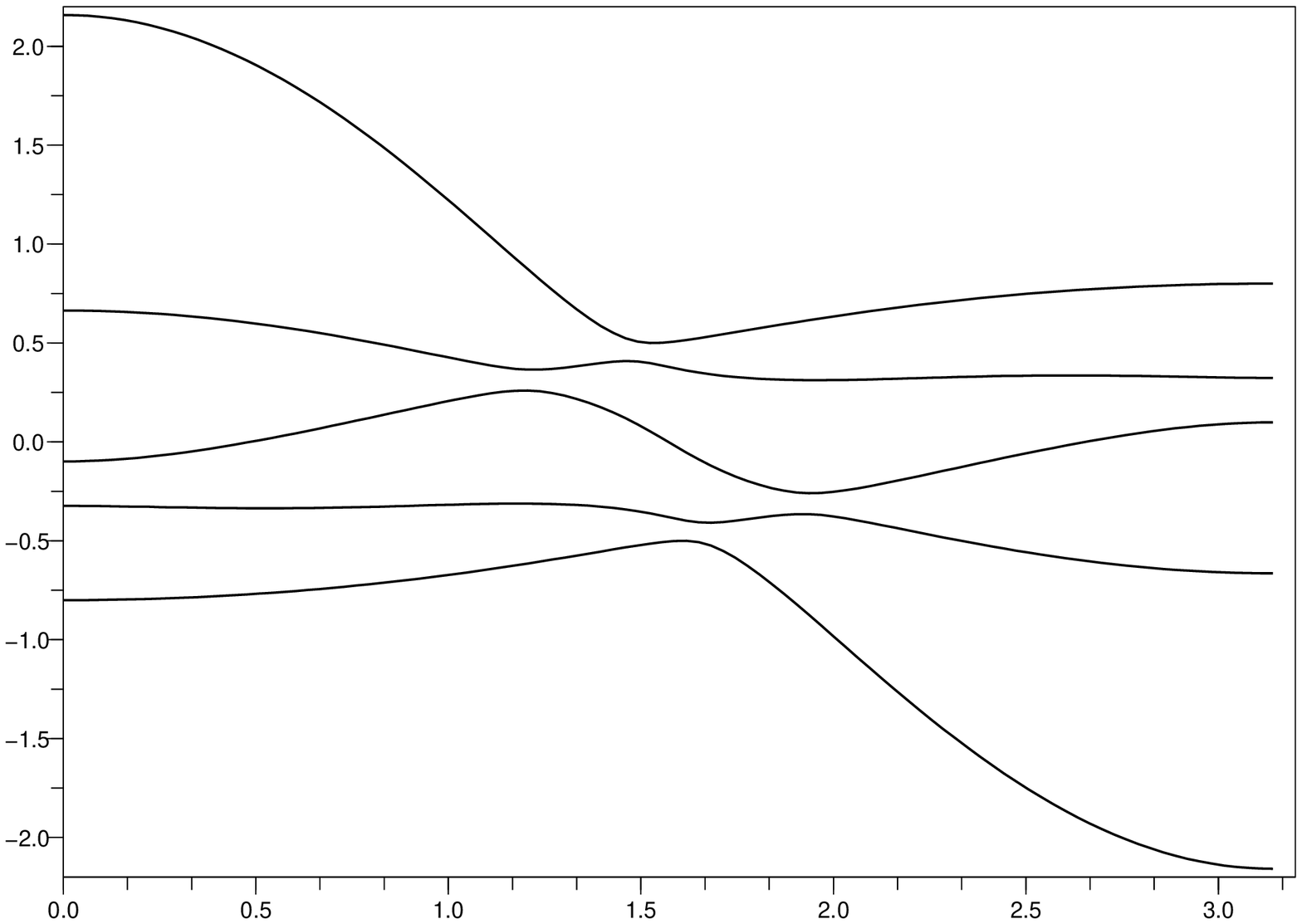} \\
\end{tabular}
\caption{\small The singular set $\Sigma_A$ (left) and the 
eigenvalues $\lambda_1(\theta),\dots,\lambda_5(\theta)$ (right) 
are represented for a typical matrix $A \in M_5(\R)$. The singular 
set consists of three closed curves, one of which (the boundary of $W(A)$) 
is smooth, and the other two have cusps. In this example, the set 
$\Pi_A$ defined in \eqref{eq:Pidef} has six connected components.}
\label{Fig2}
\end{figure}

Proposition~\ref{prop:tangents} asserts that the integer $N_J(z)$ can
only change if $z$ crosses the curve $C_J$. In fact, if the crossing
occurs at a regular point $\bar z \in C_J$, it is not difficult to
verify that $N_J(z)$ is decreased by two units if $z$ crosses $C_J$
from the convex to the concave side (i.e., in the direction of the
local center of curvature), and increased by two units if $z$ crosses
$C_J$ in the converse direction, see \cite[Section~4.1]{BY} or
Section~\ref{ss:regular} below. These simple rules give an efficient
algorithmic way to compute $N_J(z)$, and hence $N(z)$, in concrete
examples. Consider for instance Fig.\,\ref{Fig2}, where the 
singular set $\Sigma_A$ of a generic matrix $A \in M_5(\R)$ is 
represented. Here $\Sigma_A = C_A = C_1 \cup C_2 \cup C_3$, 
where $C_1$ is the boundary of $W(A)$, $C_2$ is the closed curve with 
two swallowtails, and $C_3$ is the triangle with three cusps. The
set $\C\setminus\Sigma_A$ has $11$ connected components, on which 
$N(z)$ is equal to $5$, $3$, or $1$. Using the crossing rules 
above, it is easy to verify that $N(z) = 5$ in six different regions: 
inside both swallowtails, inside the three tips of the triangle, and
outside $W(A)$. Such a result is definitely more cumbersome to obtain 
by counting directly the number of tangents to $C_A$ from a given 
point. 

\section{Qualitative properties of the numerical density}\label{s:main}

Equipped with the results of Sections~\ref{s:Radon} and \ref{s:geom},
we now derive some of the main properties of the numerical density
$f_A$ of an arbitrary matrix $A \in \bM_n(\C)$. We first establish an
explicit formula for the derivatives of order $n-2$, which allows us
to prove that $f_A$ is polynomial in some distinguished regions of the
complex plane which can be characterized geometrically. We next show
that, for a generic matrix $A \in \bM_n(\C)$, the density $f_A$ is of
class $C^{n-3}$ if $n \ge 3$. Finally, as announced in the introduction, 
we prove that the fundamental solution of the linear hyperbolic
system~\eqref{eq:hypsys} can be represented in terms of derivatives of
the numerical density $f_A$.  In particular, the lacunas of
system~\eqref{eq:hypsys} are precisely the polynomial regions
described in Corollary~\ref{cor:poly}.

\subsection{Polynomial regions}\label{ss:poly}

Most of what we know about the numerical density of nonnormal 
matrices is based on the following result:

\begin{prop}\label{prop:deriv}
Let $A \in \bM_n(\C)$, where $n \ge 2$, and let $\PP = 
P(\partial_x,\partial_y)$ be a homogeneous differential operator 
of degree $n-2$. Then, for all $z = x+iy \in \C\setminus \Sigma_A$, 
we have
\begin{equation}\label{eq:PfA}
  (\PP f_A)(z) \,=\, -\frac{(n{-}1)!}{4\pi^2} \,\fp \!\int_{S^1} 
  \frac{P(\cos\theta,\sin\theta)}{\Delta(\theta,z)}\dd \theta~,
\end{equation}
where
\begin{equation}\label{def:Delta}
  \Delta(\theta,z) \,=\, \det\Bigl(H(\theta) - \Re\,(z e^{-i\theta})
  I_n\Bigr)~, \quad \theta \in S^1~.
\end{equation}
\end{prop}

In \eqref{eq:PfA} the symbol $\fp$ denotes the finite part of 
the integral in the sense of Hadamard, but in many situations 
it is sufficient to take simply the Cauchy principal value, as 
in \eqref{eq:twod1}. This is the case in particular for generic 
matrices in the sense of Section~\ref{s:geom} (Example~1), because
all multiplicities $m_1, \dots, m_k$ are all equal to one and 
formula \eqref{eq:detfac} then shows that the map $\theta \mapsto 
\Delta(\theta,z)$ has only simple zeros on $S^1$ if $z \in 
\C\setminus\Sigma_A$. 

In fact, using the analyticity of the integrand, it is possible to
rewrite \eqref{eq:PfA} in a slightly different form which is
appropriate for further analysis. As in Section~\ref{ss:twod}, we set
$w = e^{2i\theta}$ and we observe that
\[
  P(\cos\theta,\sin\theta) \,=\, w^{1-\frac{n}{2}} \tilde P(w)~, 
  \qquad \Delta(\theta,z) \,=\, w^{-\frac{n}{2}} \tilde \Delta(w,z)~,
\]
where (if $z = x+iy$)
\[
  \tilde P(w) \,=\, P\Bigl(\frac{w+1}{2},\frac{w-1}{2i}\Bigr)~,
  \qquad
  \tilde\Delta(w,z) \,=\, \det\left(\frac{A}{2}+\frac{wA^*}{2}
  -\Bigl(x\frac{w+1}{2} + y\frac{w-1}{2i}\Bigr)I_n\right)~.
\]
Thus \eqref{eq:PfA} is equivalent to
\begin{align}\nonumber
  (\PP f_A)(z) \,&=\, \frac{i(n{-}1)!}{4\pi^2} \,\fp \!\oint_{|w|=1} 
  \frac{\tilde P(w)}{\tilde\Delta(w,z)}\dd w \\ \label{eq:PfA2}
  \,&\equiv\, \frac{i(n{-}1)!}{8\pi^2} \oint_{|w|=1-\epsilon} 
  \frac{\tilde P(w)}{\tilde\Delta(w,z)}\dd w + \frac{i(n{-}1)!}{8\pi^2} 
  \oint_{|w|=1+\epsilon} \frac{\tilde P(w)}{\tilde\Delta(w,z)}\dd w~,
\end{align}
where $\epsilon$ is any sufficiently small positive number, depending 
on $z$. In the particular case where $P = 1$ and $A \in \bM_2(\C)$ is 
given by \eqref{def:A2}, we recover \eqref{eq:twod2}. 

\bigskip\noindent{\bf Proof of Proposition~\ref{prop:deriv}.} 
We first consider the generic situation where the Hermitian matrices 
$H(\theta)$ have simple eigenvalues $\lambda_1(\theta) < 
\lambda_2(\theta) < \dots < \lambda_n(\theta)$ for all $\theta \in S^1$. 
In that case, the $B$-spline representing the numerical density of
$\mu_H(\theta)$ can be given an explicit expression using the 
divided difference formula \eqref{def:dd1}, \eqref{def:Bspline}:
\[
  B[\lambda_1(\theta),\dots,\lambda_n(\theta)](s) \,=\, 
  (n-1)\sum_{j=1}^n \frac{(\lambda_j(\theta)-s)_+^{n-2}}{
  \prod_{k\neq j} (\lambda_j(\theta)-\lambda_k(\theta))}~.
\]
In particular, differentiating $(n-1)$ times with respect to $s$, 
we find
\begin{equation}\label{eq:Bn1}
  B^{(n-1)}[\lambda_1(\theta),\dots,\lambda_n(\theta)](s) \,=\, 
  (n-1)! \,\sum_{j=1}^n \frac{\delta(s-\lambda_j(\theta))}{
  \prod_{k\neq j} (\lambda_k(\theta)-\lambda_j(\theta))}~.
\end{equation}
Let $\PP = P(\partial_x,\partial_y)$ be a homogeneous differential
operator of degree $n-2$. From the representation formula 
\eqref{eq:Radon3}, we deduce at least formally
\begin{equation}\label{eq:PfA3}
  (\PP f_A)(x+iy) \,=\, \frac{1}{4\pi} \int_{S^1} P(\cos \theta,
  \sin \theta)\,\HH B^{(n-1)}[\lambda_1(\theta),\dots,
  \lambda_n(\theta)](x\cos\theta + y\sin\theta)\dd\theta~.
\end{equation}
To evaluate the integrand in the right-hand side, we start from 
\eqref{eq:Bn1} and recall that the Hilbert transform (with respect to 
the variable $s$) of the Dirac measure $\delta(s-\lambda)$ is 
the distribution $\frac{1}{\pi}\,\pv \frac{1}{s-\lambda}$. We also
use the identity
\[
  \sum_{j=1}^n \frac{1}{\mu_j \prod_{k\neq j} (\mu_k-\mu_j)}
  \,=\, \prod_{j=1}^n \frac{1}{\mu_j}~,
\]
which holds for any collection of pairwise distinct nonzero 
complex numbers $\mu_1,\dots,\mu_n$. We thus find
\[
  \HH B^{(n-1)}[\lambda_1(\theta),\dots,\lambda_n(\theta)](s)
  \,=\, -\frac{(n{-}1)!}{\pi}\,\pv \prod_{j=1}^n 
  \frac{1}{\lambda_j(\theta)-s} 
  \,=\, -\frac{(n{-}1)!}{\pi}\,\pv \frac{1}{\det(H(\theta)-sI_n)}~. 
\]
Setting $s = x\cos\theta + y\sin\theta$ and inserting this expression
into \eqref{eq:PfA3}, we obtain \eqref{eq:PfA}. 

The calculations so far are formal, but they can be justified if 
we assume that $z = x+iy \in \C\setminus\Sigma_A$. In that case, 
we know from the proof of Proposition~\ref{prop:tangents} that 
the map $\theta \mapsto \Delta(\theta,z)$ has only simple zeros 
on $S^1$, because this is the case for each of the factors 
$f_j(\theta,z)$ in \eqref{eq:detfac} and, by assumption, the 
eigenvalues $\lambda_j(\theta)$ are all distinct for $\theta \in 
S^1$. It follows that the integral in \eqref{eq:PfA} is well-defined
in the sense of Cauchy's principal value, and depends smoothly 
on $z \in \C\setminus\Sigma_A$. This in turn implies that the 
density $f_A$ given by \eqref{eq:Radon3} is smooth on $\C\setminus
\Sigma_A$, as it should be, and that the calculations above are 
correct. 

To conclude the proof of Proposition~\ref{prop:deriv}, it remains
to verify that \eqref{eq:PfA} or \eqref{eq:PfA2} holds for an
arbitrary matrix $A \in \bM_n(\C)$. To do that, we first observe that 
the singular set $\Sigma_A \subset \C$ is an upper-semicontinuous
function of $A$ in the sense that
\[
  \delta(\Sigma_B,\Sigma_A) \,\equiv\, \sup_{z \in \Sigma_B}
  \dist(z,\Sigma_A) \,\xrightarrow[B\to A]{}\, 0~.
\]
Moreover the numerical density $f_A$, together with its derivatives,
depends continuously on $A$ in $\C\setminus\Sigma_A$. These rather
classical facts can be established using, for instance, the
representation formula \eqref{def:fA} (we omit the details). On the
other hand, it is not difficult to verify that the right-hand side of
\eqref{eq:PfA2} depends continuously on $A$ for each $z \in
\C\setminus\Sigma_A$. The crucial point here is that the polynomial
map $w \mapsto \tilde \Delta(w,z)$ keeps the same number of zeros on
the unit circle (counted with multiplicities) if the matrix $A$ is
slightly varied; in particular, we can choose the same $\epsilon > 0$
in \eqref{eq:PfA2} for all matrices in a neighborhood of $A$. This
property can be established using the factorization \eqref{eq:detfac},
Remark~\ref{rem:prod}, and general results for perturbations of
eigenvalues of Hermitian matrices, see \cite{Kat}.  Now, given $A \in
\bM_n(\C)$, there exists a sequence $\{A_\ell\}_{\ell\in\N}$ of
generic matrices converging to $A$ as $\ell \to \infty$. If $z \in
\C\setminus\Sigma_A$, we know that equation \eqref{eq:PfA2} holds for
$A_\ell$ if $\ell$ is sufficiently large, hence taking the limit $\ell
\to \infty$ and using the continuity properties mentioned above we
obtain the desired equality.  \enpr

\medskip
As a consequence of Proposition~\ref{prop:deriv}, we now establish
an important property of the numerical density in the regions of 
the complex plane where the number $N(z)$ defined in \eqref{def:N} 
takes its maximal value $n$. 

\begin{cor}\label{cor:poly} Given $A \in \bM_n(\C)$, let 
\begin{equation}\label{eq:Pidef}
  \Pi_A \,=\, \Bigl\{z \in \C \setminus \Sigma_A \,\Big|\, 
  N(z) = n\Bigr\}~.
\end{equation}
If $n \ge 3$, the numerical density $f_A$ is polynomial of degree 
at most $n-3$ on each connected component of $\Pi_A$. If $n = 2$, 
then $f_A = 0$ on $\Pi_A$. 
\end{cor}

\bepr 
We shall prove that $\PP f_A$ vanishes identically on $\Pi_A$
for any homogeneous differential operator of order $n-2$. Indeed, in
view of \eqref{eq:detfac}, \eqref{def:Delta}, the assumption $z \in
\Pi_A$ implies that the map $\theta \mapsto \Delta(\theta,z)$ has
exactly $n$ zeros (counting multiplicities) on $[0,\pi)$. Equivalently,
the polynomial $\tilde\Delta(w,z)$ has exactly $n$ zeros on the unit
circle $\{|w| = 1\}$. But since $\tilde\Delta(w,z)$ has degree $n$,
this polynomial has no zeros outside the unit circle, and using 
Cauchy's theorem we conclude that the first integral in the last 
member of \eqref{eq:PfA2} vanishes. The second integral is also zero, 
because the numerator is a polynomial of degree at most $n-2$, while the
denominator has degree exactly $n$, hence the integrand $\tilde
P(w)/\tilde \Delta(w,z)$ decays at least like $|w|^{-2}$ as $|w| \to
\infty$. Thus $\PP f_A \equiv 0$ on $\Pi_A$, and the conclusion 
follows. 
\enpr

\begin{rque} If $z = x+iy \in \sigma(A)$, then the polynomial 
$\tilde \Delta(w,z)$ has degree strictly less than $n$, and 
it follows from the above proof that $z \notin \Pi_A$. Thus 
$\sigma(A) \cap \Pi_A = \emptyset$. 
\end{rque}

The set $\Pi_A$ is never empty, because it always contains the
complement of the numerical range $W(A)$, where the density $f_A$
vanishes identically \cite{YB,BY}. Moreover, in many situations, one
or several components of $\Pi_A$ are contained in $W(A)$, in which
case Corollary~\ref{cor:poly} gives nontrivial informations on the
numerical density. For instance, if $A \in \bM_n(\C)$ is a normal
matrix whose numerical range has nonempty interior, then $N(z) = n$
for all $z \in \C\setminus \sigma(A)$, and it follows from
Corollary~\ref{cor:poly} that $f_A$ is polynomial of degree at most
$n-3$ in each connected component of $\C\setminus\Sigma_A$, in
agreement with Proposition~\ref{prop:Bspline2}.  In the same spirit,
if $A = A_1 \oplus A_2$ with $A_1 \in \bM_{n_1}(\C)$ and $A_2\in
\bM_{n_2}(\C)$, it is easily verified that $N(z) = n = n_1+n_2$ if 
$z \notin W(A_1) \cup W(A_2)$, hence $f_A$ is piecewise polynomial
outside $W(A_1) \cup W(A_2)$. Finally, Fig.\,\ref{Fig1} shows a
typical example of a matrix $A \in \bM_3(\C)$ for which $\Pi_A$ has a
component inside $W(A)$, on which the density $f_A$ is identically
constant by Corollary~\ref{cor:poly}.

\subsection{Generic regularity results}\label{ss:regular}

Our purpose here is to establish regularity results for the 
numerical density $f_A$ in the whole complex plane, and not
only outside the singular set $\Sigma_A$. For simplicity, we
assume henceforth that our matrix $A \in \bM_n(\C)$ enjoys 
the following (generic) properties:

\medskip\noindent{\bf H1:} The eigenvalues of \eqref{def:Htheta} 
satisfy $\lambda_1(\theta) < \ldots < \lambda_n(\theta)$ for all $\theta 
\in S^1$. 

\medskip\noindent{\bf H2:} For all $j \in \{1,\dots,n\}$, the function
$\lambda_j(\theta) + \lambda_j''(\theta)$ is not identically zero. 

\medskip
The second assumption guarantees that the curve \eqref{def:CJ} associated 
with $\lambda_j$ is not reduced to a single point. This of course is 
possible only if $n \ge 2$, and the two-dimensional case $n = 2$ 
is completely treated in Section~\ref{ss:twod}. So we can assume 
that $n \ge 3$, and we have the following result: 

\begin{prop}\label{prop:gen}
Assume that $A \in \bM_n(\C)$ satisfies H1, H2 above. If $n \ge 3$, 
the numerical density $f_A : \C \to \R_+$ is of class $C^{n-3}$. 
\end{prop}

We have already seen that both hypotheses H1, H2 are necessary, in
general, for the conclusion of Proposition~\ref{prop:gen} to hold.
For instance, if $A \in M_3(\C)$ is a normal matrix whose numerical
range $W(A)$ has nonempty interior, the numerical density $f_A$ is
proportional to the characteristic function of $W(A)$ and is therefore
discontinuous on $\partial W(A)$. A more subtle example is provided by
the matrix $A_3$ defined in \eqref{eq:defA3}: if $|a|^2 + |b|^2 = 4$,
we have here $\lambda_1(\theta) = -1$, $\lambda_2(\theta) = 0$,
$\lambda_3(\theta) = 1$ for all $\theta \in S^1$. Thus H1 is satisfied, 
but obviously not H2, and the explicit formula \eqref{eq:fA3} shows
that the numerical density of $A_3$ is discontinuous at the origin.

\bigskip\noindent{\bf Proof of Proposition~\ref{prop:gen}.} Assume 
that $n \ge 3$ and that $A \in \bM_n(\C)$ satisfies H1, H2. If 
$\QQ = Q(\partial_x,\partial_y)$ is a homogeneous differential 
operator of degree $n-3$, we have as in \eqref{eq:PfA3}:
\begin{equation}\label{eq:QfA1}
  (\QQ f_A)(x+iy) \,=\, \frac{1}{4\pi} \int_{S^1} Q(\cos \theta,
  \sin \theta)\,\HH B^{(n-2)}[\lambda_1(\theta),\dots,
  \lambda_n(\theta)](x\cos\theta + y\sin\theta)\dd\theta~,
\end{equation}
where
\begin{equation}\label{eq:Bn2}
  B^{(n-2)}[\lambda_1(\theta),\dots,\lambda_n(\theta)](s) \,=\, 
 -(n-1)! \,\sum_{j=1}^n \frac{\mathbf{H}(\lambda_j(\theta)-s)}{
  \prod_{k\neq j} (\lambda_k(\theta)-\lambda_j(\theta))}~.
\end{equation}
Here $\mathbf{H} : \R \to [0,1]$ denotes the Heaviside function. Taking
the Hilbert transform of \eqref{eq:Bn2} with respect to $s$ and 
using \eqref{eq:QfA1}, we thus find
\begin{align}\nonumber
  (\QQ f_A)(z) \,&=\, \frac{(n{-}1)!}{4\pi^2}\sum_{j=1}^n \int_{S^1} 
  \frac{Q(\cos \theta,\sin \theta)}{\prod_{k\neq j} (\lambda_k(\theta)
  -\lambda_j(\theta))}\,\log|f_j(\theta,z)|\dd\theta \\ \label{eq:QfA2}
  \,&\equiv\, \sum_{j=1}^n \int_{S^1} h_j(\theta)\log|f_j(\theta,z)|
  \dd\theta~, 
\end{align}
where $f_j(\theta,z) = \lambda_j(\theta) -\Re(z e^{-i\theta})$ and 
$h_j : S^1 \to \R$ is a smooth function. Note that the integrand 
in \eqref{eq:QfA2} is $2\pi$-periodic, because by assumption H1 
this is the case for all eigenvalues $\lambda_1(\theta),\dots,
\lambda_n(\theta)$.

It remains to show that each integral in the right-hand side of  
\eqref{eq:QfA2} defines a continuous function of $z \in \C$. 
Replacing $A$ by $A - zI_n$ (an operation which does not
affect the properties H1, H2), we see that it is sufficient
to prove continuity at $z = 0$. This in turn is obvious if
$\lambda_j(\theta) = f_j(\theta,0)$ does not vanish, so from now
on we focus on the case where $\lambda_j(\theta)$ has (isolated) 
zeros on $S^1$. Using a partition of unity, we can treat each 
zero separately, so it is sufficient to consider the case where
$\lambda_j(\theta)$ has a single zero of order $q \ge 1$ at 
$\theta = 0$, and $h_j(\theta)$ is localized near the origin. 
By analyticity, for $z$ close to zero we have the factorization
\[
  f_j(\theta,z) \,=\, g(\theta,z)\,\prod_{p=1}^q (\theta - \mu_p(z))~,
\]
where the (possibly complex) roots $\mu_p(z)$ depend continuously
on $z$, with $\mu_p(0) = 0$, and $g(\cdot,z)$ does not vanish 
in a neighborhood of zero. The quantity we have to study is
therefore
\[
  \sum_{p=1}^q \int_{\R} h_j(\theta) \log|\theta -\mu_p(z)| \dd\theta 
  + \int_{\R} h_j(\theta) \log|g(\theta,z)|\dd\theta~. 
\]
The last term is clearly a continuous function of $z$. In the integral
involving $\mu_p$, we make the change of variables $\theta = t + \Re
\mu_p(z)$ and observe that $\log|t| \le \log|t-i\Im \mu_p(z)| \le 0$
in a neighborhood of $(t,z) = (0,0)$, hence continuity with respect to
$z$ follows from Lebesgue's dominated convergence theorem.
\enpr

\medskip
Once continuity of the derivatives of order $n-3$ has been established, 
we can obtain further regularity results for the numerical density by 
using the representation formula \eqref{eq:PfA} or \eqref{eq:PfA2}.
In view of Lemma~\ref{lem:singsupp}, it is sufficient to study the 
density in a neighborhood of a point $\bar z \in \Sigma_A$. Since 
$\Sigma_A = C_A = C_1 \cup \ldots \cup C_k$ by assumption H1, there
exists $j \in \{1,\dots,n\}$ such that $\bar z \in C_j$, and for 
simplicity we also assume that $\bar z \notin C_p$ for all $p \neq j$.
This means that the function $\theta \mapsto f_j(\theta,\bar z)$ has a 
zero of order $m \ge 2$ at some point $\bar\theta \in S^1$, and only
simple zeros for $\theta \neq \bar\theta \mod \pi$; moreover, 
if $p \neq j$, $\theta \mapsto f_p(\theta,\bar z)$ has only simple zeros 
and does not vanish for $\theta = \bar \theta$. 

Without loss of generality, we assume from now on that $\bar z = 0$ 
and $\bar\theta = 0$, and we first consider the simplest case where $m=2$. 
This means that $\lambda_j(\theta) = \frac{\alpha}{2}\theta^2 + 
\OO(\theta^3)$ near $\theta = 0$, for some $\alpha \neq 0$. If 
$z = x+iy$ is sufficiently small, it follows that the analytic 
function $\theta \mapsto f_j(\theta) = \lambda_j(\theta) - (x\cos\theta
+y\sin\theta)$ has exactly two zeros $\theta_\pm(z)$ in a neighborhood
of the origin, which satisfy
\[
  \theta_\pm(z) \,=\, \frac{1}{\alpha}\Bigl(y \pm \sign(\alpha)
  \sqrt{y^2+2\alpha x}\Bigr) + \OO(|x| + |y|^2)~.
\]
The critical curve $C_j$, which is the set of all points $z$ for which
$\theta_+(z) = \theta_-(z)$, is therefore given by the equation $x = 
-\frac{1}{2\alpha}y^2 + \OO(y^3)$ in a small ball $B$ around the 
origin. Moreover $B \setminus C_j = B_r \cup B_c$, where $B_r$ is 
the set of all $z \in B$ for which the roots $\theta_\pm(z)$ are 
real and distinct, whereas $z \in B_c$ when $\theta_\pm(z)$ are 
complex conjugate with nonzero imaginary part. As is easily verified, 
the local center of curvature of $C_j$ is located on the side of $B_c$. 

Now, let $\PP = P(\partial_x,\partial_y)$ be a homogeneous differential
operator of degree $n-2$, and consider the expression of $(\PP f_A)(z)$ 
given by \eqref{eq:PfA2}. Assume that, in the right-hand side, the 
parameter $\epsilon > 0$ is chosen in such a way that the slit 
annulus $\AA_\epsilon = \{w \in \C\,|\, 0 < |1-|w|| < \epsilon\}$ 
contains the points $w_\pm(z) = e^{2i\theta_\pm(z)}$ for all $z \in 
B_c$, but that $\AA_{2\epsilon}$ does not contain any other root of 
the determinant $\tilde \Delta(w,z)$ for $z \in B$ (these conditions
are easily achieved by choosing first $\epsilon$ and then $B$ 
sufficiently small). Under these assumptions, the right-hand side 
of \eqref{eq:PfA2} defines a smooth function of $z \in B$, which 
coincides with $(\PP f_A)(z)$ if $z \in B_r$ but {\em not} if 
$z \in B_c$. Indeed, in the latter case, we have to consider in
addition the roots $w_\pm(z)$ of $\tilde \Delta(w,z)$ which are not 
taken into account by the fixed integration contours in \eqref{eq:PfA2}
since $|1-|w_\pm(z)|| < \epsilon$. Using Cauchy's theorem, we easily obtain
\[
  (\PP f_A)(z) = (\PP f_A)_\reg(z) + \frac{(n{-}1)!}{4\pi}
  \,\frac{\tilde Q(w_-(z),z) + \tilde Q(w_+(z),z)}{w_-(z) - w_+(z)}~,
  \quad z \in B_c~,
\]
where $(\PP f_A)_\reg(z)$ denotes the regular part of $(\PP f_A)(z)$, 
given by the right-hand side of \eqref{eq:PfA2} with fixed $\epsilon$, 
and
\[
  \tilde Q(w,z) \,=\, \frac{\tilde P(w) (w-w_+(z))(w-w_-(z))}{
  \tilde \Delta(w,z)}~.
\]
In particular, we have for all $z \in B_c$
\[
  \Bigl|(\PP f_A)(z) - (\PP f_A)_\reg(z)\Bigr| \,\le\, \frac{C}{
  |w_+(z) - w_-(z)|} \,\le\, \frac{C}{|\theta_+(z) - \theta_-(z)|}
  \,\le\, \frac{C}{\dist(z,C_j)^{1/2}}~, 
\]
and this estimate is sharp if $\tilde P(1) = P(1,0) \neq 0$, 
because in that case $\tilde Q(w,z)$ does not vanish near $w = 1$ if 
$z \in B$. Summarizing, we have reached the following 
important conclusion: Near a regular point $\bar z \in C_j$ of 
the critical set $\Sigma_A$ of a generic matrix $A \in \bM_n(\C)$, 
the derivatives of order $n-2$ of the numerical density $f_A$ are 
smooth on the convex side of the curve $C_j$, and blow up like 
$\dist(z,C_j)^{-1/2}$ on the concave side. If $n = 2$, this is 
in full agreement with the explicit formula \eqref{eq:twod3} 
obtained in Section~\ref{ss:twod}. If $n \ge 3$, we deduce 
after integrating that the numerical density $f_A$ is of class
$C^{n-5/2}$ in a neighborhood of such a point $\bar z$.

\medskip
Using the same techniques, it is also possible to study the 
regularity of the numerical density near more singular points
$\bar z \in \Sigma_A$. On the typical example represented in 
Fig.\,\ref{Fig2}, we see that the following two cases have 
to be analyzed:  

\smallskip\noindent i) {\em Crossings}, which occur when 
two critical curves $C_j$ and $C_p$ intersect transversally 
at $\bar z$. This situation can be treated exactly as before, 
except that one has to consider four distinct regions near 
$\bar z$, instead of two. 

\smallskip\noindent ii) {\em Cusps}, which arise whenever one of the
functions $\lambda_j(\theta) + \lambda_j''(\theta)$ has a simple
zero. Here we can repeat the analysis above, assuming that
$\lambda_j(\theta) = \frac{\alpha}{3}\theta^3 + \OO(\theta^4)$ near
$\theta = 0$. In a neighborhood of $\bar z = 0$, we find that
$f_j(\theta,z)$ has either three real roots (for $z$ inside a cuspidal
domain with tip at $\bar z$), or one real and two complex conjugate
roots (outside the cusp). Using the same argument as before, we
conclude that $\PP f_A$ is smooth inside the cusp, but blows up on the
other side of the critical curve $C_j$. Altogether, the numerical
density $f_A$ is of class $C^{n-8/3}$ near the cusp, see
\cite[Section~4.3]{YB} for a similar analysis of the singularities of
the fundamental solution of \eqref{eq:hypsys}.

\smallskip 
Under generic assumptions on the matrix $A$, all
intersections are transversal and all cusps are non degenerate, so
that the singular set $C_A$ is a generic curve in the sense of real
algebraic geometry, and the singularities of the numerical density
$f_A$ can be completely analyzed using the techniques described
bove. In particular, the derivatives of order $n-2$ of the numerical
density $f_A$ are locally integrable, and since we know from
Proposition~\ref{prop:gen} that $f_A \in C^{n-3}$ we conclude that the
relation \eqref{eq:PfA} holds everywhere (in the sense of
distributions), and not only on the complement of $\Sigma_A$.

\subsection{Connexion with the fundamental solution}\label{ss:fundam}

As was explained in the introduction, the numerical measure of a
matrix $A \in \bM_n(\C)$ is related to the fundamental solution of the
hyperbolic system \eqref{eq:hypsys}. This connexion can be established
rigorously by comparing the expression \eqref{eq:PfA} for the
derivatives of the numerical density $f_A$ with the representation
formulas for the fundamental solution $E(t,x)$ of \eqref{eq:hypsys},
which can be found e.g. in \cite{BY,YB}.

In what follows, we identify $\C$ with $\R^2$, and we denote by 
$x = (x_1,x_2)$ the points of the Euclidean plane. With a slight 
abuse of notation, we write $f_A(x)$ instead of $f_A(z)$, and we 
consider as subsets of $\R^2$ the various regions associated with 
$A$, such as $W(A)$ or $\Sigma_A$. To derive a representation 
formula for the fundamental solution $E(t,x)$, we take the Radon
transform of \eqref{eq:hypsys2} with respect to $x \in \R^2$, 
and obtain the one-dimensional hyperbolic system
\[
  \partial_t \tilde E(t,s,\theta) + H(\theta)\partial_s \tilde 
  E(t,s,\theta) \,=\, I_n\,\delta_{t=0} \otimes \delta_{s = 0}~,
\]
where $H(\theta)$ is given by \eqref{def:Htheta} and $\tilde 
E(t,s,\theta)$ denotes the Radon transform of $E(t,x)$. Using the 
method of characteristics, we easily find
\[
  \tilde E(t,s,\theta) \,=\, \sum_{j=1}^n \delta(s - t \lambda_j(\theta))
  P_j(\theta)~, \quad t \ge 0,
\]
where $\lambda_1(\theta),\dots,\lambda_n(\theta)$ are the eigenvalues 
of $H(\theta)$ and $P_1(\theta),\dots,P_n(\theta)$ the corresponding
spectral projections. Now, if we invert the Radon transform as 
in \eqref{eq:Radon3} and use the identity
\[
  \sum_{j=1}^n \frac{1}{s - t\lambda_j(\theta)}\,P_j(\theta) \,=\,
  (sI_n - tH(\theta))^{-1}~,
\]
we arrive at the representation formula
\begin{equation}\label{eq:Erep}
  E(t,x) \,=\, -\frac{1}{4\pi^2}\,\fp \!\int_{S^1} 
  \Bigl((x_1\cos\theta + x_2\sin\theta)I_n - tH(\theta)
  \Bigr)^{-2}\dd\theta~, \quad t \ge 0~,
\end{equation}
which coincides with Eq.~(4.4a) in \cite{BY}. Arguing as in 
Section~\ref{ss:poly}, one can show that equation \eqref{eq:Erep} 
is rigorously satisfied for all $(t,x) \in (0,\infty) \times \R^2$
with $\frac{x}{t} \in \R^2\setminus\Sigma_A$, and that $E(t,x)$ is 
smooth in that region of space-time. 

To compare the numerical density with the fundamental solution, 
it is natural to extend $f_A$ to a homogeneous function of 
space and time by setting
\begin{equation}\label{def:FAhom}
  \FF_A(t,x) \,=\, t^{n-3} f_A\Bigl(\frac{x}{t}\Bigr)~, \qquad 
  t > 0~, \quad x \in \R^2~.
\end{equation}
We then have the following result:

\begin{prop}\label{prop:compEf}
For any $A \in \bM_n(\C)$, there exists a matrix-valued homogeneous
polynomial $Q$ of degree $n-1$ such that
\begin{equation}\label{eq:relEf}
   E(t,x) \,=\, Q(\partial_t,\partial_{x_1},\partial_{x_2})\FF_A(t,x)~, 
\end{equation}
for all $(t,x) \in (0,\infty) \times \R^2$ with $\frac{x}{t} \in 
\R^2\setminus\Sigma_A$. 
\end{prop}

In particular, if $\Omega$ is a connected component of the region
$\Pi_A \subset \R^2$ defined by \eqref{eq:Pidef}, we know from
Corollary~\ref{cor:poly} that $f_A(x)$ is polynomial of degree at
most $n-3$ in $\Omega$, and \eqref{def:FAhom} then shows that
$\FF_A(t,x)$ is also polynomial of degree at most $n-3$ in the
half-cone $C_+(\Omega) = \{(t,x) \in (0,\infty) \times \R^2 \,|\,
\frac{x}{t} \in \Omega\}$. By Proposition~\ref{prop:compEf}, we
conclude that $E(t,x) = 0$ in $C_+(\Omega)$. We thus have

\begin{cor}\label{cor:lac}
The fundamental solution $E(t,x)$ of a matrix $A \in \bM_n(\C)$ 
vanishes for all $(t,x) \in (0,\infty) \times \R^2$ such that
$\frac{x}{t}$ belongs to the region $\Pi_A \subset \R^2$ defined 
by \eqref{eq:Pidef}. 
\end{cor}

In the language of partial differential equations, the {\em domain 
of influence} of the origin for system \eqref{eq:hypsys} is the half-cone
$C_+(D) = \{(t,x) \in (0,\infty) \times \R^2 \,|\, \frac{x}{t}\in D\}$, 
where $D \subset \R^2$ is the complement of the largest connected open 
region on which $E_* = E(1,\cdot)$ vanishes. If in addition $E_* = 0$ 
in some open region $L \subset D$, we say that $L$ is a {\em lacuna} 
of the hyperbolic system \eqref{eq:hypsys}. With this terminology, 
Corollary~\ref{cor:lac} asserts that each connected component of 
$\Pi_A$ either lies outside the domain of influence of the origin, 
or is a lacuna. This important geometric characterization of lacunas
is originally due to Petrovsky \cite{Pet}, and was thoroughly 
discussed in \cite{BY,YB} and from a more algebraic point of 
view in \cite{ABG1,ABG2}. Strictly speaking, it gives only a 
sufficient condition for the occurence of lacunas, but further
work allows to show that all {\em stable} lacunas satisfy this 
criterion. In other words, in generic situations, the converse 
of Corollary~\ref{cor:lac} also holds: each open region on which 
$E_*$ vanishes (in particular, any lacuna) belongs to $\Pi_A$. 

As a typical example, consider the matrix $A \in M_3(\C)$ defined by
\eqref{def:A3gen}, whose numerical density is represented in
Fig.\,\ref{Fig1}. The polynomial region $\Pi_A$ has just two
components here: the exterior of $D = W(A)$, where both $f_A$ and
$E_*$ vanish, and the interior of the cuspidal triangle, where
$E_* = 0$ and $f_A$ is identically constant, which is
therefore a lacuna. A more complicated situation is depicted in
Fig.\,\ref{Fig2}, where the polynomial region $\Pi_A$ has six
connected components, among which five correspond to lacunas. In a
different spirit, it is also interesting to consider the nongeneric
matrix \eqref{def:reduc} whose numerical density is studied in
Section~\ref{s:examples}. The domain of influence of the origin is the
union of the closed unit disk and a single point $\{a\}$, so the
numerical range $W(A) = \co(D)$ is substantially larger than $D$ if
$|a| > 1$, see Fig.\,\ref{Fig5}.  In that case, the region $\Pi_A$ has
again two components, none of which is a lacuna: the exterior of
$W(A)$, and the interior of the triangular region $W(A) \setminus D$
where $f_A$ is identically constant. We see on these examples that the
numerical density allows to distinguish between various regions where
the fundamental solution vanishes identically, and which are
nevertheless of rather different nature.

Before proving Proposition~\ref{prop:compEf}, we briefly verify 
its conclusion on a simple example. If $A \in M_2(\C)$ is defined
by \eqref{eq:defA2}, system \eqref{eq:hypsys} reduces to
\[
  \partial_t u_1 + \partial_{x_1}u_2 -i \partial_{x_2}u_2 \,=\, 0~, 
  \qquad
  \partial_t u_2 + \partial_{x_1}u_1 +i \partial_{x_2}u_1 \,=\, 0~.
\]
Combining both equations, one verifies that $\partial_t^2 u_j = 
\Delta u_j$ for $j = 1,2$, hence the fundamental solution $E(t,x)$ 
can easily be computed using Poisson's formula for the solution of the 
wave equation in two dimensions. In agreement with 
Proposition~\ref{prop:compEf}, the result is:
\[
  E(t,x) \,=\, \begin{pmatrix} \partial_t  & -\partial_{x_1}+i
  \partial_{x_2} \\ -\partial_{x_1} - i\partial_{x_2} & \partial_t
  \end{pmatrix}\FF_A(t,x)~, \qquad |x| < t~,
\]
where according to \eqref{eq:fA2}, \eqref{def:FAhom}
\[
  \FF_A(t,x) \,=\, \frac{1}{t} f_A\Bigl(\frac{x}{t}\Bigr) \,=\,
 \frac{1}{2\pi}\,\frac{1}{\sqrt{t^2-|x|^2}}\,\mathbf{1}_{\{|x| < t\}}~.
\]

\bigskip\noindent{\bf Proof of Proposition~\ref{prop:compEf}.} 
If $n = 1$ the conclusion is trivial, because both $E_*$ and 
$f_A$ vanish identically outside $\Sigma_A$ (which is reduced 
to a single point), so we assume henceforth that $n \ge 2$. 
If $\PP = P(\partial_t,\partial_{x_1},\partial_{x_2})$ is a 
homogeneous differential operator of degree $n-2$, then using 
\eqref{def:FAhom} it is straightforward to verify that
\begin{equation}\label{eq:idFA}
  (\PP \FF_A)(t,x) \,=\, \frac{1}{t} \left[P\Bigl(-\frac{x}{t}
  \cdot \nabla_\xi,\partial_{\xi_1},\partial_{\xi_2}\Bigr)f_A \right]
  \Big|_{\xi = \frac{x}{t}}~,
\end{equation}
whenever $\frac{x}{t} \in \R^2\setminus \Sigma_A$. Here $\xi \in \R^2$ 
denotes the argument of the function $f_A$, which has to be replaced by 
$\frac{x}{t}$ after differentiation. We warn the reader that equality
\eqref{eq:idFA} holds only if $P$ is of degree $n-2$. 
Applying Proposition~\ref{prop:deriv}, we deduce that
\begin{equation}\label{eq:diffit}
  P(\partial_t,\partial_{x_1},\partial_{x_2})\FF_A(t,x) \,=\,
  -\frac{(n{-}1)!}{4\pi^2t} \,\fp \!\int_{S^1} \frac{P(-\frac{x}{t}
  \cdot e_\theta,\cos\theta,\sin\theta)}{\det(H(\theta)-\frac{x}{t}
  \cdot e_\theta\,I_n)}\dd \theta~,
\end{equation}
where $e_\theta = (\cos\theta,\sin\theta)$. On the other hand, 
starting from \eqref{eq:Erep}, we observe that 
\[
  -\Bigl(\frac{x}{t}\cdot e_\theta\,I_n - H(\theta)\Bigr)^{-2}  
  \,=\, \frac{\partial}{\partial s} \Bigl((s + \frac{x}{t}\cdot e_\theta) 
  I_n - H(\theta)\Bigr)^{-1}~\Big|_{s=0}~.
\]
By Cramer's rule, the inverse of the matrix $SI_n - H(\theta) = SI_n 
- A_1\cos\theta - A_2\sin\theta$, with $S = s + \frac{x}{t}\cdot 
e_\theta$, has the following form
\[
  \Bigl(SI_n - H(\theta)\Bigr)^{-1} \,=\, 
 \frac{1}{\bar\Delta(\theta,S)}\Bigl(S\,P_0(S,\cos\theta,\sin\theta) 
 + \cos\theta\,P_1(S,\cos\theta,\sin\theta) + \sin\theta\,P_2(S,
 \cos\theta,\sin\theta)\Bigr)~,
\]
where $\bar\Delta(\theta,S) = \det(H(\theta)-SI_n)$ and 
$P_0,P_1,P_2$ are matrix-valued homogeneous polynomials of degree
$n-2$. The idea is now to insert this expansion into the right-hand 
side of \eqref{eq:Erep} and to use \eqref{eq:diffit} to express 
the result as a derivative of order $n-1$ of the numerical density 
$\FF_A$. 

We begin with the term involving $P_1$, and remark that
\[
  \cos\theta \,\frac{\partial}{\partial s} \left(\frac{P_1(s + 
  \frac{x}{t}\cdot e_\theta,\cos\theta,\sin\theta)}{\bar\Delta(\theta,
  s + \frac{x}{t}\cdot e_\theta)}\right)~\Big|_{s=0} \,=\, 
  t\,\frac{\partial}{\partial x_1}\left(\frac{P_1(\frac{x}{t}\cdot 
  e_\theta,\cos\theta,\sin\theta)}{\bar\Delta(\theta,\frac{x}{t}
  \cdot e_\theta)}\right)~.
\]
The corresponding contribution to \eqref{eq:Erep} is thus
\[
  E_1(t,x) \,=\, \frac{1}{4\pi^2 t}\,\frac{\partial}{\partial x_1}
  \,\fp \!\int_{S^1}\frac{P_1(\frac{x}{t}\cdot e_\theta,\cos\theta,\sin\theta)}
  {\bar\Delta(\theta,\frac{x}{t}\cdot e_\theta)}\dd\theta \,=\, 
  -\frac{1}{(n{-}1)!}\,\frac{\partial}{\partial x_1} 
  P_1(-\partial_t,\partial_{x_1},\partial_{x_2})\FF_A(t,x)~.
\]
Similarly, the term involving $P_2$ gives the contribution
\[
  E_2(t,x) \,=\, -\frac{1}{(n{-}1)!}\,\frac{\partial}{\partial x_2} 
  P_2(-\partial_t,\partial_{x_1},\partial_{x_2})\FF_A(t,x)~.
\]
Finally, to treat the expression containing $P_0$, we observe that
\[
  \frac{\partial}{\partial s} \left(\frac{(s + 
  \frac{x}{t}\cdot e_\theta)P_0(s + \frac{x}{t}\cdot e_\theta,
  \cos\theta,\sin\theta)}{\bar\Delta(\theta,s + \frac{x}{t}\cdot 
  e_\theta)}\right)~\Big|_{s=0} \,=\, (1+x\cdot\nabla_x) \left(\frac{
  P_0(\frac{x}{t}\cdot e_\theta,\cos\theta,\sin\theta)}{\bar\Delta(
  \theta,\frac{x}{t}\cdot e_\theta)}\right)~,
\]
and using \eqref{eq:idFA}, \eqref{eq:diffit} we obtain the contribution
\[
  E_0(t,x) \,=\, -\frac{1}{(n{-}1)!\,t}\,(1 + x\cdot\nabla_x)
  P_0(-\partial_t,\partial_{x_1},\partial_{x_2})\FF_A(t,x) \,=\, 
  \frac{1}{(n{-}1)!} \,\frac{\partial}{\partial t} P_0(-\partial_t,
  \partial_{x_1},\partial_{x_2})\FF_A(t,x)~.  
\]
Recalling that $E(t,x) = E_0(t,x) + E_1(t,x) + E_2(t,x)$, we 
arrive at \eqref{eq:relEf}. 
\enpr

\section{Three-dimensional examples}\label{s:examples}

To illustrate the results of the previous sections, we consider here  
four concrete examples which, according to \cite[Section~I.7]{Kip}, 
give a complete picture of what can happen for three-dimensional 
matrices. The two-dimensional case, which is much simpler, was already 
treated in Section~\ref{ss:twod}, and the references \cite{ABG2,BY,JAG}
include a detailed study of the singular set $\Sigma_A$ for a few
higher-dimensional examples.

\begin{figure}
\begin{tabular}{cc}
\includegraphics[width=8.0cm,height=5cm]{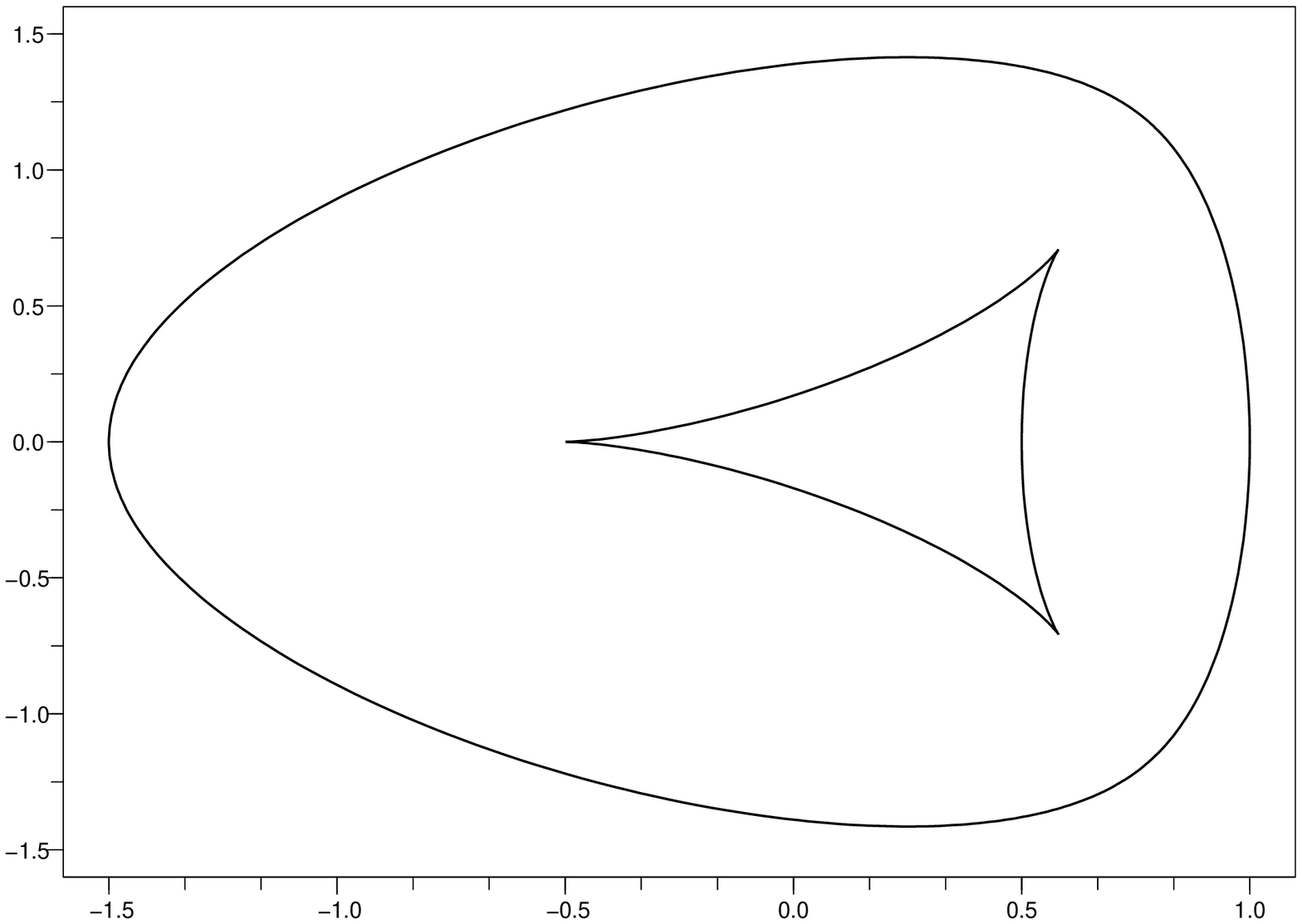} & 
\includegraphics[width=8.0cm,height=5cm]{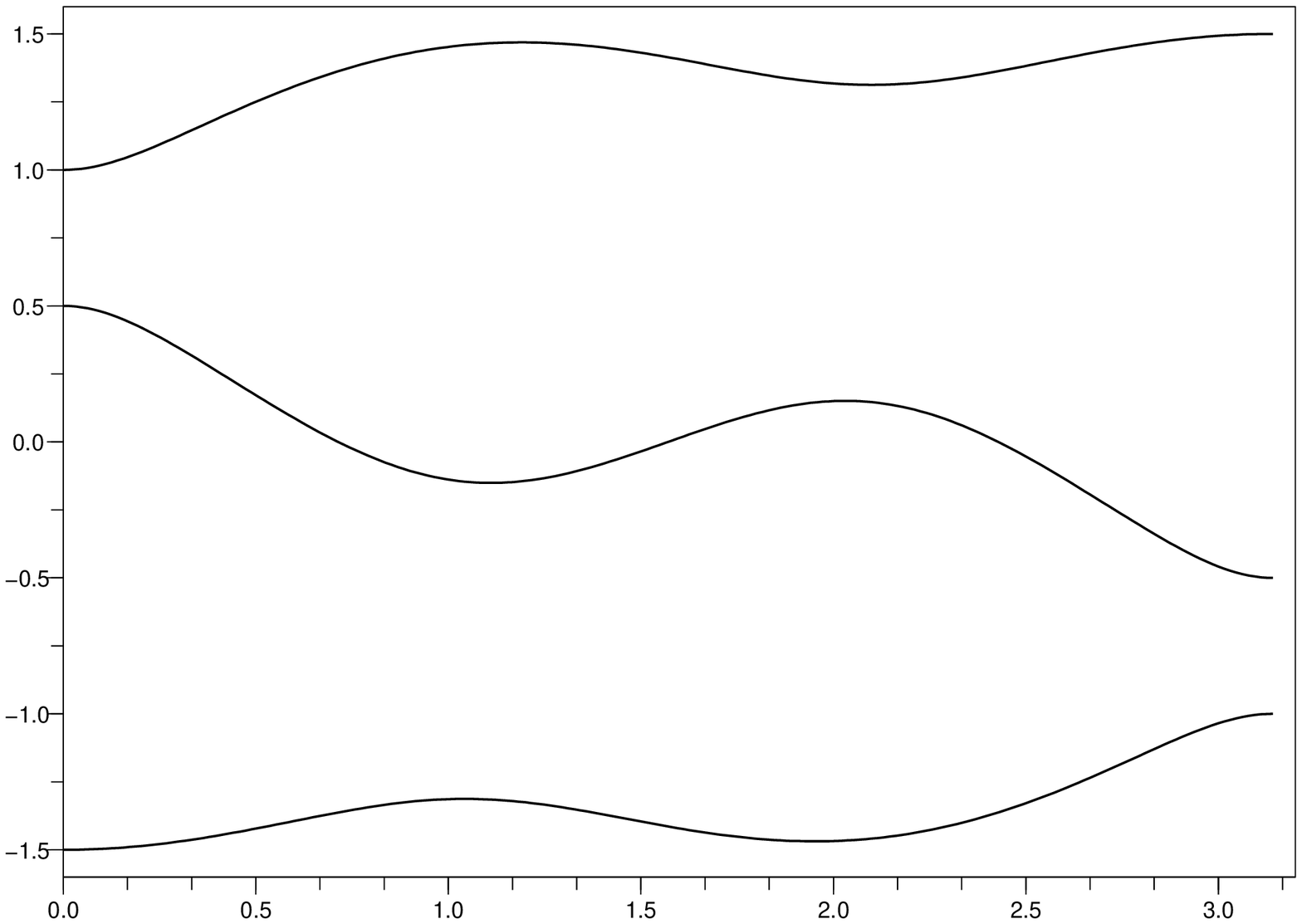} \\
\end{tabular}
\caption{\small The critical set $\Sigma_A$ (left) and the eigenvalues 
$\lambda_j(\theta)$ (right) are represented for the matrix 
\eqref{def:A3gen}.}
\label{Fig3}
\end{figure}

\medskip\noindent{\bf Example~1.} We first consider the $3\times 3$ matrix 
\begin{equation}\label{def:A3gen}
  A \,=\, \begin{pmatrix} -1.5 & 1 & 0 \\ -1 & 1 & 1 \\ 0 & -1 & 0.5
  \end{pmatrix}~,
\end{equation}
which is generic in the sense that hypotheses H1, H2 in 
Section~\ref{ss:regular} are fulfilled. The eigenvalues of the 
associated Hermitian matrix $H(\theta)$ satisfy $\lambda_1(\theta) 
< \lambda_2(\theta) < \lambda_3(\theta)$ for all $\theta \in [0,\pi]$, 
as can be seen from Fig.\,\ref{Fig3} (right). The critical set 
$\Sigma_A = C_A = C_1 \cup C_2$ is an algebraic curve of degree $6$ 
consisting of a smooth ovate curve $C_1$ enclosing a cuspidal triangle 
$C_2$ (Fig.\,\ref{Fig3}, left). The component $C_1$ corresponds to the
eigenvalues $\lambda_1(\theta), \lambda_3(\theta)$ while $C_2$ is 
associated with $\lambda_2(\theta)$. All multiplicities are equal to 
one, and the number $N(z)$ defined by \eqref{def:N} is equal to $3$ 
outside $C_1$ and inside $C_2$, and to $1$ in the intermediate region. 
The numerical density $f_A$ is continuous, identically constant inside
$C_2$, and vanishes outside $C_1$. Moreover $f_A$ is H\"older 
continuous with exponent $1/2$ across $C_1$ and $C_2$, except at the
cusps. The level lines of $f_A$ are represented in Fig.\,\ref{Fig1}. 

\begin{figure}[!h]
\begin{tabular}{cc}
\includegraphics[width=8.0cm,height=5cm]{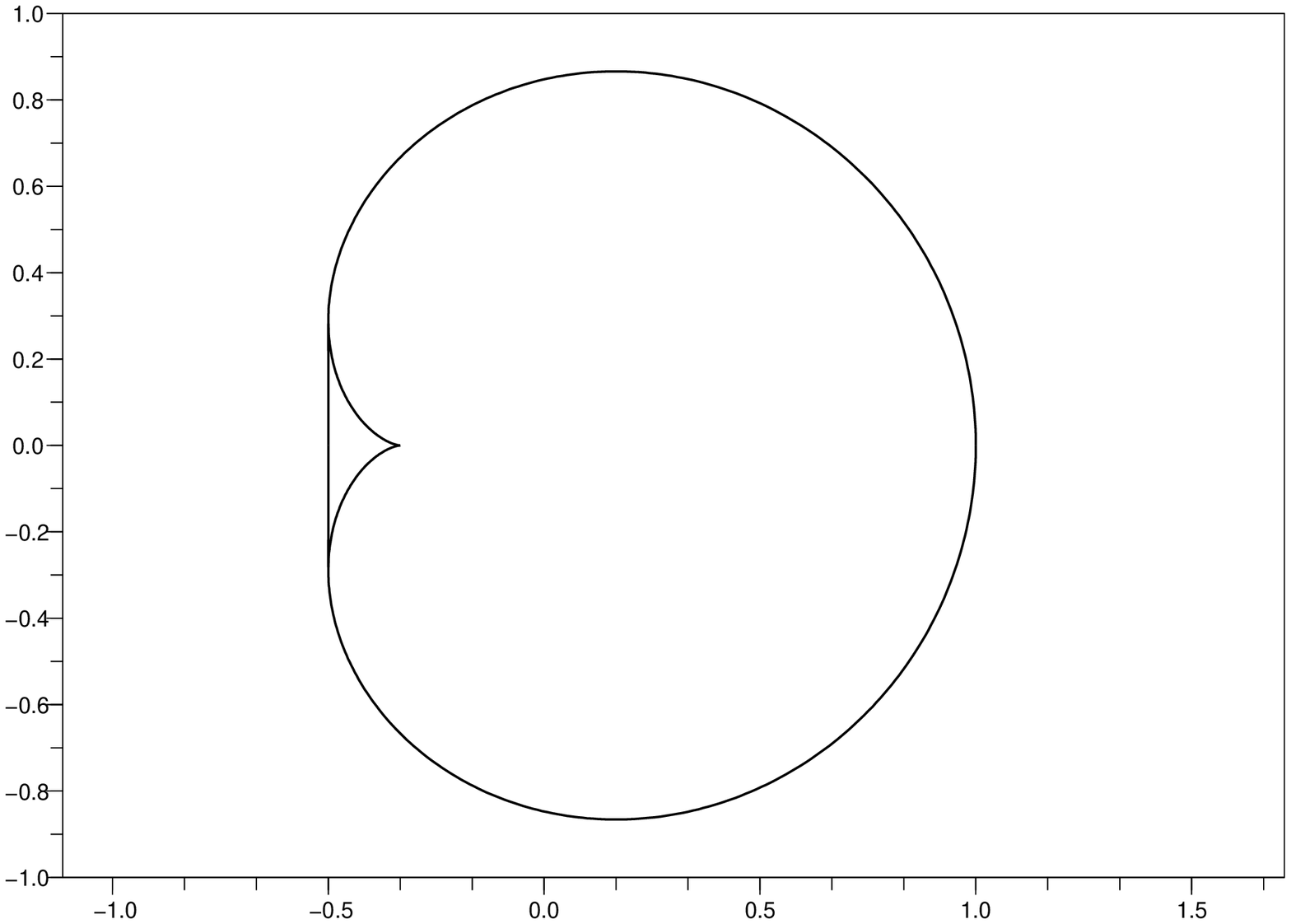} & 
\includegraphics[width=8.0cm,height=5cm]{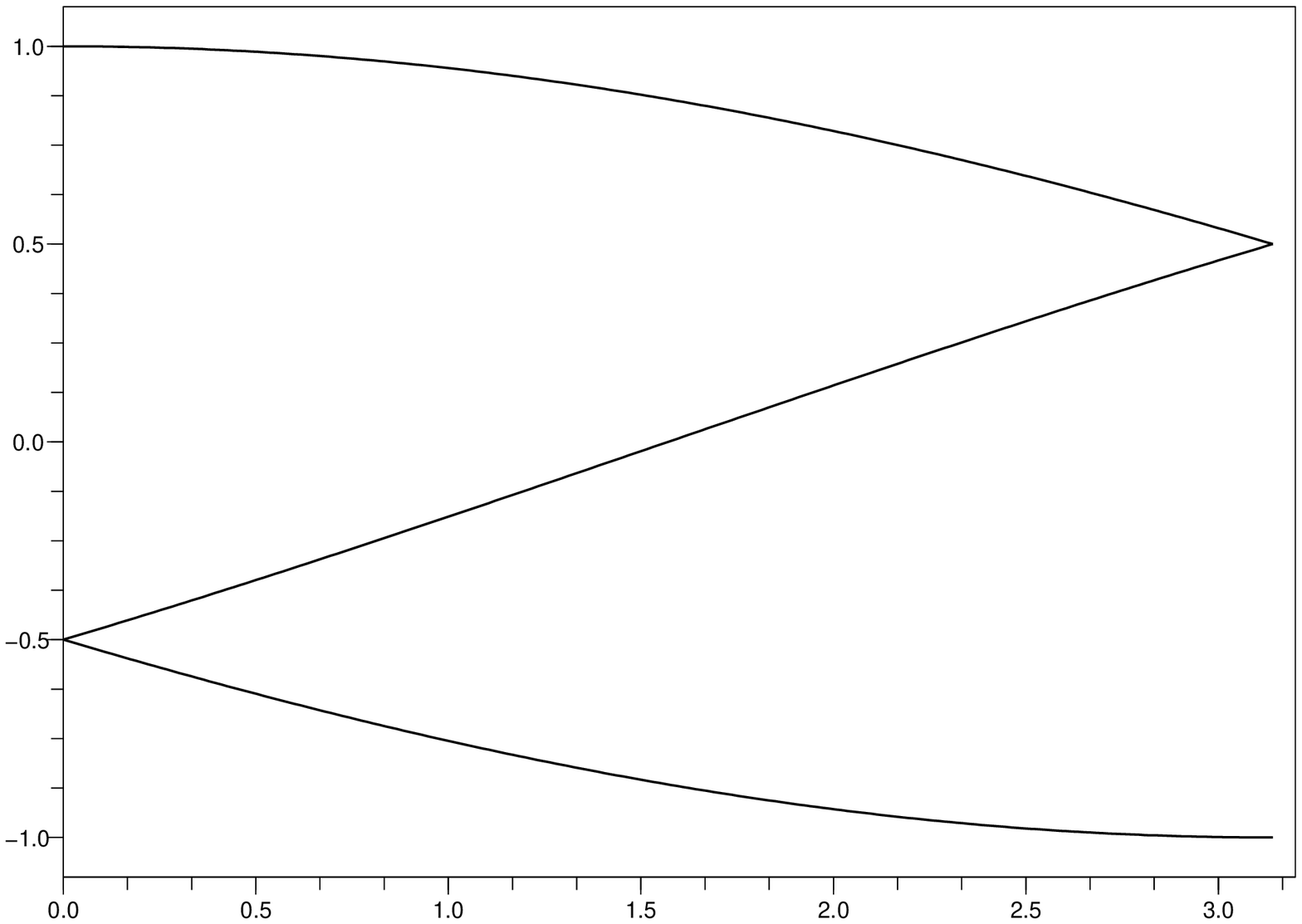} \\
\end{tabular}
\caption{\small The critical set $\Sigma_A$ (left) and the eigenvalues 
$\lambda_j(\theta)$ (right) are represented for the matrix 
\eqref{def:cardio}.}
\label{Fig4}
\end{figure}

\medskip\noindent{\bf Example~2.} We next consider a nongeneric matrix 
\begin{equation}\label{def:cardio}
  A \,=\, \begin{pmatrix} 0 & 1 & 1 \\ 0 & 0 & 1 \\ 0 & 0 & 0
  \end{pmatrix}~,
\end{equation}
for which the critical set $\Sigma_A$ can be computed exactly. Indeed, 
if $H(\theta)$ is the Hermitian matrix \eqref{def:Htheta}, it is easy 
to verify that $\det(\lambda I_3 - H(\theta)) = \lambda^3 -\frac34\lambda
-\frac14\cos\theta$, hence
\[
  \lambda_1(\theta) \,=\, \cos\Bigl(\frac{\theta+2\pi}{3}\Bigr)~, \quad 
  \lambda_2(\theta) \,=\, \cos\Bigl(\frac{\theta-2\pi}{3}\Bigr)~, \quad 
  \lambda_3(\theta) \,=\, \cos\Bigl(\frac{\theta}{3}\Bigr)~.
\]
Thus the permutation \eqref{def:tau} is just a cycle $\tau =
(1\,2\,3)$, and applying \eqref{eq:zj} we easily find that the
critical curve $C_A$ is the cardioid defined by $C_A = \{\frac13 (2
e^{i\phi} + e^{2i\phi})\,|\, \phi \in S^1\}$. Since $\lambda_1(0) =
\lambda_2(0)$, the bitangent set $C_A'$ is not empty and consists of
the line segment joining the points $-1/2\pm i/(2\sqrt{3})$, see
Fig.\,\ref{Fig4} (left).  Altogether we have $\Sigma_A = C_A \cup
C_A'$, and we observe that $\Sigma_A$ encloses a convex region of the
complex plane which is of course the numerical range $W(A)$. The index
$N(z)$ defined by \eqref{def:N} is equal to $3$ outside $C_A$ and to
$1$ inside. The numerical density $f_A$ vanishes outside $W(A)$ and is
equal to a nonzero constant inside the cuspidal region, in agreement
with Propositions~\ref{prop:pos} and \ref{prop:deriv}. In particular,
$f_A$ is discontinuous along the line segment $C_A'$.

\medskip\noindent{\bf Example~3.} The matrices considered so far 
were unitarily irreducible. In contrast, the matrix 
\begin{equation}\label{def:reduc}
  A \,=\, \begin{pmatrix} 0 & 2 & 0 \\ 0 & 0 & 0 \\ 0 & 0 & a
  \end{pmatrix}~, \qquad \hbox{where } a \in \C~,
\end{equation}
is the direct orthogonal sum of the two-dimensional Jordan 
block \eqref{eq:defA2} and the one-dimensional matrix $(a)$. The 
numerical density $f_A$ can therefore be computed using 
Proposition~\ref{prop:AplusB}, and without loss of generality 
we can assume that $a \ge 0$. However, we have to distinguish 
between three cases:

\medskip\noindent{i)} If $0 \le a < 1$, the numerical range 
$W(A)$ is the closed unit disk, and the numerical density has the
following expression:  
\[
  f(z) \,=\, \frac{1}{\pi\sqrt{1-a^2}}\,\mathrm{argch}\left(\frac{1-az_1}
  {\sqrt{(1-az_1)^2 - (1-|z|^2)(1-a^2)}}\right)~, \qquad |z| < 1~,
\]
which reduces to \eqref{eq:fA3} when $a = 0$. In particular $f_A$ 
vanishes on the unit circle, has a logarithmic singularity 
at the point $\{a\}$, and is otherwise smooth. The singular 
set $\Sigma_A$ is the union of the unit circle and the point 
$\{a\}$. 

\medskip\noindent{ii)} In the limiting case $a = 1$, the numerical
range is still the closed unit disk, but the formula
\[
  f(z) \,=\, \frac{1}{\pi}\frac{\sqrt{1-|z|^2}}{1-z_1}~, \quad
  \qquad |z| < 1~,
\]
shows that the numerical density has now an algebraic singularity at 
the boundary point $z = 1$. 

\begin{figure}[!h]
\begin{tabular}{cc}
\includegraphics[width=8.0cm,height=5cm]{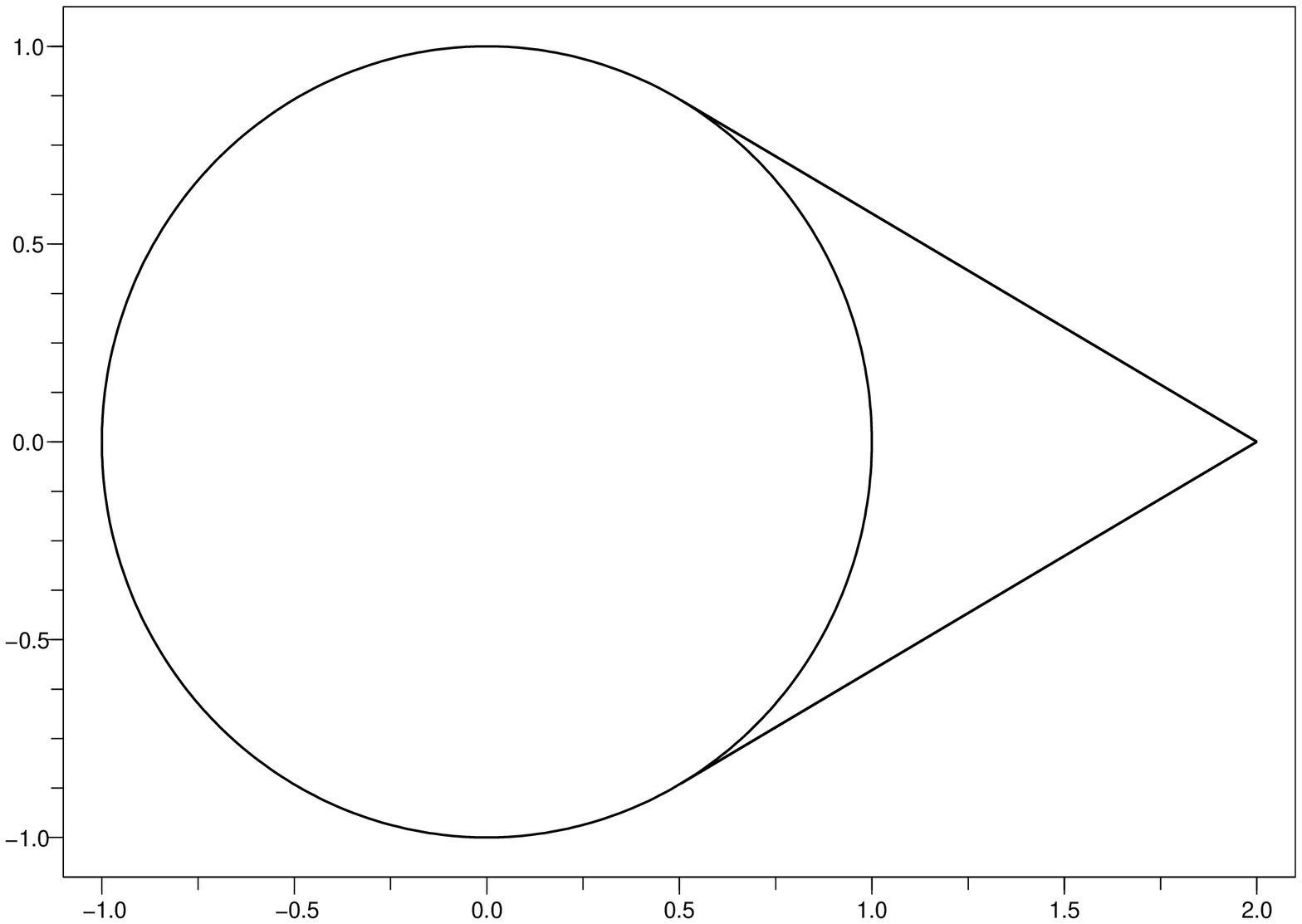} & 
\includegraphics[width=8.0cm,height=5cm]{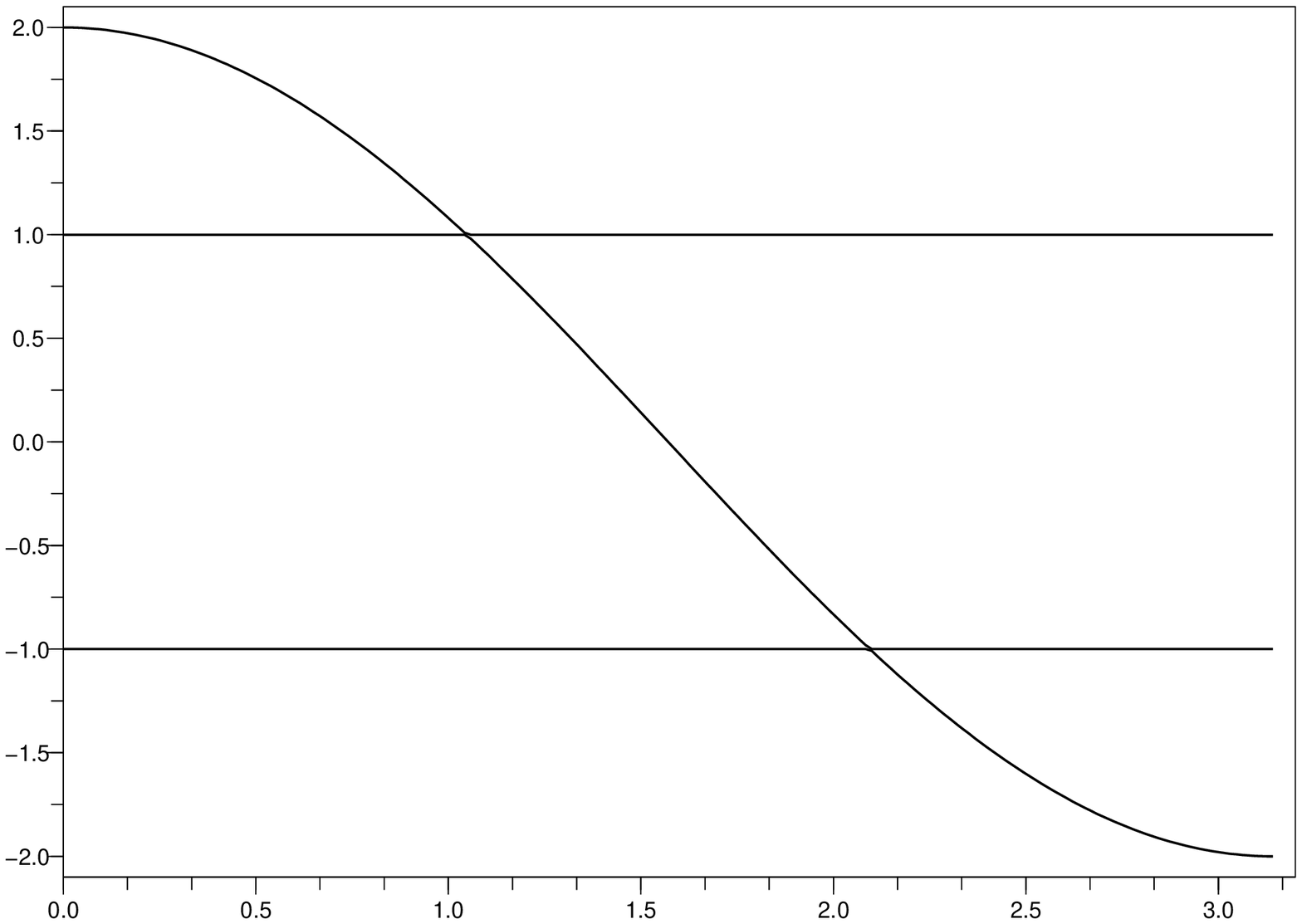} \\
\end{tabular}
\caption{\small The critical set $\Sigma_A$ (left) and the eigenvalues 
$\lambda_j(\theta)$ (right) are represented for the matrix 
\eqref{def:reduc} with $a = 2$.}
\label{Fig5}
\end{figure}

\medskip\noindent{iii)} When $a > 1$, the numerical range $W(A)$ 
is the convex hull of the union of the unit disk and the exterior 
point $\{a\}$. Within this region, the numerical density satisfies
\[
  f(z) \,=\, \frac{1}{\pi\sqrt{a^2-1}}\,\arccos\left(\frac{1-az_1}
  {\sqrt{(1-az_1)^2 + (1-|z|^2)(a^2-1)}}\right)~, \qquad 
  \hbox{when} \quad |z| < 1~,
\]
and $f(z) = 1/\sqrt{a^2-1}$ when $|z| > 1$. As is easily verified, 
the eigenvalues of the Hermitian matrix $H(\theta)$ are 
$\lambda_1(\theta) = -1$, $\lambda_2(\theta) = 1$, and $\lambda_3(\theta) 
= a\cos\theta$, see Fig.\,\ref{Fig5} (right). The algebraic curve 
$C_A$ consists of the unit circle (associated with $\lambda_1,\lambda_2$)
and the point $\{a\}$ (corresponding to $\lambda_3$), but the 
bitangent set $C_A'$ is not empty and consists of two line segments, 
see Fig.\,\ref{Fig5} (left). 

\medskip\noindent{\bf Example~4.} As a final example, we consider
the case of a normal matrix $A \in M_3(\C)$ whose eigenvalues 
$\lambda_1, \lambda_2, \lambda_3$ are not colinear. Then the 
numerical range is the triangle with vertices $\{\lambda_1,\lambda_2,
\lambda_3\}$ and the numerical density is a multiple of the 
characteristic function of $W(A)$. In that situation $\Sigma_A = 
C_A \cup C_A'$, where $C_A$ is the set of all vertices and 
$C_A'$ the set of all edges of the triangle. 

\section{Statistical properties of the numerical measure}\label{s:stat}

In this section, we study the numerical measure from a statistical
point of view, and we establish various convergence results which show
that the measure $\mu_A$ of a large matrix $A \in \bM_n(\C)$ is 
concentrated in a neighborhood of the barycenter of the spectrum 
$\sigma(A)$. 

\subsection{Concentration phenomena}

\begin{prop}\label{prop:stat}
For any $A\in \bM_n(\C)$, the first moment of the probability measure
$\mu_A$ is the normalized trace of $A$:
\[
  \overline{\mu}_A \,=\, \int_\C z\dd \mu_A(z) \,=\, \frac{1}{n}\,\Tr(A)~,
\]
and the variance of $\mu_A$ is given by 
\begin{equation}\label{eq:var}
  \Var(\mu_A) \,=\, \int_\C |z-\overline{\mu}_A|^2\dd\mu_A(z) \,=\, 
  \frac{1}{n+1}\,\Bigl(\frac{1}{n}\,\Tr(A^*A)-\Bigl|\frac{1}{n}\,
  \Tr(A)\Bigr|^2\Bigr)~.
\end{equation}
\end{prop}

\bepr
Let $a_{jk}$ denote the entries of $A$. Applying definition \eqref{def:muA}
with $\phi(z) = z$, we have to compute the average of $\sum_{jk}a_{jk}
\overline{x_j}x_k$ over the unit sphere. By symmetry, the average 
of $\overline{x_j}x_k$ is equal to zero if $j\ne k$ and to $1/n$ 
if $j=k$. Thus
\[
  \overline{\mu}_A  \,=\, \frac{1}{n}(a_{11}+\dots+a_{nn}) \,=\, 
  \frac{1}{n}\,\Tr(A)~.
\]
To compute the second moment, we take $\phi(z) = |z|^2$ and proceed 
in exactly the same way. By symmetry, the average of $\overline{x_j}x_k
\overline{x_\ell}x_m$ is zero unless $j = k$ and $\ell = m$, 
or $j = m$ and $\ell = k$. Moreover, it is easy to verify that the
average of $|x_j|^2|x_\ell|^2$ is equal to $r_n$ if $j = \ell$ and 
to $s_n$ if $j \ne \ell$, where
\[
  r_n \,=\, \int_{\partial\B^n}|x_1|^4 \dd\bar\sigma(x) \,=\, 
  \frac{2}{n(n+1)}~, \qquad
  s_n \,=\, \int_{\partial\B^n}|x_1|^2 |x_2|^2 \dd\bar\sigma(x) 
  \,=\, \frac{1}{n(n+1)}~.
\]
Thus
\begin{align*}
  \int_\C |z|^2\dd\mu_A(z) \,&=\, \int_{\partial\B^n} 
  |\langle Ax,x\rangle|^2 \dd\bar\sigma(x) \,=\,
  \sum_{j,k,\ell,m} \overline{a_{kj}}
  a_{\ell m} \int_{\partial\B^n}\overline{x_j}x_k \overline{x_\ell}x_m
  \dd\bar\sigma(x) \\
  \,&=\, s_n \sum_{j \ne \ell}(\overline{a_{jj}}a_{\ell\ell}+|a_{\ell j}|^2) 
  + r_n \sum_{j=1}^n |a_{jj}|^2\\
  \,&=\, s_n\sum_{j,\ell=1}^n(\overline{a_{jj}}a_{\ell\ell}+|a_{\ell j}|^2)
  \,=\, s_n(|\Tr A|^2+\Tr(A^*A))~.
\end{align*}
This gives the desired result, since $\Var(\mu_A) = \int_\C |z|^2\dd\mu_A(z) 
- |\overline{\mu}_A|^2$. 
\enpr

\bigskip
Now we consider a sequence of matrices $\{A_n\}_{n \ge 1}$ such that
$A_n\in \bM_n(\C)$ for each $n \ge 1$. As is well known, we have
\[
  \frac{1}{n}\,\Tr(A_n^*A_n) \,\le\, \|A_n\|^2 \,\le\, \Tr(A_n^*A_n)~,
\]
where $\|A_n\| = \sup\{\|A_n x\|\,|\, x \in \partial \B^n\}$. As a 
consequence, if we suppose that $\|A_n\|^2 = o(n)$ as $n \to \infty$, 
it follows from \eqref{eq:var} that the variance of $\mu_{A_n}$ converges
to zero as $n \to \infty$. This gives:

\begin{cor}\label{c:weakconv}
Assume that $A_n\in \bM_n(\C)$ and that $\|A_n\|^2/n \to 0$ as 
$n \to \infty$. Then the measure $\mu_{A_n} - \delta_{\overline{\mu}_{A_n}}$
converges weakly to zero as $n \to \infty$. 
\end{cor}

We recall that the {\em numerical radius} of a matrix $A \in \bM_n(\C)$ 
is defined by
\[
  R(A) \,=\, \sup\{|z|\,|\, z \in W(A)\} \,=\, \sup\{|\langle Ax,x\rangle 
  |\,|\, x \in \partial\B^n\}~,
\]
and satisfies $R(A) \le \|A\| \le 2R(A)$ \cite{HJ}. Thus, a sequence 
of matrices $A_n \in \bM_n(\C)$ is uniformly bounded (in the operator
norm) if and only if the numerical ranges $W(A_n)$ are all contained
in a bounded region of the complex plane. Under this assumption,  
Proposition~\ref{prop:stat} shows that the variance of $\mu_{A_n}$ 
is $\OO(1/n)$ as $n \to \infty$, so that the numerical measure is 
asymptotically concentrated in a disk of radius $\OO(1/\sqrt{n})$ 
around the mean $\overline{\mu}_{A_n}$. 

\bigskip In the introduction, we have observed that the numerical
measure $\mu_{A_n}$ is the distribution of the random variable
$\langle A_n X_n, X_n\rangle \in \C$ when the vector $X_n$ is uniformly
distributed on the unit sphere $\partial\B^n$. With this interpretation, 
Corollary~\ref{c:weakconv} is reminiscent of the weak law of large 
numbers in probability theory. Under slightly stronger assumptions,
it is also possible to obtain a pointwise convergence result in the
spirit of the strong law of large numbers. Without loss of generality, 
we assume from now on that $\Tr(A_n) = 0$ for all $n \ge 1$, so that 
the measure $\mu_{A_n}$ is centered at the origin. 

\begin{prop}\label{prop:sconv}
Assume that $A_n \in \bM_n(\C)$ satisfies $\Tr(A_n) = 0$ for 
all $n \ge 1$ and
\begin{equation}\label{eq:scond}
  \sup_{n \ge 1}\frac{(\log n)\,\|A_n\|}{n^{1/2}} \,<\, \infty~.
\end{equation}
If for each $n \ge 1$ the random variable $X_n$ is uniformy distributed
on the unit sphere $\partial\B^n$, then $\langle A_n X_n, X_n\rangle$ 
converges almost surely to zero as $n \to \infty$. 
\end{prop}

\bepr
It is clearly sufficient to prove the result for Hermitian matrices
$A_n$, because the general case then follows by considering the 
real and imaginary parts of $\langle A_n X_n, X_n\rangle$. We thus 
assume that $A_n = A_n^*$ for all $n \ge 1$, and we denote by 
$\lambda_{n,1},\dots,\lambda_{n,n}$ the eigenvalues of $A_n$. 
For each $n \ge 1$, let $Y_{n,1},\dots,Y_{n,n}$ be independent
and identically distributed complex random variables with 
density function $f_Y(z) = \pi^{-1}e^{-|z|^2}$, $z \in \C$. In 
particular, we have $E(|Y_{n,m}|^{2k}) = k!$ for each $k \in \N$. 
Since the Euclidean measure on $\partial\B^n$ is the projection on 
the unit sphere of the standard Gaussian measure in $\C^n$, we 
obtain a uniformly distributed random variable on $\partial\B^n$
by setting $X_n = U_n Y_n/\|Y_n\|$, where $Y_n = (Y_{n,1},\dots,
Y_{n,n})^\top$ and $U_n \in \bU_n(\C)$ is a unitary matrix such that 
$U_n^* A_n U_n = \diag(\lambda_{n,1},\dots,\lambda_{n,n})$. Thus
\[
  \langle A_n X_n, X_n\rangle \,=\, \frac{\lambda_{n,1}|Y_{n,1}|^2 + \dots 
  + \lambda_{n,n}|Y_{n,n}|^2}{|Y_{n,1}|^2 + \dots + |Y_{n,n}|^2} 
  \,=\, \frac{P_n}{Q_n}~,
\]
where
\[
  P_n \,=\, \frac1n \sum_{m=1}^n \lambda_{n,m}|Y_{n,m}|^2 \,=\, 
  \frac1n \sum_{m=1}^n \lambda_{n,m}(|Y_{n,m}|^2-1)~, \qquad 
  Q_n \,=\, \frac1n \sum_{m=1}^n |Y_{n,m}|^2~.
\]
By the strong law of large numbers, the denominator $Q_n$ converges 
almost surely to $1$ as $n \to \infty$, hence it remains to show that 
the numerator $P_n$ converges almost surely to zero. But this follows
from classical theorems on the limiting behavior of weighted sums 
of independent random variables, see \cite{CL,Tei}. Since for 
each $n \ge 1$ the random variables $X_{n,m} = |Y_{n,m}|^2-1$ 
($1 \le m \le n$) are independent, have zero mean and finite 
second order moment, and since by \eqref{eq:scond} the coefficients 
$a_{n,m} = n^{-1}\lambda_{n,m}$ satisfy
\[
  \max_{1 \le m \le n}|a_{n,m}| \,\le\, \frac{C}{n^{1/2}\log n}~, 
\]
the results of \cite[Section~3]{Tei} imply that $P_n = \sum_{m=1}^n
a_{n,m}X_{n,m}$ converges almost surely to zero as $n \to \infty$.  
\enpr

\subsection{Central limit theorems}

The results established so far show that for a sequence of traceless
matrices $A_n \in \bM_n(\C)$ the numerical measure $\mu_{A_n}$ tends
to concentrate on the origin as $n \to \infty$. Under stronger 
assumptions, we now prove that the rescaled measure $\mu_{\sqrt{n}A_n}$
converges to a Gaussian distribution, as in the classical central
limit theorem. We first consider the Hermitian case, which is 
somewhat simpler. 
 
\begin{prop}\label{p:CLH}
Let $\{A_n\}_{n\ge1}$ be a sequence of Hermitian matrices such that 
$A_n \in \bM_n(\C)$ and $\Tr(A_n) = 0$. We assume that 
\begin{equation}\label{eq:lcond}
  \frac{(\log n)\,\|A_n\|}{n^{1/2}} \,\xrightarrow[n\to\infty]{}\, 0~, 
  \qquad \hbox{and}\quad 
  \frac{1}{n}\,\Tr(A_n^2) \,\xrightarrow[n\to\infty]{}\, \sigma^2 > 0~.
\end{equation}
Then the rescaled numerical measure $\mu_{\sqrt{n}A_n}$ converges 
weakly to the normal distribution $\NN(0,\sigma^2)$ as $n \to \infty$. 
\end{prop}

\bepr
We use the same notations as in the proof of Proposition~\ref{prop:sconv}.
For each $n \ge 1$, the numerical measure $\mu_{A_n}$ is the distribution 
of the random variable
\[
  Z_n \,=\, \frac{\lambda_{n,1}|Y_{n,1}|^2 + \dots + \lambda_{n,n}
  |Y_{n,n}|^2}{|Y_{n,1}|^2 + \dots + |Y_{n,n}|^2} 
  \,=\, \frac{P_n}{Q_n}~,
\]
where $\lambda_{n,1},\dots,\lambda_{n,n}$ denote the eigenvalues of 
$A_n$ and $Y_{n,1},\dots,Y_{n,n}$ are independent complex random 
variables with density function $f_Y(z) = \pi^{-1}e^{-|z|^2}$. 
Since $Q_n$ converges almost surely to $1$ as $n \to \infty$, we 
have to show that $n^{1/2}P_n$ converges in law to 
$\NN(0,\sigma^2)$. 

To do that, we use the Lindeberg-Feller theorem for triangular 
arrays of random variables \cite[Section~2.4.b]{Dur}. Let
\[
  X_{n,m} \,=\, \frac{1}{\sqrt{n}}\,\lambda_{n,m}\,(|Y_{n,m}|^2-1)~,
  \qquad 1 \le m \le n~,
\]
so that $n^{1/2} P_n = X_{n,1} + \dots + X_{n,n}$.  For each fixed 
$n \ge 1$, the random variables $X_{n,m}$ are independent and satisfy 
$E(X_{n,m}) = 0$ for $m = 1,\dots,n$. Moreover, 
\[
  \sum_{m=1}^n E(|X_{n,m}|^2) \,=\, \frac{1}{n} \sum_{m=1}^n 
  \lambda_{n,m}^2 \,=\, \frac{1}{n}\,\Tr(A_n^2) 
  \,\xrightarrow[n\to\infty]{}\, \sigma^2 > 0~.
\]
Finally, for any $\epsilon > 0$, we have
\[
  E(|X_{n,m}|^2\,;\, |X_{n,m}| \ge \epsilon) \,=\, 
 \frac{\lambda_{n,m}^2}{\pi n}\int_{D_{n,m,\epsilon}}
 (|z|^2-1)^2\,e^{-|z|^2}\dd z~, 
\]
where $D_{n,m,\epsilon} = \{z \in \C\,|\, \lambda_{n,m}^2 (|z|^2{-}1)^2 
\ge n \epsilon^2\}$. Thus, using the first assumption in 
\eqref{eq:lcond}, we obtain by a direct calculation
\[
  \sum_{m=1}^n E(|X_{n,m}|^2\,;\, |X_{n,m}| \ge \epsilon)
  ~\le~ n \sup_{1\le m \le n}E(|X_{n,m}|^2\,;\, |X_{n,m}| \ge 
  \epsilon) \,\xrightarrow[n\to\infty]{}\, 0~.
\]
Invoking the Lindeberg-Feller theorem, we conclude that $n^{1/2}P_n
= X_{n,1} + \dots + X_{n,n}$ converges in law to $\NN(0,\sigma^2)$, 
which is the desired result. 
\enpr

\bigskip\noindent{\bf Remark.} Since the numerical density of a
Hermitian matrix is a $B$-spline, Proposition~\ref{p:CLH} shows under
very general assumptions that $B$-splines of degree $n$ satisfy a
central limit theorem in the limit $n \to \infty$. In the particular
case of uniform $B$-splines, this result was obtained by Unser
{\sl et.\thinspace al.} in \cite{Uns}. 

\bigskip
Before considering more general matrices, we would like to mention
an alternative proof of Proposition~\ref{p:CLH} which has its own 
interest. The starting point is a very nice formula for the 
moments of the numerical measure of a Hermitian matrix $A \in 
\bM_n(\C)$. Fix $k \in \N$ and let $\lambda_1,\dots,\lambda_n \in \R$ 
denote the eigenvalues of $A$. Using \eqref{eq:normrep} or 
\eqref{eq:normrep2} with $\phi(x) = x^k$, we find
\[
  \int_\R x^k \dd\mu_A(x) \,=\, (n-1)! \int_{D_{n-1}} 
  (t_1 \lambda_1 + \dots + t_n \lambda_n)^k \dd t_1 \dots 
  \dd t_{n-1}~,
\]
where $D_{n-1}$ is the $(n{-}1)$-dimensional simplex defined in 
\eqref{def:Dk} and $t_n = 1 - (t_1 + \dots + t_{n-1})$. To 
evaluate the right-hand side, we apply the multinomial formula
\[
  (X_1+\cdots+X_n)^k \,=\, \sum_{|\alpha|=k} \frac{k!}
  {\alpha!}\,X^\alpha~,
\]
where the sum runs over all multi-indices $\alpha \in \N^n$ of 
order $|\alpha| = \alpha_1+\cdots+\alpha_n = k$. Here we use 
the standard notations $X^\alpha = X_1^{\alpha_1}\cdots X_n^{\alpha_n}$ 
and $\alpha! = (\alpha_1!)\cdots(\alpha_n!)$. Now, it is not difficult
to verify that
\[
  \int_{D_{n-1}} t^\alpha\dd t_1\dots \dd t_{n-1} \,=\, 
  \frac{\alpha!}{(n+k-1)!}~,
\]
for any $\alpha \in \N^n$ with $|\alpha| = k$. We thus obtain the 
following identity
\begin{equation}\label{eq:momk}
  \int_\R x^k \dd\mu_A(x) \,=\, \frac{k!\,(n-1)!}{(n+k-1)!}
  \sum_{|\alpha|=k}\lambda^\alpha~,
\end{equation}
which shows that the $k$-th moment of the numerical measure $\mu_A$ is
the complete symmetric homogeneous polynomial of degree $k$ in the
variables $\lambda_1,\dots,\lambda_n$, divided by the combinatorial
factor $\binom{n+k-1}{k}$ which is just the number of terms in the
sum.  

Using the Newton identities, the right-hand side of
\eqref{eq:momk} can be decomposed into as a sum of products of 
elementary symmetric polynomials of the form $p_\ell = \lambda_1^\ell 
+ \dots + \lambda_n^\ell$, see \cite[Eq.~(2.14')]{McD}. If we assume
that $\Tr(A) = 0$, then $p_1 = 0$ and the number of nonzero terms 
in the sum is considerably reduced. Using these remarks, it is not 
difficult to show that, under the assumptions of 
Proposition~\ref{prop:sconv}, the $k$-th moment of the rescaled numerical 
measure $\mu_{\sqrt{n}A_n}$ satisfies
\[
  \int_\R x^k \dd\mu_{\sqrt{n} A_n}(x) 
  ~\xrightarrow[n\to\infty]{}~ \begin{cases}
  0 & \hbox{if } k \hbox{ is odd}~, \\
  2^{-k/2}\sigma^k\frac{k!}{(k/2)!} & \hbox{if } k \hbox{ is even}~. 
 \end{cases}
\]
Since the moments in the right-hand side are those of the normal law, 
we conclude that $\mu_{\sqrt{n}A_n}$ converges weakly to $\NN(0,\sigma^2)$ 
as $n \to \infty$.

\bigskip
We now consider general matrices $A_n\in{\bf M}_n(\C)$, and 
obtain a central limit theorem by applying Proposition~\ref{p:CLH} 
to the Radon transform of $A_n$. 

\begin{thm}\label{th:limgen}
Let $\{A_n\}_{n\ge1}$ be a sequence of matrices satisfying $A_n \in 
\bM_n(\C)$, $\Tr(A_n) = 0$, and $(\log n) \|A_n\|/n^{1/2} \to 0$ as 
$n \to \infty$. 
We assume that 
\begin{equation}\label{eq:clmhyp}
  \frac{1}{n}\,\Tr(A_n^* A_n) \,\xrightarrow[n\to\infty]{}\, a > 0~,
  \qquad \frac{1}{n}\,\Tr(A_n^2) \,\xrightarrow[n\to\infty]{}\, b\in[0,a)~.
\end{equation}
Then the rescaled numerical measure $\mu_{\sqrt{n}A_n}$ converges 
weakly to the Gaussian measure $f_\infty (z)\dd z$ as $n \to 
\infty$, where
\begin{equation}\label{def:fin}
  f_\infty(x +iy) \,=\, \frac{1}{\pi\sqrt{a^2-b^2}}
  \,e^{-\frac{x^2}{a+b}-\frac{y^2}{a-b}}~.
\end{equation}
\end{thm}

\bepr
For each $\theta \in S^1$, the Hermitian matrices $H_n(\theta)
= \frac12 (e^{-i\theta}A_n + e^{i\theta}A_n^*)$ satisfy the assumptions 
of Proposition~\ref{p:CLH}, with
\[
  \sigma^2 \,=\, \sigma(\theta)^2 \,=\, \lim_{n\to \infty}\frac{1}{2n}
  \Bigl(\Tr(A_n^*A_n) + \Re(\Tr(A_n^2)e^{-2i\theta})\Bigr)
  \,=\, \frac12(a+b\cos(2\theta))~.
\]
Let us denote by $\mu_n$ and $\mu_{n,\theta}$ the numerical measures of 
$\sqrt{n}\,A_n$ and $\sqrt{n}\,H_n(\theta)$, respectively. Since 
$\mu_{n,\theta}$ is the Radon transform of $\mu_n$, the two-dimensional
Fourier transform $\hat\mu_n(\xi)$ for $\xi = r e^{i\theta}$ is precisely
the one-dimensional Fourier transform of $\mu_{n,\theta}$ evaluated
at $r$ \cite{Hel}. Thus, applying Proposition~\ref{p:CLH}, we find
\[
  \hat\mu_n(r e^{i\theta}) \,=\, \int_\R e^{-ixr}\dd\mu_{n,\theta}(x)
  \,\xrightarrow[n\to\infty]{}\, \frac{1}{\sigma(\theta)\sqrt{2\pi}}
  \int_\R e^{-ixr} e^{-x^2/(2\sigma(\theta)^2)}\dd x \,=\, 
  e^{-\sigma(\theta)^2r^2/2}~,
\]
for any $r \ge 0$, $\theta \in S^1$. This shows that $\mu_n$ converges 
weakly as $n \to \infty$ to the measure $\mu_\infty$ on $\C$ defined by
\[
  \hat\mu_\infty(\xi) \,=\, e^{-\frac{1}{2}|\xi|^2\sigma(\theta)^2}
  \,=\, e^{-\frac{a+b}4\Re(\xi)^2}e^{-\frac{a-b}4\Im(\xi)^2}~, \quad 
  \xi \in \C~.
\]
Inverting the Fourier transform, this gives ${\rm d}\mu_\infty = 
f_\infty(z)\dd z$ with $f_\infty$ as in \eqref{def:fin}. 
\enpr

\bigskip\noindent{\bf Remarks.}\\
{\bf 1.} If we assume for simplicity that $\|A_n\|$ is uniformly 
bounded, we can suppose (up to extracting a subsequence) that
$\frac1n\Tr(A_n^*A_n)$ converges as $n \to \infty$ to some 
$a \ge 0$. However, we have to assume in \eqref{eq:clmhyp} that 
$a > 0$ in order to get a universal Gaussian limit.\\[1mm]
{\bf 2.} Similarly, the first assumption in \eqref{eq:clmhyp}
implies that $\frac1n\Tr(A_n^2)$ converges, after extracting
a subsequence, to some $b\in\C$ with $|b| \le a$. Multiplying $A_n$ 
by a unit complex number, we can assume that $0 \le b \le a$, but in 
the borderline case where $b = a$ the limiting measure $f_\infty(z)\dd z$ 
should be replaced $(2\pi a)^{-1/2} e^{-x^2/(2a)}\dd x \otimes\delta_{y=0}$.

\bigskip\noindent{\bf Example.} Let us consider the Jordan block 
of size $n$:
\[
  A_n \,=\, \begin{pmatrix} 
  0 & 1 & 0 & \cdots & 0 \\ 
  \vdots & \ddots & \ddots & \ddots & \vdots \\ 
  \vdots &  & \ddots & \ddots & 0 \\  
  \vdots &  &  & \ddots & 1 \\  
  0 & \cdots & \cdots & \cdots & 0 \end{pmatrix}~.
\]
Since $\Tr(A_n^*A_n) = n-1$ and $\Tr(A_n^2)=0$, the assumptions of
Theorem~\ref{th:limgen} are satisfied with $a = 1$ and $b = 0$.
Thus the numerical measure of $\sqrt{n}A_n$ converges weakly to the 
normal distribution $\pi^{-1}e^{-|z|^2}\dd z$ as $n \to \infty$. 
In this example, the measure $\mu_{A_n}$ has in fact a radially 
symmetric density for all $n \ge 1$, see Section~\ref{ss:rad}.

\section{Perspectives}\label{s:persp}

As a conclusion, we briefly mention  a natural extension
of our work, which is left for future investigation. Recall that
a homogeneous polynomial $P \in \R[X_0,X_1,\dots,X_d]$ of
total degree $n$ is {\em hyperbolic} in the direction $e_0 =
(1,0,\dots,0)$ if, on the one hand, it has partial degree $n$
with respect to the first variable $X_0$, and on the other hand
the $n$ roots of the univariate polynomial $t\mapsto P(t,y)$
are real for every vector $y\in\R^d$. Hyperbolic polynomials
arise as principal symbols of hyperbolic differential operators of
order $n$ in $d$ space variables, see \cite{Gar}. As an example,
if $A_1,\dots,A_d \in \bM_n(\C)$ are Hermitian matrices, the
polynomial
\begin{equation}\label{eq:exhyp}
   P(X_0,X_1,\dots,X_d) \,=\, \det(X_0I_n-X_1A_1-\dots-X_dA_d)
\end{equation}
is hyperbolic. In the particular case where $d = 2$, it has been
conjectured in \cite{Lax}, and proved in \cite{HV}, that all monic
hyperbolic polynomials are of the form \eqref{eq:exhyp}.
This is no longer true if $d \ge 3$.

It might be argued that a large part of our work is not really
about matrices, but rather concerns hyperbolic polynomials in
$2+1$ variables. Indeed, if $A \in \bM_n(\C)$ and $A_1$, $A_2$
are as in \eqref{def:Htheta}, the eigenvalues $\lambda_1(\theta),
\dots,\lambda_n(\theta)$ of the Hermitian matrix $H(\theta)$ are
the solutions of the equation
\[
   P_A(\lambda,\cos\theta,\sin\theta) \,=\, 0~, \qquad \theta
   \in S^1~,
\]
where $P_A(X_0,X_1,X_2) = \det(X_0 I_n - X_1 A_1 - X_2 A_2)$ is
the hyperbolic polynomial associated with $A_1, A_2$. As was
shown in Section~\ref{s:Radon}, the numerical measure $\mu_A$
is entirely determined by the eigenvalues $\lambda_j(\theta)$,
hence by the polynomial $P_A$.

This in turn suggests a natural way to associate to any hyperbolic
polynomial $P$ of degree $n$ in $d+1$ variables a probability
measure $\mu_P$ on $\R^d$, which coincides with the numerical
measure $\mu_A$ when $d = 2$ and $P = P_A$. Given a unit vector
$\omega \in S^{d-1} \subset \R^d$, let $\lambda_1(\omega),\dots,
\lambda_n(\omega) \in \R$ be the roots of the polynomial equation
$P(\lambda,\omega) = 0$, and let $B_\omega(s) = B[\lambda_1(\omega),
\dots,\lambda_n(\omega)](s)$ be the normalized $B$-spline with
knots $\lambda_1(\omega),\dots,\lambda_n(\omega)$. The
``numerical measure'' of $P$ is then the unique probability
measure $\mu_P$ on $\R^d$ whose Radon transform satisfies
\[
   (\RR \mu_P)(\omega,\dd s)  \,=\, B_\omega(s)\dd s~,
\]
where, by definition, $(\RR \mu_P)(\omega,I) = \mu_P(\{x\in \R^d \,|\, 
x\cdot\omega \in I\})$ for any interval $I\subset \R$.

In this generalized setting, the counterpart of normal matrices is 
the case where the polynomial $P$ split into linear factors
\[
  P(X) \,=\, \prod_{k=1}^n(X_0-v^k\cdot(X_1,\ldots,X_d))~, \qquad
  \hbox{with}~v^k \in \R^d~.
\]
When $n\ge d+1$ and the vectors $v^k$ span the affine space $\R^d$, the 
density of $\mu_P$ with respect to the Lebesgue measure is the 
multivariate $B$-spline in $d$ variables, whose nodes are the $v^k$'s. 
It is piecewise polynomial of degree $n-d-1$. In the generic situation 
where any $(d+1)$-uplet of vectors $v^k$ is an affine basis, it is of 
class $C^{n-d-2}$, see \cite{Dah}. Again the density is $\log$-concave 
in this case.

As in the two-dimensional case, the measure $\mu_P$ can be expressed
in terms of $B_\omega$ using the back-projection method. The example 
above suggests that, as the space dimension increases, the measure $\mu_P$
becomes more singular. Then, the inversion formula has to be understood
in the sense of distributions. In the three-dimensional case $d = 3$,
we arrive at the simple expression
\begin{equation}\label{eq:invTd}
   \mu_P = -\frac{1}{8\pi^2}\,\Delta_x \int_{S^2}B_\omega(x\cdot\omega)
   \dd \sigma(\omega)~,
\end{equation}
where $\dd\sigma$ denotes the Euclidean measure on the unit sphere
$S^2$.  As an example, in the particular situation where $P =
X_0^2-X_1^2-X_2^2-X_3^2$, which corresponds to the differential
operator $\partial_t^2-\Delta$ of the wave equation, we obtain $\mu_P
= \frac{1}{4\pi}\dd\sigma(x)$.  This shows that, when $d \ge 3$, the
support of $\mu_P$ does not need to be convex. Because the wave
equation satisfies the Huyghens Principle, this example suggests that
the link between polynomial regions of the density and lacunas of
differential operators persists in higher dimensions.

\end{document}